\documentclass[leqno,oneside]{amsart}
\usepackage{latexsym,amssymb,amsmath,graphics}

\def\gg{{\mathcal{G}}}
\def\kk{{\mathcal{K}}}

\def\B{{\mathbb{B}}}
\def\C{{\mathbb{C}}}
\def\E{{\mathbb{E}}}

\def\N{{\mathbb{N}}}
\def\P{{\mathbb{P}}}

\def\R{{\mathbb{R}}}

\def\Z{{\mathbb{Z}}}

\newcommand{\8}{\infty}
\renewcommand{\d}{\delta}
\renewcommand{\a}{\alpha}
\renewcommand{\b}{\beta}
\newcommand{\D}{\Delta}

\renewcommand{\l}{\lambda}
\renewcommand{\L}{\Lambda}
\newcommand{\p}[1]{\P\big[#1\big]}
\newcommand{\e}[1]{\E\big[#1\big]}
\newcommand{\eps}{\varepsilon}

\newcommand{\g}{\gamma}
\newcommand{\s}{\sigma}

\newcommand{\ov}{\overline}
\newcommand{\wt}{\widetilde}

\renewcommand{\wt}{\widetilde}
\newcommand{\wh}{\widehat}
\newcommand{\ro}{\rho}
\renewcommand{\t}{\theta}
\renewcommand{\th}{\theta}
\renewcommand{\d}{\delta}
\renewcommand{\e}{\varepsilon}
\newcommand{\tel}{\th,\e,\l}
\newcommand{\is}[2]{\langle #1,#2\rangle}
\newcommand{\iss}[1]{\langle #1\rangle}
\newcommand{\bis}[2]{\Big\langle #1,#2\Big\rangle}
\newcommand{\bbis}[2]{\bigg\langle #1,#2\bigg\rangle}
\newcommand{\suf}[1]{\lceil #1\rceil}

\newcommand{\psitv}{\psi_{t,v}}

\newcommand{\xm}{{x_-}}
\newcommand{\xp}{{x_+}}
\newcommand{\ym}{{y_-}}
\newcommand{\yp}{{y_+}}
\newcommand{\vm}{{V_{\frac \a2,-}}}
\newcommand{\vp}{{V_{\frac \a2,+}}}

\newcommand{\dt}{b_n}

\newtheorem{thm}[equation]{Theorem}
\newtheorem{cor}[equation]{Corollary}
\newtheorem{lem}[equation]{Lemma}

\newtheorem{mthm}[equation]{Main Theorem}

\newtheorem{prop}[equation]{Proposition}
\theoremstyle{definition}

\numberwithin{equation}{section}

\begin{document}

\title[ Convergence to stable laws]{Convergence to stable laws
for a class of multidimensional stochastic recursions}
    \dedicatory{\normalsize In memory of Andrzej Hulanicki}
\author[D. Buraczewski, E. Damek, Y. Guivarc'h]
{Dariusz Buraczewski, Ewa Damek, Yves Guivarc'h}
\address{D. Buraczewski, E. Damek \\
University of Wroclaw\\
Institute of Mathematics\\
pl. Grunwaldzki 2/4\\
 50-384 Wroclaw\\
 Poland}
\email{dbura@math.uni.wroc.pl\\ edamek@math.uni.wroc.pl}
\address{Y. Guivarc'h\\ IRMAR, Universit\'e de Rennes 1\\
Campus de Beaulieu\\
35042 Rennes cedex, France} \email{yves.guivarch@univ-rennes1.fr}

\thanks{
This research project has been partially supported by
 Marie Curie Transfer of Knowledge Fellowship {\it Harmonic Analysis,
Nonlinear Analysis and Probability} (contract number MTKD-CT-2004-013389).
D.~Buraczewski and E.~Damek were also supported by KBN
grant N201 012 31/1020.}

\begin{abstract}
We consider a Markov chain $\{X_n\}_{n=0}^\8$ on $\R^d$ defined by
the stochastic recursion $X_{n}=M_n X_{n-1}+Q_n$, where $(Q_n,M_n)$
are i.i.d. random variables taking values in the affine group $H=\R^d\rtimes
{\rm GL}(\R^d)$. Assume that $M_n$ takes values in the similarity
group of $\R^d$, and the Markov chain has a unique stationary
measure $\nu$, which has unbounded support.
We denote by $|M_n|$ the expansion coefficient of $M_n$ and we assume
$\E |M|^\a=1$ for some positive $\a$.
We show that the partial
sums $S_n=\sum_{k=0}^n X_k$, properly normalized, converge to a
normal law  ($\a\ge 2$) or to an infinitely divisible law, which is stable in a
natural sense ($\a<2$). These laws are fully nondegenerate, if $\nu$ is not
supported on an affine
 hyperplane.  Under a natural hypothesis,
 we prove also a local limit theorem for the sums $S_n$.
If $\a\le 2$, proofs are based on the homogeneity at infinity
of $\nu$ and
on a detailed spectral analysis of
a family of
Fourier operators $P_v$ considered as   perturbations of the transition
operator $P$ of the chain $\{ X_n \}$. The characteristic function of the limit
law has a  simple expression in terms of    moments of $\nu$ ($\a > 2$)
or of the tails of  $\nu$ and of  stationary measure for an associated Markov operator ($\a\le 2$).
We extend the results to the situation where $M_n$ is a random generalized
similarity.
\end{abstract}

\maketitle
\section{Introduction and main results}
We consider the vector space $V=\R^d$ endowed with the scalar product
$\is xy=\sum_{i=1}^d x_i y_i$ and the norm $|x|=\sum_{i=1}^{d}|x_i|^2$.
Let $H=V\rtimes {\rm GL} (V)$ be the affine group of $V$
i.e. $H$ is a semi-direct product of the linear group ${\rm GL} (V )$
and the group of translations of $V$. The action of $h=(b,g)$,
$b\in V$, $g\in {\rm GL} (V )$ on $x\in V$ is
$$
hx=gx +b. $$
We denote by $u^*$ the adjoint operator of $u\in {\rm End } V$.

\medskip

 Given a probability measure $\mu $ on $H$ and
$x\in V $, we consider the recurrence relation with random coefficients
\begin{equation}
\label{rekursja}
\begin{split}
X_0^x&=x,\\
X_n^x&=M_nX_{n-1}^x +Q_n,
\end{split}
\end{equation}
where the random pairs $(Q_n, M_n)\in H$ are independent and distributed
according to  $\mu$. We assume that a
unique stationary law $\nu$ for this recursion exists and has unbounded support.
We denote by $\ov\mu$ the projection of $\mu$ on ${\rm GL}(V)$ and by $G_{\ov\mu}$
the closed subgroup generated by supp$\ov \mu$.

\medskip

We are interested in
the limiting behavior of the sum $S_n^x = \sum_{k=0}^n X_k^x$ of the
non independent random variables $X_k^x$ ($0\le k\le n$).
Such
a problem  was
considered in \cite{K}, and convergence to stable laws
 for sums like $S_n^x$,
but with i.i.d. increments, was stated there.
Under some conditions,  the homogeneity at infinity of stationary laws was proved  and was an
essential aspect of the limit theorems.
For general information on limit theorems for analogous situations see \cite{AD,BL,BDP,HH}.
For a study of homogeneity   of tails
in    closely related contexts see \cite{G,GL2}.
For motivations to consider such affine recursions see e.g. \cite{DF}.

For the main part of the paper we will assume that $M_n$ belongs to the similarity group $G$
of $V$, i.e. the group of elements $g$ of ${\rm GL }(V)$ satisfying
$$
|g v| = |g||v|,
$$ for every $v\in V$. In this case under some moment conditions, including
$\E|M_n|^\a=1$ for some $\a>0$,
 a detailed
study of the (unique) finite stationary measure for \eqref{rekursja}
and of its tail $\L$ is available (see \cite{BD}).
  For  a study of tails in  closely related  1-dimensional models see \cite{Gr1,Gol,GL2}.
  We observe that $\L$ is homogeneous of degree $\a$ with respect to $G_{\ov\mu}$.
  In contrast to the general case of recursion \eqref{rekursja}, we observe that here, modulo a compact
  subgroup, $G_{\ov\mu}$ is isomorphic to $\R$ or $\Z$. If $\a\le 2$, this
  fact will be reflected in the form of the limit laws.
   If $G_{\ov\mu}$  contains $\R^*_+$ or if $\a > 2$, then $\nu$ belongs to the domain of attraction of a stable law. More generally,
  if $\a<2$
  the concept of semistability and normalization along a subsequence of integers is relevant (see \cite{L}).
   We show that the limiting law of the properly normalized sum $S_n^x$
exists, is infinitely divisible and stable in a natural sense.
If ${\rm supp} \mu$ has no invariant affine subspace, this law is fully
nondegenerate. If $\a\le 2$, the tail $\L$ enters as an essential component
in the description of the limit law.
 If $\a\ge 2$, this law is normal and if $\a>2$ its covariance form is a
simple modification of the covariance form of $\nu$.
In particular, if $\mu$ varies continuously and satisfies very general moment
conditions, one passes from  Gaussian asymptotics $(\a>2)$
to non Gaussian ones $(\a<2)$. This is analogous
to a  phase transition, as in statistical physics (see \cite{S,DLNP}).

If $\a\le 2$, in particular in the non normal case,
the description of the  parameters of the limit law for $S_n^x$ involves
another family of Markov chains  and stationary measures.
For any fixed nonzero $v\in V$, the Markov chain on $V$ defined
by the recursion
\begin{equation}
\label{rekursja2}
\begin{split}
W_0&=0,\\
W_{n}& = M_{n}^*(W_{n-1}+v),
\end{split}
\end{equation}
has also a finite stationary measure $\eta_v$. It turns out that the  tails of family $\eta_v$
enter in the expression of the limit law for $S_n^x$.
We observe, that in most cases of convergence to non normal stable laws
for functionals of Markov chains, which are considered in the literature,
the Birkhoff sums have the same limiting behavior as if the increments
were i.i.d. with law equal to the stationary measure of the chain
(see for example \cite{GLJ}, for the case of continuous
fraction expansion).
  This is not the case here and the limit law has a tail, which depends
  linearly of the tails of the stationary laws $\nu$ and $\eta_v$.

\medskip

In order to state our main results, we need some notations.
For $v\in V$ we write $\chi_v(x) = e^{ i \is vx}$,
$v^*(x) = \is vx$ and the characteristic function of a probability
measure $\th$ on $V$ will be written $\wh \th(x) = \int_V\chi_x(y) \th(dy)$.

We will say that $\mu$ or recursion \eqref{rekursja} satisfies hypothesis {\bf H} if
\begin{itemize}
\item No point of $V$ is invariant under the action of ${\rm supp} \mu$.
\item There exists $\a>0$ with $\E |M|^\a=1$.
\item $m_\a = \E[|M|^\a\log|M|]<\8$ and $\E|Q|^\a<\8$.
\end{itemize}
Hypothesis {\bf H} implies $\E [\log|M|] <0$, hence (see \cite{Bra})
 the Markov chain defined
by \eqref{rekursja} has a unique stationary measure $\nu$
and  the support of $\nu$ is unbounded.
The affine  subspace
generated by ${\rm supp} \nu$
is ${\rm supp} \mu$ - invariant and, if useful, we can assume
that there is no proper    ${\rm supp} \mu$ - invariant affine subspace.
The transition operator
of the chain $\{X_n^x\}$ will be denoted by $P$,
hence
$$
P\phi(x) = \E[\phi(X_1^x)] = \int_H \phi(gx+b)\mu(dh).
$$
Also the series $Q_0+\sum_{k=1}^\8 M_0\ldots M_{k-1}Q_k$ converges
$\P$-a.e. to  a $V$-valued random variable $R$ and $\nu$ is the law of $R$.
Similar properties are valid for the Markov chain \eqref{rekursja2}  associated to the transition
operator $T_v$ ($v\in V$) given by
$$
T_v\phi(x) = \int_G\phi(g^*(x+v))\ov\mu(dg)
$$
and we denote by $\eta_v$ its unique stationary measure, i.e. the law of $\big(\sum_1^\8 M_0^*\cdots M_{k-1}^*\big)v$.

\medskip
The group $G$ is the direct product of $\R^*_+$ and the
 orthogonal group $K=O(V)$. We denote by $R_{\ov\mu}$ the projection of $G_{\ov\mu}$
 on $\R_+^*$.
  The center of $G_{\ov\mu}$ will be denoted by $Z_{\ov\mu}$.
Let $K_{\ov\mu}=G_{\ov\mu}\cap K$. Since $\P[|M_n|=1]<1$,  and $R_{\ov\mu}$ is closed
we have
\begin{equation*}
R_{\ov \mu}=\R _+^*\ \ \mbox{or}\ \ R_{\ov \mu}=\langle p\rangle
=\{p^n:n\in \Z \}\ \mbox{for a}\ p>1.
\end{equation*}
There exists a closed subgroup $A_{\ov\mu}\subset G_{\ov\mu}$
such that the projection $g\mapsto |g|$ defines an isomorphism
of $A_{\ov\mu}$ onto $R_{\ov\mu}$, and $G_{\ov\mu}=A_{\ov\mu}\ltimes K_{\ov\mu}$
is the semidirect product of $A_{\ov\mu}$ and $K_{\ov\mu}$. Furthermore,
$A_{\ov\mu}$ can be chosen to contain a central subgroup of $G_{\ov\mu}$
as a finite index subgroup.
In particular the center $Z_{\ov\mu}$ of $G_{\ov\mu}$ is the product of $Z_{\ov\mu}\cap K$
by a subgroup isomorphic to $\R$ or $\Z$.
Below, the elements of $Z_{\ov\mu}$ (if $0<\a\le 2$) or $\R^*_+$ (if $\a>2$) will be used
to normalize the sums $S_n^x$. If $G_{\ov\mu}\supset \R^*_+$ or if $\a>2$, the normalization is as
usual, by positive numbers.
See  Appendix \ref{app1} for some further discussion
on the structure of $G_{\ov\mu}$ and $Z_{\ov\mu}$.

We denote by $\Sigma_1$ the fundamental domain of $A_{\ov\mu}$ on $V\setminus\{0\}$ given by:
$\Sigma_1 = \{ x\in V;\; 1\le |x|<p\}$  if  $R_{\ov\mu}=\langle p\rangle$,
$\Sigma_1 = S_1$,  the unit sphere of $V$, if  $R_{\ov\mu}=\R_+^*$.
Then we write $x=a(x)\ov x$ with $a(x)\in A_{\ov\mu}$ and $\ov x\in \Sigma_1$.
Then $r(x)=|a(x)|\le|x|$ takes values in $R_{\ov\mu}$, and $r(x)=|x|$ if $R_{\ov\mu} = \R_+^*$.

\medskip

It is shown in \cite{BD} that under hypothesis {\bf H}, the following
$G_{\ov\mu}$-homogeneous Radon measure $\L$ is well defined by the following
weak convergence on $V\setminus\{0\}$
\begin{equation}
\label{tailm}
\L = \lim_{|g|\to 0, g\in G_{\ov\mu}}|g|^{-\a} g\nu.
\end{equation}
Then $\L$ is called the tail measure (or tail)  of $\nu$ and the support of $\L$ is studied
there under natural conditions. Here we need the fact that $\Lambda$
is nonzero and this is a consequence of hypothesis ${\bf H}$ only.
In the case $d=1$ and $R_{\ov\mu}=\R^*_+$ the measure $\L$ is defined by
$$
\L(dx) = C_+{\bf 1}_{(0,\8)}(x)\frac{dx}{x^{\a+1}} +C_-{\bf 1}_{(-\8,0)}(x)\frac{dx}{|x|^{\a+1}}.
$$
In general $\L$ has a product form. Let $l$ be the Haar measure on $A_{\ov\mu}$ i.e. either $l(da)=\frac{d|a|}{|a|}$
if $R_{\ov\mu}=\R^*_+$ or $l$ is the counting measure multiplied by $\log p$, if $R_{\ov\mu}=\langle p \rangle$.
Define $l^\a(da)=|a|^{-\a}l(da)$, then there exists a finite measure $\s$ on the fundamental domain $\Sigma_1$
such that $\L$ can be written as the product of $l^\a$ and $\s$:
\begin{equation}
\label{radial}
\int_{V\setminus\{0\}} f(x)\L(dx) = \int_{A_{\ov\mu}}\int_{\Sigma_1}f(aw)l^{\a}(da)\s(dw).
\end{equation}
 Also, if $\a>1$, we will denote by $m=\int_V x\nu(dx)$
the mean of $\nu$, by $q(x,y)=\int_V \is{x}{\zeta-m}\is{y}{\zeta-m}\nu(d\zeta)$,
the covariance form  of $\nu$.
We write also $z=\E[M_n]$ for the averaged operator of $M_n$.

\medskip

It will be shown that $\eta_v$ has also a tail $\D_v$ given by
$$
\D_v =\lim_{|g|\to 0,g\in G_{\ov\mu}} |g|^{-\a} g^*\eta_v.
$$
We define the $G_{\ov\mu}^*$-homogeneous continuous function $\wt \L$ on $V$ by
\begin{eqnarray*}
\wt\L(y) &=&
\int_V (\chi_y(x)-1)\L(dx),\quad\mbox{if }0<\a<1,  \\
\wt\L(y) &=& \int_V \bigg(\chi_y(x)-1 -\frac{i\is xy}{1+|y|^2|x|^2}\bigg)\L(dx),\quad\mbox{if }\a=1,\\
\wt\L(y) &=& \int_V \big(\chi_y(x)-1-i\is xy \big)\L(dx),\quad\mbox{if }1<\a<2,\\
\wt\L(y) &=& -\frac 14\int_{\Sigma_1}\is yx^2\s(dx),\quad\mbox{if } \a=2.
\end{eqnarray*}
We denote by $\wt\L^1$ the function on $V$ defined by
$$
\wt\L^1(y) = \wt \L (\ov y) {\bf 1}_{[1,\8)}(r(y)).
$$

Given a closed subgroup $U$ of ${\rm GL}(V)$ and a continuous homomorphism
$\ov \a$ of $U$ in $\R^*_+$, we will say that a probability measure
$\th$ on $V$ is $(U,\ov \a)$ stable if $\th$ belongs to a one parameter
convolution semigroup $\th^t$ ($t\ge 0$) and for
every $u\in U$, there exists $\b(u)\in V $ such that
$$
u(\th) = \th^{\ov\a(u)}*\b(u).
$$
This equation implies that $\wh\th(\l)$ ($\l\in V$) do not vanish and if
$\phi(\l) = \log \wh\th (\l)$, then for any $u\in U$, $\phi(u^* \l)=\ov \a(u)
\phi(\l) + i\is{\b(u)}{\l}$. Conversely these conditions imply the $(U,\ov\a)$
stability of $\th$, and in particular $\th$ belongs to a well defined one parameter convolution semigroup.

If $U\subset G$, the structure of $U$ implies that $\ov\a$ is of the form  $\ov\a(u)=|u|^\a$
with $\a>0$, and if $\a\not=1$ the stability relation can be reduced, using
translations, to $u(\th)=\th^{\ov \a(u)}$.
In the case $d=1$, $U=\R_+^*$, $U$-stability coincides with stability in the classical
sense. If $d\ge 1$, $U\supset \langle p \rangle $ ($p>1$) (resp. $U\supset \R^*_+$)  and $\a\not=1$, $U$-stability
coincides with semi-stability (resp. stability) in the sense of \cite{L}.

\begin{mthm}
\label{mthm}
Assume that the probability measure $\mu$ on $H$ satisfies hypothesis
{\bf H}. Then for any $x\in V$.
\begin{enumerate}
\item If $\a>2$, $\frac 1{\sqrt n} (S_n^x-nm)$ converges in law
to the normal law with the Fourier transform
$$
\Phi_{2+}(v)=\exp\big( - q(v,v)/2 - q(v,(I-z^*)^{-1}z^* v)\big).
$$
\item If $\a\in(0,2)$, assume $c_n\in {Z_{\ov\mu}}$ is related to $n\in\N$ by $\big[|c_n|^{-\a}\big]=n$
and define $d_n=0,=n\xi(c_n),=nc_nm$, resp. if $ \a<1,=1,>1$,
where $\xi(c)= \int_V\frac{cx}{1+|cx|^2}\nu(dx)$ for $c\in Z_{\ov\mu}$.
Then $c_nS_n^x -d_n$ converges
in law towards the $(Z_{\ov\mu},\ov\a)$-stable law with
the Fourier transform $\Phi_\a(v)=\exp C_\a(v)$, with
\begin{eqnarray*}
C_\a(v) &=&\a m_\a \D_v(\wt\L^1), \quad \mbox{if } \a\not=1,\\
C_\a(v) &=& m_1 \D_v(\wt\L^1) + i\g(v), \mbox{if } \a=1,
\end{eqnarray*}
where $\g(v)\in\R$, if $R_{\ov\mu}=\R_+^*$. If $R_{\ov\mu}=\langle p \rangle$, the same formulas are valid, where $\a \D_{v}(\wt\L^1)$
is replaced by $\frac{1-p^{-\a}}{\log p}\D_v(\wt \L^1)$.
Furthermore if $\a=1$, then for some constant $I(v)>0$:
\begin{eqnarray*}
|\xi(c)|&\le& I(v)|c| |\log|c||,\quad \mbox{ if } |c|< 1/2,\\
|\xi(c)|&\le& I(v)|c|,\quad \quad \qquad \mbox{ if } |c|\ge 1/2.
\end{eqnarray*}
\item
If $\a=2$, assume $c_n\in Z_{\ov\mu}$ satisfy $\lim_{n\to\8} |c_n|\sqrt{n\log n}=1$, then
$c_n(S_n^x-nm)$ converges in law to the normal law with the Fourier transform
$ \Phi_2(v)= \exp(C_2(v))$, with
$$ C_2(v)=  -\frac 14 \int_{\Sigma_1}\Big( \is vw^2 + 2\is vw\eta_v(w^*)\Big)\s(dw)= 2 \D_v(\wt\L^1), $$
if $R_{\ov\mu}=\R_+^*$. If $R_{\ov\mu}=\langle p \rangle$, the same formula is valid
with $\frac {1-p^{-2}}{\log p}$ instead of $2$. In both cases $C_2(c^*v)=|c|^2C_2(v)$ if $c\in Z_{\ov\mu}$
\end{enumerate}
If  no affine subspace of $V$
is ${\rm supp} \mu$ invariant, then the limit laws are fully
nondegenerate i.e., their supports are not contained in a proper
subspace of $V$.
\end{mthm}
{\bf Remarks}

a) If
 $G_{\ov\mu}\supset   \R^*_+$ or if $\a>2$, $\Phi_\a(v)$ is the characteristic
function of a multidimensional  stable law in the  sense of  \cite{L} (p. 213-224).


b) In case $d=1$ and $G_{\ov\mu}=\R^*_+$, the analogue of
Theorem \ref{mthm} has been proved in \cite{GL}.
For another proof of assertion (1) in Theorem \ref{mthm} in a more general
context and
under a moment condition of order 4, see \cite{HH1}.

c) If $R_{\ov\mu}=\langle p\rangle, \a<2$ the sequence $c_n$ given in Theorem \ref{mthm}
is lacunary, hence also the sequence of integers defined there. However
the limit  law is infinitely divisible; in general the tail of
$\nu$ has a nontrivial periodic  multiplicative part, hence $\nu$ do not belong to the domain
of attraction of a stable law (see \cite{F}, p. 577), then the limit law is only semistable in
the sense of \cite{L}. If $\a=2$ and $R_{\ov\mu}=\langle p\rangle$, the sequence $c_n$ is also lacunary but the limit law is normal.

d) If $\a=2$, since $C_2(c^*v)=|c|^2 C_2(v)$ if $c\in Z_{\ov\mu}$ and $C_2(v)$ is a quadratic
form, the corresponding normal law is invariant under the subgroup of $K$, which is the projection of $Z_{\ov\mu}$
on $K$.

e) As in \cite{GL} the proofs follow the Fourier analytic approach of \cite{GH}
(see also \cite{BDP,HH}).
However, here the dominant eigenvalue of the Fourier operator is not analytic
and even not differentiable if $\a<2$.
Thus, an important point is to get explicit asymptotic fractional expansions.
This is based on the homogeneity
at infinity of stationary measures, studied in \cite{BD} and a remarkable
intertwining relation.
Moreover, instead of the analytic perturbation theory used in \cite{GH}, we
need to use here the operator perturbation theorem of \cite{KL}.

\begin{mthm}[Local Limit Theorem]
\label{mthmloc}
Assume that $R_{\ov\mu} = \R^*_+$, hypothesis {\bf H} is satisfied,
no affine subspace of $V$ is ${\rm supp}\mu$-invariant and $\a \notin\{1,2\}$.
Then for every $v\in V$ and domain $I\subset \R^d$ with negligible boundary
$$
\lim_{n\to\8} {n^{\chi}} \P\big[  S_n^x-d_n\in I \big] =
p_\a(0){\lambda(I)},
$$
where
\begin{itemize}
\item $\chi=\frac d{\a }, \frac d{2}$ if $\a<2,>2$, resp.
\item $d_n=0,=nm$ if $\a<1,>1$, resp.
\item
$p_\a$ is the density of the corresponding  limit law in Theorem \ref{mthm};
\item  $\lambda(I)$ denotes the Lebesgue measure of $I$.
\end{itemize}
\end{mthm}

{\bf Remark.} This theorem can be interpreted as a local limit theorem
for a random walk defined by $\mu$ on a homogeneous space $\widetilde V$
of a larger group $\wt H$ (see Section \ref{section-local}). Then we see that the exponent $\chi$
of the corresponding local limit asymptotics  is determined by the geometry
of $(\wt H,\wt V)$ if $\a>2$, while it depends strongly of $\mu$ if $\a<2$. Such a situation,
in  case of Lie groups,  was considered in \cite{V}.

\medskip

In order to get an idea of what happens in general case we
 consider also the  more general situation of generalized
similarities. We will say that $g\in {\rm GL} (V)$ is a componentwise  similarity
if $V$ is an orthogonal direct sum $V=\oplus_{j=1}^lV_j$ and $g$ acts on $V_j$
through a similarity $g_j$, i.e. for any $x_j\in V_j$, $gx_j\in V_j$ and $|gx_j|=|g_j||x_j|$.
We write $x=\sum_{j=1}^l x_j$, $g=(g_1,\ldots,g_l)$. Here we fix positive
numbers $1=\l_1<\l_2\cdots<\l_l$ and an orthogonal direct sum $V=\oplus_{j=1}^l V_{\l_j}$.
We consider a 'homogeneous norm' $\tau$, i.e. $\tau(x)=\sum_{j=1}^l |x_{j}|^{\frac 1{\l_j}}$
and we observe that if $a>0$ and $\g_a\in {\rm GL}(V)$ is given by $\g_a(x_j)=a^{\l_j}x_j$,
then $\g_a$ is a componentwise similarity, which satisfies $\tau(\g_a x)= a\tau(x)$.
We denote $D=\{ \g_a; a\in \R^*_+ \}$, $|g|=\sup_{\tau(x)=1}\tau(gx)$, $G=\{ g\in {\rm GL}(V):\; \tau(gx)=|g|\tau(x), \forall x\in V \}$.
Then any $g\in G$ is a componentwise similarity, with $V_{\l_j}=V_j$.
If $g\in G$ we call $g$ a $\tau$-similarity.
 If $l=1$, we are back in the situation of similarities. Here we will use
 the same notations; their meaning will be clear from the context. We also denote
 $K = \{ g\in G;\; |g|=1 \}$. Then, if $K_j=K\cap {\rm GL}(V_j)$, $G_j=G\cap {\rm GL}(V_j)$,
 we have: $G_j = \R_+^*\times K_j$, $K=\prod_{j=1}^l K_j$, $G = D \times K$, where $K_j$ is
 identified with a subgroup of $G$.
 For $\g\ge 1$, we define  subspaces of $V$:
$V_{\g,-}=\oplus_{\l_j<\g}V_{\l_j}$, $V_{\g,+}=\oplus_{\l_j>\g}V_{\l_j}$. Moreover for
$\g_1 < \g_2$ we define  $V_{\g_1,\g_2} = V_{\g_1,+}\cap V_{\g_2,-}$.
For $x\in V$,  $x_{\g,+}$, $x_{\g,-}$, $x_{\g_1,\g_2}$ will denote the
projections of $x$ onto the corresponding subspaces.

Here we will assume that $M_n\in G$, hence $G_{\ov\mu}\subset G$. See Appendix for more information
on the structure of $G_{\ov\mu}$ and in particular for the fact that $G_{\ov\mu}$
has a finite index subgroup, which is the
product of $G_{\ov\mu}\cap K$ by a subgroup isomorphic to $\R$ or $\Z$. Also the center $Z_{\ov\mu}$
of $G_{\ov\mu}$ has the same form. Here $R_{\ov\mu}$ is defined as the projection of $G_{\ov\mu}$ on $D$, modulo $K$.
Moreover, the action of $G_{\ov\mu}$ on $V$ is reducible
and non isotropic. This property is reflected in the mixture of Gaussian and non Gaussian asymptotics
in the theorem below.

If $\a>1$, we define the mean of $\nu$ as above, i.e. $m_{\a,-}=\int_{V_{\a,-}}x\nu(dx)$.
Also  if $\a>2$ we define the averaged operator of $M_n$ by $z=\E\big[M_n|_{V_{\frac \a2,-}}\big]$
and the covariance form $q$ on $V_{\frac \a2,-}$ by
$q(x,y)=\int_{V_{\frac \a2,-}}\is{x}{\zeta-m}\is{y}{\zeta-m}\nu(d\zeta)$. For a description
of $\L$ in this situation see \cite{BD}, Appendix.

\begin{mthm}
\label{mthme} Assume that the probability measure $\mu$ on $H$ satisfies hypothesis {\bf H}.
Let $\{c_n\}$ be a sequence of elements  of $Z_{\ov\mu}$
such that $\big[ |c_n|^{-\a} \big]=n$ and put
$d_n= 0, = n\xi_1(c_n),=n\xi_2(c_n)$, resp. if $\a<1, \a\in [1,2), \a>2$, resp.,
where
$\xi_1(c) = cm_{\a,-}+\int_V\frac{c x_\a}{1+|cx_\a|^2}\nu(dx)$
and
$\xi_2(c) =  cm_{\frac \a2,\a} +\int_{V} \frac{c x_\a}{1+|cx_\a|^2} \nu(dx)$.
\begin{enumerate}
 \item
If $\a\in(0,2)$,
 then $c_n S_n^x - d_n$ converges in law to the $(Z_{\ov\mu},\ov \a)$ stable
law with Fourier transform
\begin{eqnarray*}
 \Phi_\a(v)&=& \exp\bigg[\int_V (\chi_v(x)-1)\wh\eta_v(x)\Lambda(dx)
\bigg], \qquad \mbox{ if } \a<1,\\
 \Phi_1(v)&=& \exp\bigg[\int_V\Big( (\chi_v(x)-1)\wh\eta_v(x)
-\frac {i\is v{x_1}}{1+|x_1|^2}
\Big)\Lambda(dx)
\bigg], \qquad \mbox{ if } \a=1,\\
 \Phi_\a(v)&=& \exp\bigg[\int_V\Big( (\chi_v(x)-1)\wh\eta_v(x)
-i\is v{x_{\a,-}} -\frac {i\is v{x_\a}}{1+|x_\a|^2}
\Big)\Lambda(dx)
\bigg], \qquad \mbox{ if } 1< \a<2.
\end{eqnarray*}
\item If $\a>2$ and $V_{\frac \a 2}=\{0\}$, then
$\frac 1{\sqrt n}(S_n^x-nm)_{\frac \a2,-}+(c_nS_n^x-d_n)_{\frac \a2,+}$
converges in law to the direct product of a normal law on $V_{\frac \a2,-}$
and a $\big(Z_{\ov\mu},\ov \a\big)$ stable law on $V_{\frac \a 2,+}$ with Fourier transforms
$$\Phi_{2+}(v)=
 \exp\Big( -q(v,v)/2 -q(v,(I-z^*)^{-1}z^*v) \Big)$$
 and
$$\Phi_{\frac \a2,+}= \exp\bigg[ \int_{V_{\frac \a 2,+}}\Big( (\chi_v(x)-1)\wh\eta_v(x)
-i\is v{x_{\a,-}} -\frac {i\is v{x_\a}}{1+|x_\a|^2}
\Big)\Lambda(dx)
\bigg].$$
\end{enumerate}
Moreover in all cases, if  no affine subspace of $V$
is ${\rm supp} \mu$ invariant, then the limiting laws are fully
nondegenerate i.e., their supports are not contained in  proper
subspaces of $V$.
\end{mthm}

{\bf Remarks}

a) If $\a\in (1,2)$ and  $V_\a=\{0\}$, the formulas for $\Phi_\a(v)$ simplify.
In this case, they extend  the formulas of the stable or semistable laws (see \cite{L}, p. 213-224).

b) If $\a>2$, then use of different normalizations  depending of the components allows to get fully
nondegenerate laws. Furthermore, the result allows to predict the value of the exponent
$\chi$ in the local limit asymptotics, as in Theorem \ref{mthmloc}: $\chi=\frac 12{\rm dim} V_{\a,-}+\frac 1\a{\rm dim } V_{\a,+}$.
The product form of the limit law is  reminiscent  of the results of \cite{GLJ} and
\cite{BaP}.

c) Here, modulo $Z_{\ov\mu}\cap K$, the normalization operators $a_n=(\frac 1{\sqrt n},c_n)$ are suitable powers of a single matrix,
which is a componentwise similarity. For a general approach to normalizations by linear
operators and limit laws of iterated convolutions see \cite{JM}.
It turns out that if $R_{\ov\mu}=D$, then the limit law in Theorem \ref{mthme} is
'operator stable' as defined in \cite{JM}, but its parameters are different
from those of  the limit law corresponding to $\nu^{*n}$.
 In the non normal case considered here
detailed information (see Appendix) on $G_{\ov\mu}$, $Z_{\ov\mu}$ is needed
 for the construction of the normalization operators.

\section{Stochastic recursions and some properties of their stationary measures}
\label{sectio-measure}
In sections \ref{sectio-measure} - \ref{section3}, we assume that $V$ is equipped
with a homogeneous norm $\tau$ and we study recursion  \eqref{rekursja}, if
$M_n$ is a $\tau$-similarity.

\medskip

Here we will describe some further properties of stationary measures
$\nu$ and $\eta_v$ of recursions \eqref{rekursja} and \eqref{rekursja2}, respectively,
that will be used in the remaining part of the paper.
If $M_n\in G$, recursion \eqref{rekursja} is studied  in \cite{BD} and proofs of all its
properties listed below can be found there. For general information on
recursion \eqref{rekursja} see \cite{Bra}.

\medskip

We define $\kappa(s)=\E|M|^s$. Under  hypothesis {\bf H}, the
function $\kappa$ is well defined for $s\in[0,\a]$ and it is strictly
convex, hence $\kappa(s)<1$ for $s<\a$.
It is known that the sequence $\{X_n^x\}_{n=0}^\8$
converges in distribution to a random variable $R$ with law $\nu$, and
finite $\t$-moments for $\th<\a$:
\begin{equation}
\label{moments}
\nu(\tau^\t)=\E\big[\tau(R)^\t\big] \le \sup_n \E\big[\tau(X_n)^\t\big] < \8.
\end{equation}
Furthermore the tail of the stationary measure $\nu$ is well
understood i.e. there exists a $G_{\ov\mu}$ -homogeneous Radon measure $\L$
on $V\setminus\{0\}$ such that
\begin{equation}
\label{has}
\lim_{|g|\to 0 , g\in G_{\ov\mu}} |g|^{-\a} g\nu(f) =
\lim_{|g|\to 0 , g\in G_{\ov\mu}} |g|^{-\a} \int_Vf (gx)\nu(dx) = \L(f)
\end{equation}
and the convergence is valid for every function $f$ such that the set of
discontinuities of $f$ has $\L$ measure 0 and for some $\eps>0$
\begin{equation}
\label{beta}
\sup_{x\not= 0 }\Big( \tau(x)^{-\a}|\log\tau(x)|^{1+\eps}|f(x)| \Big) < \8.
\end{equation}
${G_{\ov\mu}}$-homogeneity of $\Lambda$ means that for every $g\in{G_{\ov\mu}}$
\begin{equation}
\label{hom} \Lambda(f\circ g) = |g|^\a\L(f).
\end{equation}
In particular $\L$ is $K_{\ov\mu}$-invariant.

\begin{lem}
\label{2.5}
Assume $\mu$ satisfies hypothesis {\bf H}. Then $\nu$ has no atom.
If furthermore there is no proper ${\rm supp}\mu$-invariant affine subspace,
then $\nu$ gives zero measure to every affine subspace.
\end{lem}
\begin{proof}
The first assertion is a special case of Proposition 2.4
in \cite{BD}. We give a simple proof, for the sake of completeness.
Let $X$ be the set of atoms of $\nu$. Since $\sum_{x\in X}\nu(x)\le 1$, $\nu(x)$
reaches its maximum value $a$ and $X_0=\{y\in X; \nu(y)=a \}$ is finite.
On the other hand $\mu$-stationarity of $\nu$
 implies that $X_0$
is ${\rm supp }\mu$-invariant. Hence the barycenter of $X_0$ is
${\rm supp}\mu$-invariant. This contradicts the first condition
in hypothesis {\bf H}. It follows $X_0=\emptyset$, hence $\nu$ has
no atom.

\medskip

For the second assertion we can repeat the first part of the above argument.
Thus we consider the set $W$ of affine subspaces $L$ of minimal dimension
such that $\nu(L)>0$. From the definition of $W$: $\nu(L\cap L')=0$ if
$L,L'\in W$ and $L\not=L'$. Hence $\sum_{L\in W}\nu(L)\le 1$, and there
exists $N\in W$ with  $a'=\nu(N)=\sup_{L\in W}\nu(L)$. Let
$W_0=\{ N\in W;\nu(N)=a' \}$. Then as above $W_0$ is finite
and ${\rm supp}\mu$-invariant, hence $H_{\mu}$-invariant,
where $H_{\mu}$ is the closed subgroup of $H$ generated by
${\rm supp \mu}$. It follows
that $H_{\mu}'=\{ h\in H:  hN=N, \forall N\in W_0\ \}$ is a finite index
subgroup of $H_{\mu}$. Let  $h=(b,g)$ be an element of $H_{\mu}'$
with $|g|<1$, and $h^+$ its unique fixed point. Then for every
$v\in V$, $\lim_{n\to\8}h^n v=h^+$. In particular, let $v_i$
be a point of $N_i\in W_0$ ($1\le i \le p$). Then $\lim_{n\to\8} h^n v_i = h^+\in V_i$
($1\le i \le p$). It follows $h^+\in\bigcap_{N_i\in W_0}N_i$, in particular
$\bigcap_{N_i\in W_0}N_i \not= \emptyset$ is a ${\rm supp}\mu$-invariant affine subspace.
This contradicts the hypothesis, hence $W=\emptyset$, i.e. $\nu(L)=0$
for every affine subspace $L$ of $V$.
\end{proof}

 We complete
the  result of \cite{BD} concerning nondegeneracy of the tail measure
and following methods described in \cite{Gr2,Gol,B} we prove it under
hypothesis {\bf H}, without any further assumptions.
\begin{prop}
\label{2.6}
The tail measure $\L$ is nonzero. In particular, if $\mu$ satisfies hypothesis
{\bf H}, there exists $k>0$ with $\P\big[ |R|>t \big]\ge k t^{-\a}$ for $t$ large enough.
\end{prop}
\begin{proof}
Define the backward process $ \breve  R_n$:
\begin{eqnarray*}
 \breve R_0&=& 0,\\
 \breve R_n &=& \pi_V \big( ( Q_1, M_1)\cdot\ldots\cdot( Q_n, M_n) \big)=
 Q_1 +  \breve \Pi_1  Q_2 + \cdots +  \breve\Pi _{n-1}  Q_n,
\end{eqnarray*}
where
$ \breve \Pi_k =  M_1 \cdot\ldots\cdot  M_k$.
Recall that $\breve  R_n$ converges pointwise to $ R$, and
$ R =  \breve R_n  + \breve  \Pi_n  \breve R^{n}$,
where
$  \breve R^{n} = \sum_{k=n+1}^\8 \big(  M_{n+1} \cdot\ldots\cdot  M_{k-1} \big)
  Q_k
$,
hence  for any $n$, $ \breve R^{n}$ and $ R$ have the same distribution.

Fix two positive numbers $\eta$ and $\d$ and a point $u\in{\rm supp \nu}$. For any
ball $U$ of center $u$ and radius $\d$, $\eps = \p{
R\in U}$ is positive. We have, using independence of $\breve R^n$ and $(\breve R_i, \breve \Pi_i )$ for $i<n$,
\begin{multline*}
\p{\inf_{x\in U}| \breve R_n +  \breve \Pi_n x|>t \mbox{ for some }n}
= \sum_n\p{\max_{i<n}\inf_{x\in U}| \breve R_i +  \breve \Pi_i x|\le t \mbox{ and }\inf_{x\in U}| \breve R_n +  \breve \Pi_n x| > t}\\
=\frac 1{\eps} \sum_n\p{\max_{i<n}\inf_{x\in U}| \breve R_i +  \breve \Pi_i
x|\le t \mbox{ and }\inf_{x\in U}| \breve R_n +  \breve \Pi_n x| > t}
\p{ \breve R^{n}\in U}\\
\le \frac 1{\eps} \sum_n\p{\max_{i<n}\inf_{x\in U}|\breve  R_i + \breve  \Pi_i
x|\le t \mbox{ and }\inf_{x\in U}|\breve  R_n + \breve  \Pi_n x| > t
\mbox{ and } | R| > t} \\
\le \frac{1}{\eps}\p{ | R| > t}.
\end{multline*}
Since $\P\big[ Mu+Q=u\big]<1$,
there exist a  positive number $\eta$ such that
$$\th=\p{| Q + ( M -I)u| > 2\eta}>0.$$
Moreover  there is a large number $N$ such that
$$\p{| M|\ge N}\le \frac {\theta}2.$$
Choose  $\d=\frac{\eta}{N+1}$ and define
$$
U_n =  \breve R_n +  \breve \Pi_n u - ( \breve R_{n-1} +  \breve \Pi_{n-1} u ) = \breve \Pi_{n-1}(
Q_n + ( M_n-I)u).
$$
Then
\begin{eqnarray*}
\p{| R|>t}&\ge&\eps \p{\inf_{x\in U}|\breve  R_n +  \breve \Pi_n x|>t \mbox{ for some }n}\\
&\ge&\eps\p{| \breve R_n+\breve \Pi_n u|-| \breve \Pi_n|\d>t \mbox{ for some }n}\\
&\ge& \eps\p{|U_n| -(| \breve \Pi_n|+ | \breve \Pi_{n-1}|)\d>2t \mbox{ for some }n}\\
&=& \eps\p{|\breve \Pi_{n-1}|(| Q_n+( M_n-I)u| - (| M_n|+1)\d) >2t \mbox{ for some }n}
\end{eqnarray*}
We define for $n\ge 0$, $Y_n= | Q_n+( M_n-I)u| - (| M_n|+1)\d$
and we observe that, using independence
\begin{eqnarray*}
\p{|\breve \Pi_{n-1}|Y_n  >2t \mbox{ for some }n}
&=& \p{| \breve\Pi_{n-1}|Y_0 >2t \mbox{ for some }n}\\
&\ge & \p{| \breve\Pi_{n-1}|>2t/\eta \mbox{ for some }n}\p{Y_0 >\eta}\\
&\ge&  \p{\max_{n\ge 1}| \breve\Pi_{n-1}|>2t/\eta }\p{Y_0 >\eta}
\end{eqnarray*}
On the other hand
\begin{eqnarray*}
\p{Y_0>\eta} &\ge& \p{|Q+Mu-u|>2\eta \mbox{ and } |M|<N}\\
&\ge& \p{|Q+Mu-u|>2\eta} - \p{|M|\ge N}\\
&\ge& \frac {\th}2.
\end{eqnarray*}
Since $\E\big[\log |M_1|\big]<0$, $\E\big[|M_1|^\a\big]=1$, we can use Cramer estimate
of ruin (see \cite{F}, p. 411) for $\p{\max_{n\ge 1}|\breve\Pi_{n-1}|>2t/\eta}$. This gives the
existence of $C>0$ (depending of $\mu_r$ only) such that $\p{\max_{n\ge 1}|\breve\Pi_{n-1}|>2t/\eta}\ge C \eta^\a t^{-\a}$.
Finally
$$\p{| R|>t}\ge \frac{\theta \eps C\eta^\a}2t^{-\a}.$$
Hence we can take $k= \frac {\th \eps C\eta^\a}2$. By definition of $\Lambda$, $\Lambda\not=0$.
\end{proof}
\begin{cor}
\label{cor-2.7}
 The function $x\mapsto |x|^\a$ is not $\nu$-integrable.
\end{cor}
\begin{proof}
The relation $\int_V|x|^\a \nu(dx) =\8$ follows from $\P\big[|R|>t \big]\ge kt^{-\a}$.
\end{proof}
\begin{cor}
\label{2.9}
Assume furthermore that there is no ${\rm supp}\mu$-invariant affine subspace. Then, for every affine subspace
$W$ of $V$, $\Lambda(W)=0$.
\end{cor}
\begin{proof}
We use the formula for $\Lambda$ obtained in \cite{BD}, Theorem 1.6:
$$\Lambda=\frac 1{m_\a}\int_V g(\nu-\ov\mu*\nu)\l^\a(dg),$$
where $\l^\a$ is a Radon measure on $G_{\ov\mu}$ equivalent to the Haar measure  of $G_{\ov\mu}$ and $m_\a$
was defined in hypothesis {\bf H}. Let $W$ be an
affine subspace of $V$ and $X\subset W$ a compact subset with $\Lambda(X)>0$. Then
$$\int_V (\nu-\ov\mu*\nu)(g^{-1}X)\l^\a(dg)>0.$$
Hence, for some $g\in G_{\ov\mu}$:
$$ (\nu-\ov\mu*\nu)(g^{-1}X)>0$$
and
$$ \nu(g^{-1}X) = \ov\mu*\nu(g^{-1}X) + (\nu-\ov\mu*\nu)(g^{-1}X)>0.
$$
In particular $\nu(g^{-1}W)>0$, which contradicts Lemma \ref{2.5}. The conclusion follows.
\end{proof}

The properties of $\eta_v$ that will be useful are contained in
the following Lemma
\begin{lem}
\label{3.6'}
Assume that $\mu$ satisfies hypothesis {\bf H}. Then the sequence
\begin{equation}
\label{zn}
Z_n=\sum_{k=1}^n M_{k-1}\ldots M_0
\end{equation}
converges $\P$-a.e.
to $Z=\sum_{k=1}^\8 M_{k-1}\ldots M_0$. For any $v\in V$, the law
of $Z^*v$ is the unique stationary measure $\eta_v$ of the Markov chain
on $V$ defined by \eqref{rekursja2}.

If $v\not= 0 $  then for any
$x\in V$, $\P[M^*(x+v)=x]<1$. In particular the recursion \eqref{rekursja2} satisfies
hypothesis {\bf H}, the measure $\eta_v$ has no atoms and
 has all moments smaller than $\a$, i.e.
$\eta_v(\tau^\t)<\8$ for $\t<\a$ and $\eta_v(\tau^\a)=\8$.

Moreover for every $c\in Z({\ov\mu})$, the centralizer of $G_{\ov\mu}$ in $G$, $\eta_{c^*v}(f)=\eta_v(f\circ c^*)$
for $f\in C_b(V)$.
\end{lem}
\begin{proof}
If suffices to show the convergence of $\is{Z_n(\omega)x}y$
for any $x,y\in V$. But  $\is{Z_n(\omega)x}y=\is x{Z_n^*(\omega)y}$
and since $\E[\log|M^*|]<0$
$$
Z^*(\omega) = \lim_{n\to\8}\Big( \sum_{k=1}^n M_0^*(\omega)\ldots M_{k-1}^*(\omega) \Big)
$$ exists $\P$-a.e.
and also the existence and uniqueness of $\eta_v$ is clear
(see \cite{BD} for some further explanations).

 If $x\in V$ satisfies $M^*(x+v)=x$ $\P$-a.e. then for any $g\in{\rm supp}\ov\mu$,
 $v=(g^*)^{-1}x-x$. Therefore putting into the last equation two arbitrary
 elements belonging to the support of $\ov\mu$, say $g$ and $g'$,
 we obtain $(g^*)^{-1}x=(g'^*)^{-1}x$.
  If $x\not=0$ this implies
 $|g|=|g'|$, which contradicts hypothesis {\bf H}, since $|g|$ ($g\in{\rm supp}\ov \mu$) takes at
 least two different values. Thus,  hypothesis {\bf H} is valid and by Lemma \ref{2.5}, $\eta_v$ has no atoms and $\eta_v(\tau^{\th})<\8$ if
 $\th<\a$. Also $\eta_v(\tau^\a)=\8$, by Corollary \ref{cor-2.7}

 For the last assertion notice that if $f={\bf 1}_U$ for some $U\subset V$
 and  $c\in Z({\ov\mu})$, then
 $$
 \eta_{c^*v}({\bf 1}_U) = \P\big[ Z^*c^*v\in U \big] =
 \P\big[ Z^* v\in (c^*)^{-1}\cdot U \big] = \eta_v({\bf 1}_U\circ c^*).
 $$
\end{proof}

\section{Fourier operators and their spectral properties}
\label{section2}

\subsection{Analysis  of the  Fourier operators}
On continuous  functions on $V$ we introduce as in \cite{LP}
the seminorm
$$
[f]_{\e,\l}= \sup_{x\not= y} \frac{|f(x)-f(y)|}{\tau(x-y)^\e(1+\tau(x))^\l(1+\tau(y))^\l}
$$
and the
 two norms
\begin{equation}
\begin{split}
|f|_{\th} &= \sup_x \frac{|f(x)|}{(1+\tau(x))^{\th}}\\
\|f\|_{\th,\eps,\l} &= |f|_{\th}+ [f]_{\eps,\l}\\
\end{split}
\end{equation}
Notice, that if $\l+\e\le\t$ (that will be  always assumed), then $[f]_{\e,\l}<\8$ implies $|f|_\t<\8$.
\medskip

Define Banach spaces
\begin{equation}
\label{banach-space}
\begin{split}
\C_{\th}&=\Big\{ f:\; |f|_{\th} < \8 \Big\},\\
\B_{\tel}&=\Big\{ f:\; \|f\|_{\tel}  < \8 \Big\}.
\end{split}
\end{equation}
On $\C_{\th}$ and $\B_{\tel}$ we  consider the transition operator
\begin{equation}
Pf(x) =
\E\Big[ f(Mx+Q)\Big] = \int_H f(gx+b)\mu(dh),
\end{equation} where $(Q,M)$ is a random variable distributed according to the measure $\mu$.
We consider also the Fourier operator $P_v$ defined by
\begin{equation}
\label{pt}
P_{v}f(x) =
\E\Big[ e^{i\is{v}{Mx+Q}}f(Mx+Q)\Big] = P(\chi_v f)(x),
\end{equation}
where $v\in V$  and $\chi_v(x)=e^{i\is vx}$.
Notice $P_{0}=P$.
We will prove later (Lemma \ref{2e4} and Proposition \ref{u2}) that the operators $P_v$ are bounded
on $\B_{\tel}$ for appropriately chosen parameters $\tel$.
It follows from the inequality in Proposition \ref{u2} below and
 the Theorem of Ionescu Tulcea and Marinescu \cite{ITM} that
 all the operators $P_{v}$  have at most  finitely many eigenvalues
of modulus $1$, they have finite multiplicity, and the rest of the spectrum
is contained in a ball centered at the origin of radius less than $1$.
Moreover, for $|v|$ small, the perturbation theorem of Keller and Liverani \cite{KL} provides
uniform control of these spectrums.  Namely the spectrum and spectral
properties of $P$ can be approximated in an appropriate  way by the corresponding features of operators $P_{v}$.
All the details will be given below. For an operator $A$ we denote by $\sigma(A)$ its spectrum
and by $r(A)$ its spectral radius.
 After a few lemmas we will apply \cite{KL} to our
situation.

\medskip

For random variables $\{X_n^x\}$ defined in \eqref{rekursja} we consider partial sums $S_n^x = \sum_{k=1}^n X_k^x$.

The following simple lemma is the basis of the use
of spectral methods in limit theorems for functionals of Markov chains.
\begin{lem}
\label{1e4} We have
\begin{equation*}
P_v^n f(x) = \E\Big[ \chi_v(S_n^x)f(X_n^x)\Big]
\end{equation*}
\end{lem}
\begin{proof}
If $n=1$, then the formula above coincides with  definition \eqref{pt}. Assume the result holds for $n$.
If $(Q,M)$ is independent of $S_n^x$ we write
\begin{eqnarray*}
P_v^{n+1}f(x)&=& \E\Big[\chi_v(Mx+Q)P_v^n f(Mx+Q)\Big] \\
&=& \E\Big[ \chi_v(Mx+Q) \chi_v(S_n^{ Mx+Q})f(X_n^{Mx+Q})\Big]
= \E\Big[ \chi_v(S_{n+1}^x
) f(X_{n+1}^x)\Big],
\end{eqnarray*}
that completes the proof.
\end{proof}

We will need the following inequality, valid for any $\b\in [0,1]$:
\begin{equation}
\label{h4}
|e^{i\langle x,y\rangle}-1|\leq 2\tau(x)^{\b}\tau(y)^{\b}, \qquad \mbox{for every $0<\b \le 1$.}
\end{equation}

\begin{lem}
\label{2e4}
For every $v\in V$, $n\in \N$ and $\th<\a$. We have
$$ |P_v^n f|_\t\le D|f|_\t $$
with $D=3^\t\big[2+\sup_n\E[\tau(X_n)^\t]\big]<\8$.
\end{lem}
\begin{proof}
Notice first that
\begin{equation}
\label{skr}
X_n^x=X_n^y+\Pi_n(x-y),
\end{equation}
where $\Pi_n=M_n M_{n-1}\ldots M_1$.
Therefore by Lemma \ref{1e4}, for every $x\in V$ we have
\begin{eqnarray*}
\frac{|P_v^n f(x)|}{(1+\tau(x))^\t}&\le& \E\bigg[ \frac{|f(X_n^x)|}{(1+\tau(X_n^x))^\t}\cdot \frac{(1+\tau(X_n^x))^\t}{(1+\tau(x))^\t}\bigg]\\
&\le& |f|_\t\E\bigg[ \frac{(1+\tau(X_n)+ |\Pi_n|\tau(x))^\t}{(1+\tau(x))^\t}\bigg]
\le 3^\t|f|_\t \big( 1+ \E\tau(X_n)^\t  + \kappa^n(\th)\big).
\end{eqnarray*}
Since $\th<\a$,
in view of \eqref{moments}
the factor of $|f|_\t$ above is bounded by  $D=3^\t\big[2+\sup_n\E[\tau(X_n)^\t]\big]<\8$ and the lemma follows.
\end{proof}
\begin{prop}
\label{u2}
Assume  $2\l+\e < \a$, $\e<1$ and $\t <2\l$.
Then there exist constants
$C_1,C_2$ and $\ro<1$ independent of $v$  such that for every $n\in\N$, $f\in \B_{\tel}$, $v\in V$
$$ [P_v^nf]_{\e,\l} \le C_1\ro^n [f]_{\e,\l} + C_2 \tau(v)^\e|f|_\t. $$
\end{prop}
\begin{proof} We have
$$P_v^nf(x) - P_v^nf(y)
= \E \Big[ \chi_v({S_n^y})(f(X_n^x)-f(X_n^y))\Big]
+ \E \Big[\big(\chi_v({S_n^x}) -\chi_v({S_n^y}) \big)f(X_n^x)\Big]
$$
Let us denote these two functions above by $\D_1$ and $\D_2$, respectively, and estimate first $\D_1$
\begin{multline*}
\frac{\big|\D_1(x,y) \big|}{\tau(x-y)^\e(1+\tau(x))^\l(1+\tau(y))^\l}
\le [f]_{\e,\l}\cdot\E\bigg[\frac{\tau(X_n^x - X_n^y)^\e (1+\tau(X_n^x))^\l (1+\tau(X_n^y))^\l
}{\tau(x-y)^\e(1+\tau(x))^\l(1+\tau(y))^\l}\bigg]\\
\le  [f]_{\e,\l}\E\bigg[ \frac{|\Pi_n|^\e (1+\tau(X_n)+ |\Pi_n|\tau(x))^\l
(1+\tau(X_n)+ |\Pi_n|\tau(y))^\l
}{(1+\tau(x))^\l(1+\tau(y))^\l}\bigg]\\
\le 3^{2\l} [f]_{\e,\l}\cdot\E \Big[|\Pi_n|^\e (1+\tau(X_n)^\l+ |\Pi_n|^\l)^2 \Big]
\end{multline*}
Expanding the expression in brackets we obtain a sum of 6 factors of the form
$
\E\big[ |\Pi_n|^\b \tau(X_n)^\g\big]
$
for  $\b+\g \le \e+2\l<\a$.
Applying the H\"older inequality with parameters $p=\frac{\b+\g}\b$,
$q=\frac{\b+\g}\g$, in view of  \eqref{moments}, we have
$$
\E\Big[ |\Pi_n|^\b \tau(X_n)^\g\Big] \le
\kappa^{\frac np}(\b+\g)\Big(\E\big[  \tau(X_n)^{\b+\g}\big]\Big)^{\frac1q}
= C_{\b,\g} \ro_{\b,\g}^n,
$$ for  $\ro_{\b,\g}=\kappa(\b+\g)^{\frac 1p}$,
 strictly smaller than 1. Therefore
if $C_1 = 3^{2\l+2}\sup_{\b,\g}C_{\b,\g}$ and $\ro = \sup_{\b,\g}\ro_{\b,\g}<1$, then
\begin{equation}
\label{pa31}
\frac{\big|\D_1(x,y) \big|}{\tau(x-y)^\e(1+\tau(x))^\l(1+\tau(y))^\l}
 \le  C_1 \ro^n [f]_{\e,\l}.
\end{equation}
Now we are going to estimate $\D_2$. Define the random variable $B_n = 1+|\Pi_1|+\cdots+|\Pi_n|$,
then for $\d<\min\{1,\a\}$ we obtain
$$
\E\big[ B_n^\d\big] = \E \big[ 1+|M_n| B_{n-1} \big]^\d \le 1+\kappa(\d)\E\big[ B_{n-1}^\d\big] \le \sum_{j=0}^n \kappa^j(\d).
$$
Therefore
$$
\sup_n \E B_n^\d = \frac 1{1-\kappa(\d)}<\8.
$$
Assume $\tau(y)\ge \tau(x)$.
Applying \eqref{h4}, we write
$$
\big| \chi_v(S_n^x)-\chi_v(S_n^y) \big| \le 2\tau(v)^{\eps}\tau(Z_n)^{\eps}\tau( x-y )^{\eps}
\le 2\tau(v)^{\eps} B_n^{\eps} \tau(x-y)^{\eps},
$$ therefore
\begin{multline*}
\frac{|\D_2(x,y)|}
{\tau(x-y)^\e(1+\tau(x))^\l (1+\tau(y))^\l}
 \le 2 \tau(v)^{\e} |f|_\t \E\bigg[\frac{B_n^{\e}\tau(x-y)^{\e} (1+\tau(X_n^x))^\t}
{\tau(x-y)^\e(1+\tau(x))^\l (1+\tau(y))^\l} \bigg]\\
\le 2 \tau(v)^{\e} |f|_\t \E\bigg[\frac{B_n^{\e}
(1+\tau(X_n)+|\Pi_n|\tau(x))^\t}
{(1+\tau(x))^\l (1+\tau(y))^\l} \bigg]
\le 2\cdot 3^\t \tau(v)^{\e} |f|_\t \E \Big[ B_n^{\e}  \tau(X_n)^\th + B_n^{\e} \big( 1+ |\Pi_n|^\t\big)\Big]
\end{multline*}
Applying as before the H\"older inequality we prove that the expression in brackets above is
bounded and
\begin{equation}
\label{pa32}
\frac{|\D_2(x,y)|} {\tau(x-y)^\e(1+\tau(x))^\l (1+\tau(y))^\l} \le C_2 \tau(v)^{\e} |f|_\t,
\end{equation}
with $C_2 = 2\cdot 3^\t \sup_n\E\big[ B_n^{\eps}\tau(X_n)^\t + B_n^{\eps}(1+|\Pi_n|^\t)\big]$.
Finally combining \eqref{pa31} and \eqref{pa32} we obtain the Lemma.
\end{proof}
\begin{lem}
\label{p9}
If  $\l+2\e<\t<\a$ and $\d\le\eps$, then there exists a  constant $C$,   such that for every
 $\g$ satisfying $\l+2\e\le \g\le \th$ and $v,w\in V$:
  $$ |(P_{v}-P_{w})f|_\g \le C\tau(v-w)^\d  \|f\|_{\tel}.  $$
\end{lem}
\begin{proof}
 Using \eqref{h4} we have
\begin{multline*}
\bigg| \frac{(P_{v}-P_{w})f(x)}{(1+\tau(x))^\g} \bigg|
\le \bigg| \frac{(P_{v}-P_{w})(f-f(0))(x)}{(1+\tau(x))^\g}\bigg| +\bigg| \frac{(P_{v}-P_{w})(f(0))(x)}{(1+\tau(x))^\g}\bigg|\\
\le\E\bigg[ \Big|1 - e^{i\is{v-w}{Mx+Q}}\Big|\cdot\frac{|f(Mx+Q)-f(0)|}{(1+\tau(x))^\g}\bigg]
 + |f(0)|\cdot \E \bigg[\frac{|1 - e^{i\is{v-w}{Mx+Q}}|}{(1+\tau(x))^\g}\bigg]\\
\le 2[f]_{\e,\l} \tau(v-w)^\d \cdot\E \bigg[ \tau(Mx+Q)^\d  \cdot
\frac{(\tau(Mx+Q)^{\e} (1+\tau(Mx+Q))^\l}{(1+\tau(x))^\g}\bigg]\\
+2 \tau(v-w)^\d |f(0)|\cdot \E \bigg[\frac{\tau(Mx+Q)^\d}{(1+\tau(x))^\g}
\bigg] \le C \tau(v-w)^\d \|f\|_{\tel},
\end{multline*}
where $C=2\sup_x\E\big[(\tau(Mx+Q)^{\l'}+1)/(1+\tau(x))^{\l'}\big]$ for $\l'=\l+\eps+\d<\th$.
\end{proof}
\begin{lem}
\label{3.12}
The unique eigenvalue of modulus 1 for
 $P$ acting on $\C_\t$ is 1. The corresponding eigenspace is $\C 1$
and the projection on $\C 1$ along
the hyperplane ${\rm Ker} \nu = \{f\in \C_\t:\;\nu(f)=0 \}$
is given by the map $f\mapsto \nu(f)$.
\end{lem}
\begin{proof}
Of course constant functions are eigenfunctions of $P$ with eigenvalue 1
and $P$ acts on $\C_\th$ in view of Lemma \ref{2e4}.
 In fact there are no other elements of $\C_\t$ satisfying $Pf=f$. Indeed, let $f$
 be such a function, then
 for every $x\in V$,
 $\lim_{n\to\8}P^n f(x)= f(x)$. On the other hand,
$\lim_{n\to\8}P^n f(x)= \lim_{n\to\8} \E f(X_n^x)= \nu(f)$ (recall that
the law of   $\{X_n^x\}$ tends in distribution to $\nu$).
Hence $f(x) = \nu(f)$ for every $x\in V$, and $f$ must be a constant.
Furthermore, we observe that $\C_\t={\rm Ker}\nu \oplus \C 1$ and $f=(f-\nu(f))1+\nu(f)1$.
The assertion for the projection of $f$ follows.

\medskip

To prove that there are no other
eigenvalues of modulus 1 we proceed similarly. Assume that for some $z$ of modulus 1 and a nonzero
function $f\in \C_{\th}$ we have $Pf=z f$. Then $\lim_{n\to\8}P^n f(x)= \nu(f)$, but
$P^n f(x)= z^n f(x)$ and if $\eta$ would be different than 1, the sequence $z^n f(x)$, for every $x$
such that $f(x)\not=0$, couldn't converge to a constant. This implies $z=1$ and finishes the proof.
\end{proof}

\begin{lem}
\label{3.13}
Assume that no affine subspace of $V$ is supp$\mu$-invariant. Then
for every $v\not=0$,  the equation $P_v f = zf$, $|z|=1$, $f\in \B_{\tel}$
implies $f=0$. In particular  the spectral radius of $P_v$ is strictly smaller than 1.
\end{lem}
\begin{proof}
We proceed as in \cite{GL}.  Assume that
\begin{equation}
\label{dr1}
P_v f = z f
\end{equation}
for some  nonzero $f\in \B_{\tel}$
and $|z|=1$. Then the function $f$ is bounded. Indeed for every $n$
$$
|f(x)| = |z^n f(x)| \le P^n(|f|)(x)
$$
hence
$$
|f(x)|\le \lim_{n\to\8} P^n(|f|)(x) =\nu(|f|).
$$
Next observe $\is{\nu}{\nu(|f|)-|f|}=0$, therefore,
since $f$ is continuous, on the support
of $\nu$ the function $|f|$ is equal to its maximum  and without
any loss of generality we may assume that this maximum is 1. A
convexity argument, Lemma \ref{1e4} and \eqref{dr1} imply that for
every $n$ and $x\in \mbox{supp}\nu$.
$$
z^n f(x) = e^{i\is{v}{S_n^x}}f(X_n^x)\qquad \P \mbox{ a.s.}
$$ Hence for any $x,y\in \mbox{supp} \nu$
\begin{equation}
\label{sr1} \frac{f(x)}{f(y)} e^{i\is{v}{Z_n(y-x)}} =
\frac{f(X_n^x)}{f(X_n^y)},
\end{equation}
where $Z_n$ was defined in \eqref{zn}.
Observe that taking $p=\frac{2\l+\e}{\e}$ and $q=\frac{2\l+\e}{2\l}$, by the H\"older inequality we have
\begin{multline*}
\limsup_{n\to\8} \E\bigg| \frac{f(X_n^x)}{f(X_n^y)}-1 \bigg| \le
[f]_{\e,\l} \limsup_{n\to\8} \E\Big[ \tau(X_n^x-X_n^y)^\e (1+\tau(X_n^x))^\l(1+\tau(X_n^y))^\l \Big]\\
= [f]_{\e,\l} \limsup_{n\to\8} \E\Big[ \tau(M_n\cdots M_1(x-y))^\e (1+\tau(X_n^x))^\l(1+\tau(X_n^y))^\l \Big]\\
\le [f]_{\e,\l}  \tau(x-y)^\e \limsup_{n\to\8}\Big(\E|M_n\ldots M_1|^{2\l+\e} \Big)^{\frac 1p}\cdot
\limsup_{n\to\8} \bigg(  \E \Big[ (1+\tau(X_n^x))^{\l+\frac{\e}{2}}(1+\tau(X_n^y))^{\l+\frac{\e}{2}} \Big] \bigg)^{\frac 1q}.
\end{multline*}
In view of \eqref{moments} the last term is finite and since $2\l+\e<\a$,  $\lim_{n\to\8}\kappa^{\frac np}(2\l+\e)=0$,
hence
$$
\lim_{n\to\8} \E\bigg| \frac{f(X_n^x)}{f(X_n^y)}-1 \bigg| = 0.
$$
Therefore for  $\P$ a.e. trajectory $\omega$ there exists a sequence $\{n_k\}$ such that
$$
\lim_{n_k\to\8}\frac{f(X_{n_k}^x)}{f(X_{n_k}^y)} = 1.
$$ Notice that, in view of Lemma \ref{3.6'}
$\lim_{n\to\8} Z_n(\omega)=Z(\omega) $ exists $\P$-a.e. Hence passing
with $n$ to infinity in \eqref{sr1} we obtain $ \frac{f(x)}{f(y)}=
e^{-i\is{v}{Z(\omega)(x-y)}} $ for $\omega\in\Omega_0$ with $\P(\Omega_0)=1$. Then for every
$\omega\in\Omega_0$,  $ e^{-i\is{Z^*(\omega)v}{(x-y)}} \frac{f(x)}{f(y)}=1$.
We are
going to prove that this leads to a contradiction whenever $v\neq
0$. We choose  points $x_j,y_j\in\mbox{ supp
}\nu$, $j=1,\ldots,d$ with $v_j=x_j-y_j$ spanning  $V$ as a vector space.
Such points exist because the support of $\nu$ as a set invariant
under the action of supp$\mu$ is not  contained in some
proper affine subspace of $V$.
Let
$\eta_v$ be the law of
 $W(\omega) = Z^*(\omega)v$.
Then for every $j$ the support of $\eta_v$ is contained in the union of affine hyperplanes $\bigcup_{n\in\Z}\{H_j+ns_jv_j\}$, where $H_j$ is some hyperplane
orthogonal to $v_j$ and $s_j$ are appropriately chosen constants. Taking intersection of all such sets
defined for every $j$ we conclude that
supp$\eta_v$ is contained is some discrete  set of points, hence
${\rm supp}\eta_v$ is discrete. This contradicts Lemma \ref{2.5}.

\medskip

For the last
assertion we observe that in view of Theorem of Ionescu Tulcea and Marinescu \cite{ITM},
if $z$ belongs to the spectrum of $P_v$ and $|z|=1$ then $z$ is an eigenvalue of $P_v$.

\end{proof}
\subsection{A perturbation theorem}
 For  $c\in Z(\ov\mu)\cup\{0\}$
we denote $P_{c, v}=P_{c^*v}$ and we write $c\to 0$
for $|c|\to 0$.
We observe that $Z(\ov\mu)$, the centralizer of $G_{\ov\mu}$ in $G$, contains $\R^*_+,\ Z_{\ov\mu}$ and
$$P_{c,v}f(x)=\int_H \chi_v(c(gx+b))f(gx+b)\mu(dh).$$

\medskip

In view of Lemmas \ref{2e4},  \ref{p9}  and Proposition \ref{u2} we may use the perturbation theorem
of Keller and Liverani \cite{KL} for the family $P_{c,v}$
(the hypothesis concerning the essential radius is fulfilled by a result of Hennion \cite{H}, Corollary 1).
Their result is stated for the case when the parameter is real, but  what they really
use is the H\"older continuity in   Lemma \ref{p9}, which is valid also in our more general settings. Then we obtain the following.
\begin{prop}
\label{kel} Assume $\e<1$, $\l+2\e < \th < 2\l$, $2\l+\e < \a$, $v\in V$ is fixed then there exist $t_0>0$,
$\d>0$, $\rho<1-\d$ such that for every $c\in Z(\ov\mu)\cup\{0\}$ with $|c|\le t_0$:
\begin{itemize}
\item The spectrum of  $P_{c,v}$ acting on $\B_{\tel}$
 is contained in ${\mathcal D}= \{ z:\; |z|\le \ro \}\cup \{ z:\; |z-1|<\d \} $.
\item The set $\s(P_{c,v})\cap \{ z:\; |z-1|<\d \}  $ consists of exactly one eigenvalue $k(c,v)$,
the corresponding eigenspace is one dimensional and moreover $\lim_{c\to 0 }k(c,v) = 1$.
\item If $\pi_{c,v}$ is the projection of $P_{c,v}$ onto the mentioned above eigenspace, then
there exists an operator $Q_{c,v}$ such that $\|Q_{c,v}\|\le \rho$, $\pi_{c,v} Q_{c,v} = Q_{c,v} \pi_{c,v} = 0$ and for every $n$
$$ P_{c,v}^n f = k(c,v)^n\pi_{c,v}(f) + Q_{c,v}^n(f),\qquad f\in \B_{\tel}. $$
\item For any $z$ belonging to the complement of ${\mathcal D}$:
 $$\| (z-P_{c,v})^{-1}f \|_{\tel}\le D\|f\|_{\tel},$$ for some constant $D$ independent of $c$.
\end{itemize}
\end{prop}
Define for small $|c|$ the function $g_{c,v} = \pi_{c,v}(1)$.
 Then for every function $f$ belonging to $\B_{\tel}$,
 we define
$\nu_{c,v}(f)\in \C_\t$  by $\pi_{c,v}(f) = \nu_{c,v}(f)g_{c,v}$.
\begin{prop}
\label{lis22} Assume additionally that $\l+3\e < \th$, $2\l+ 3\e
<\a$. The identity embedding of $\B_{\tel}$ into
$\B_{\th,\e,\l+\e}$ is continuous and the decomposition
$P_{c,v}=k(c,v)\pi_{c,v}+Q_{c,v}$ coincides on both spaces. Moreover there
exist constants $D$ and $t_1$ such that for $|c|\le t_1$, $c\in Z(\ov\mu)\cup\{0\}$  we have, if $\tau(v)\le 1$
\begin{eqnarray*}
&&\mbox{i)\ \ }\|(P_{c,v}-P)f\|_{\t,\e ,\l+\e} \le D |c|^{\e}\|f\|_{\tel};\\
&&\mbox{ii)\ \ } \|(k(c,v)\pi_{c,v} - \pi_0)f\|_{\t,\e ,\l+\e} \le D |c|^{\e}\|f\|_{\tel};\\
&&\mbox{iii)\ \ } \|(\pi_{c,v} - \pi_0)f\|_{\t,\e ,\l+\e} \le D |c|^{\e}\|f\|_{\tel};\\
&&\mbox{iv)\ \ } \|(Q_{c,v}-Q)f\|_{\t,\e ,\l+\e} \le D |c|^{\e}\|f\|_{\tel};\\
&&\mbox{v)\ \ } \|g_{c,v} - 1\|_{\tel}\le D |c|^{\e}; \\
&&\mbox{vi)\ \ } |k(c,v)-1|\le D |c|^{\e}\\
&&\mbox{vii)\ \ $\nu_{c,v}$ is bounded on $\B_{\tel}$ and  }
|\nu_{c,v}(f)-\nu(f)|\le D |c|^{\e}\|f\|_{\tel}.\qquad\qquad\qquad\qquad\qquad
\end{eqnarray*}
\end{prop}
\begin{proof}
The triple $(\th,\e,\l+\e)$ satisfies assumptions of Proposition \ref{kel}, and of course
$$\|\cdot\|_{\th,\e,\l+\e}\le \|\cdot\|_{\th,\e,\l}, $$ therefore considering the family
of operators $\{P_{c,v}\}$ on both Banach spaces $\B_{\th,\e,\l+\e}$ and $\B_{\th,\e,\l}$ we obtain the same decomposition
of $P_{c,v}$.

To prove i), in view of Lemma \ref{p9}, it is enough to estimate
\begin{multline*}
 \big[(P_{c,v}-P)(f)\big]_{\e,\l+\e}
=\sup_{x\not= y} \frac{|(P_{c,v}-P)f(x) - (P_{c,v}-P)f(y)|}{\tau(x-y)^{\e }(1+\tau(x))^{\l+\e}(1+\tau(y))^{\l+\e}}\\
\le\sup_{x\not= y} \frac{\E\big[|\chi_v(cX_1^x)-1||f(X_1^x) - f(X_1^y)|\big]}{\tau(x-y)^{\e }(1+\tau(x))^{\l+\e}(1+\tau(y))^{\l+\e}}
+ \sup_{x\not= y} \frac{\E\big[ |\chi_v(cX_1^x)-\chi_v(cX_1^y)||f(X_1^y)|\big]}{\tau(x-y)^{\e }(1+\tau(x))^{\l+\e}
(1+\tau(y))^{\l+\e}} \\ = \D_1+\D_2.
\end{multline*}
Next  we have, using \eqref{h4}
\begin{eqnarray*}
\D_1 &\le& 2 |c|^{\e} \|f\|_{\tel}
 \sup_{x\not= y} \E\Bigg[ \frac{\tau(X_1^x)^{\e}\tau(X_1^x - X_1^y)^\e(1+\tau(X_1^x))^\l(1+\tau(X_1^y))^\l}
{\tau(x-y)^{\e }(1+\tau(x))^{\l+\e}(1+\tau(y))^{\l+\e}} \Bigg]\\
&\le& 2|c|^{\e} \|f\|_{\tel}
\E\Big[(|M_1|+\tau(Q_1))^{ \e}|M_1|^\e(1+|M_1|+\tau(Q_1))^{2\l} \Big]
\end{eqnarray*}
Reasoning as in Proposition \ref{u2}, one can prove that the expected value above is finite.
Similarly we estimate $\D_2(x,y)$:
\begin{eqnarray*}
\D_2 &\le&
 \sup_{x\not= y} \E\Bigg[ \frac{ |\chi_v(cX_1^x)-\chi_v(cX_1^y)|\cdot\big(|f(X_1^y) - f(0)|+|f(0)|\big)}{\tau(x-y)^{\e }(1+\tau(x))^{\l+\e}
(1+\tau(y))^{\l+\e}}\Bigg]\\
&\le&  |c|^{\e}  \|f\|_{\tel}  \sup_{x\not= y}  \E\Bigg[\frac{\tau(X_1^x-X_1^y)^{\e}
\big( \tau(X_1^y)^\e(1+\tau(X_1^y))^\l+1 \big) }{\tau(x-y)^{\e }(1+\tau(x))^{\l+\e}
(1+\tau(y))^{\l+\e}}\Bigg]\\
&\le& D|c|^{\e}\|f\|_{\tel}\E\Big[ |M_1|^\e \big(1 + |M_1| + \tau(Q_1)\big)^{\l+\e} \Big].
\end{eqnarray*}
Again arguments from  Proposition \ref{u2} prove that the foregoing value is finite. Similarly we prove
$|(P_{c,v}-P)(f)|_\t \le 2 |c|^\e |f|_\t$ that gives i)

\medskip

In order to prove ii) and iii) we will use the fact that both $\pi_{c,v}$ and
$Q_{c,v}$ can be expressed  in terms of the resolvent of $P_{c,v}$. We follow arguments
in \cite{LP2}
\begin{equation*}
\begin{split}
k(c,v)\pi_{c,v} &= \frac1{2\pi i}\int_{|z-1|=\d'}z(z-P_{c,v})^{-1}dz,\\
\pi_{c,v} &= \frac1{2\pi i}\int_{|z-1|=\d'}(z-P_{c,v})^{-1}dz,\\
Q_{c,v} &= \frac1{2\pi i} \int_{|z|=\ro'}z(z-P_{c,v})^{-1}dz,
\end{split}
\end{equation*}
for appropriately chosen constants $\d'$ and $\ro'$. Then combining
the formulas above with $$ (z-P)^{-1} - (z-P_{c,v})^{-1} =
(z-P)^{-1}(P-P_{c,v})(z-P_{c,v})^{-1},$$ the point i) and  estimates of
the norm of resolvent (Proposition \ref{kel}) we conclude ii),
iii) and iv). v) is an immediate consequence of ii). To
prove vi) we write $$(k(c,v)-1)\pi _0=k(c,v)\pi_{c,v}-\pi _0-k(c,v)(\pi_{c,v}-\pi
_0),$$ apply $(k(c,v)-1)\pi _0$ to $1$ and we use ii) and iii).
Similarly, writing
$$(\nu_{c,v}(f)-\nu _0(f))1 =  \pi_{c,v}(f)-\pi _0(f) -  \nu_{c,v}(f)(g_{c,v}-1)$$
and applying iii) and v) we obtain vii), that finishes the proof.

\end{proof}

\begin{lem}
\label{cont}
For fixed $v$ the function
 $c\mapsto r(P_{c,v})$
 defined on $Z(\ov\mu)\cup\{0\}$ is continuous at 0. Moreover
 for any $f\in \B_{\tel}$, $c_0\in Z(\ov\mu)$ with $r_f(P_{c,v})=\limsup_n\|P_{c,v}^n f\|^{\frac 1n}$,
we have $\limsup_{c\to c_0}r_f(P_{c,v})=r_f(P_{c_0,v})$.
\end{lem}
\begin{proof}
From Proposition \ref{kel}, $r(P_{c,v})=|k(c,v)|$ for $c$ small, hence the continuity of $r(P_{c,v})$
follows from Proposition \ref{lis22}. Using again Proposition \ref{lis22}, for any
fixed $f$ and $n$, $\|P_{c,v}^n f\|^{\frac 1n}$ depends continuously on $c$.
Hence $\limsup_{n\to\8}\|P_{c,v}^n f\|^{\frac 1n}=r_f(P_{c,v})$ is upper semicontinuous in $c$.
\end{proof}
\subsection{Eigenfunctions of $P_{c,v}$}
 Proposition \ref{lis22} says that the dominant eigenvalues $k(c,v)$ of $P_{c,v}$ tend to 1 with rate at least
$|c|^\e$. However to prove our Main Theorem we will need more precise information concerning
the asymptotic expansion of $k(c,v)$,
that will be described in Theorem \ref{asymptkt}. For this purpose, following ideas
of  \cite{GL} we will express the eigenfunction corresponding to $k(c,v)$ in a more explicit
way.

\medskip

For $c\in Z(\ov\mu)\cup\{0\}$,
 let us define a family of operators on the Banach space $\B_{\tel}$:
$$T_{c,v} f(x) =
\E \Big[ e^{i\is{Q}{c^*(x+v)}}f(M^*(x+v))\Big] = \int_H\chi_b(c^*(x+v))f(g^*(x+v))\mu(dh).
$$
Then  $T_v=T_{0,v}$ and we have $T^*_v\eta_v = \eta_v$, where
$\eta_v$ is the stationary measure of the Markov chain  $\{W_n\}$ defined in \eqref{rekursja2}
(see also Lemma \ref{3.6'})

It turns out that the family $\{T_{c,v}\}$ satisfies assumptions of the perturbation theorem of \cite{KL}
and also the analogue of Lemma \ref{3.12} is valid in these settings, i.e. the set of peripherical
eigenvalues consists of one element 1 and $\eta_v$ is the projection onto the corresponding
one dimensional eigenspace.

We omit the argument, because this can be proved exactly in the same way as for the family $\{P_{c,v}\}$.
Therefore for small values
of $|c|$, the spectrum of $T_{c,v}$ intersected with some neighbourhood of 1, consists of exactly one
point $k'(c,v)$, which is the dominant eigenvalue of $T_{c,v}$. Let us denote the corresponding
projection by $\pi'_{c,v}$ and as before define, for every $f\in \B_{\tel}$, $\eta_{c,v}(f)$
to be the unique number such that $\pi'_{c,v}(f) =\eta_{c,v}(f)\pi'_{c,v}(1)$, hence
$\eta_{c,v}(1)=1$ and $T_{c,v}^*\eta_{c,v}=k'(c,v)\eta_{c,v}$.
Reasoning as in Proposition \ref{lis22} one can prove
\begin{prop}
\label{lis23} There exist constants $t_2$, $C'$ and $D'$ such that for $|c|\le t_2$, $c\in Z(\ov\mu)$,
 and every $f\in \B_{\tel}$ we have
$$
\big|\eta_{c,v}(f) - \eta_v(f)\big|\le C'|c|^\e \|f\|_{\tel}.
$$
In particular
$$ \|\eta_{c,v}\|_{\tel} \le 1+C'|c|^{\eps}\le D',$$
if $|c|\le t_2$.
\end{prop}
One  easily verifies that
for any $x$ the function $\chi_x$ is an element of $\B_{\tel}$ and moreover
\begin{equation}
\label{p15}
 \| \chi_x \|_{\tel} \le 1+2\tau(x)^{\e}
\end{equation}
It follows that for any $\eta\in \B_{\tel}$ the Fourier transform
$\wh{\eta}(x) = \eta(\chi_x)$ is well-defined and
$$|\wh{\eta}(x)|\le 1+2\tau(x)^{\eps}.$$

The following intertwining relation between $P_{c^*v}$ and $T^*_{c,v}$ plays an
essential role in the calculation of the expansion of $k(c,v)$ ($c\in Z(\ov\mu)$).
\begin{lem}
\label{3.20'}
For any $c\in Z(\ov\mu)\cup\{0\}$, $v\in V$, $\eta\in \B^*_{\tel}$
$$
P_{c^*v}(\wh\eta\circ c) = \big(\wh{ T^*_{c,v}\eta}\big)\circ c.
$$
\end{lem}
\begin{proof}
We observe that, since $c\in Z(\ov\mu)\cup \{0\}$
$$T_{c,v}(\chi_{cx})(y) = \int_H \chi_{cb}(y+v)\chi_{gcx}(y+v)\mu(dh)
 = \int_H \chi_{c(gx+b)}(y+v) \mu(dh)$$
On the other hand
\begin{eqnarray*}
P_{c^*v}(\wh \eta\circ c)(x)&=& \int_H\int_V \chi_{c^* v}(gx+b)\chi_y(c(gx+b))\eta(dy)\mu(dh)\\
&=& \int_H\int_V \chi_{y+v}(c(gx+b))\eta(dy)\mu(dh)\\
&=& \eta\big( T_{c,v} (\chi_{cx})\big)
= \big(\wh{ T^*_{c,v}\eta}\big) (cx).
\end{eqnarray*}
The lemma follows.
\end{proof}
\begin{lem}
\label{3.21}
Suppose $\e<1/2$. There exists $t_3$ such that for $|c|\le t_3$, $c\in Z(\ov\mu)\cup\{0\}$ the function
$$\psi_{c,v} = \wh{\eta}_{c,v}\circ c$$
is a nonzero element of $\B_{\tel}$ and is an eigenfunction of $P_{c,v}$
corresponding to the eigenvalue $k(c,v)$, i.e.
\begin{equation}
\label{ptv}
P_{c,v} (\psi_{c,v}) = k(c,v)\psi_{c,v}.
\end{equation}
Moreover $k'(c,v) = k(c,v)$ and
$$
(k(c,v)-1)\nu(\psi_{c,v}) = \nu(\psi_{c,v}(\chi_{c^* v}-1)).
$$
\end{lem}
\begin{proof}
First we shall prove that $\psi_{c,v}$ is an element of $\B_{\tel}$.
In view of Proposition \ref{lis23} and \eqref{p15} we have
$$ |\psi_{c,v}|_\t = \sup_x
{\frac{\is{\eta_{c,v}}{\chi_{cx}}}{(1+\tau(x))^\t}}\le \|\eta_{c,v}\|\cdot \sup_x
\frac{\|\chi_{cx}\|_{\tel}}{(1+\tau(x))^\t}  \le
(1+C|t_2|^{\eps})(1+2|t_2|^{\eps}).
 $$ To estimate
$[\psi_{c,v}]_{\e,\l}$ we define the function $$ g_{c,x,x'}(y) =
\frac{\chi_{cx}(y)-\chi_{cx'}(y)}{\tau(x-x')^\e(1+\tau(x))^\l(1+\tau(x'))^\l}.
$$ Then
we have
 $$\sup_{x\not= x'}|g_{c,x,x'}|_\t \le 2\cdot
\sup_y\frac{|c|^\e \tau(x-x')^\e \tau(y)^\e} {\tau(x-x')^\e(1+\tau(x))^\l(1+\tau(x'))^\l(1+\tau(y))^\t} \le 2|c|^\e. $$ Next
we have, since $2\eps\le \{1,\l\}$
\begin{eqnarray*}
\sup_{x\not= x'}[g_{c,x,x'}]_{\e,\l}
&=& \sup_{x\not= x'} \sup_{y\not= y'}\frac{| \chi_{cx}(y)-\chi_{cx'}(y)
-\chi_{cx}(y')+\chi_{cx'}(y')|}
{ \tau(x-x')^{\e}\tau(y-y')^{\e}(1+\tau(y))^\l(1+\tau(y'))^\l(1+\tau(x))^\l(1+\tau(x'))^\l }\\
&\le&  \sup_{\tau(y-y')\le \tau(x-x')}\frac{|\chi_{cx}(y-y') - 1 |+| \chi_{cx'}(y-y') -1|
}{ \tau(y-y')^{2\e}(1+\tau(y))^\l(1+\tau(y'))^\l(1+\tau(x))^\l(1+\tau(x'))^\l }\\
&\le&  \sup_{\tau(y-y')\le \tau(x-x')}\frac{2(|c|^{2\e}\tau(x)^{2\e}+\tau(x')^{2\e})
\tau(y-y')^{2\e}}
{ \tau(y-y')^{2\e}(1+\tau(y))^\l(1+\tau(y'))^\l(1+\tau(x))^\l(1+\tau(x'))^\l }\\
&\le& 2 |c|^{2\e},
\end{eqnarray*}
that proves $\sup_{x\not= x'}\|g_{c,x',x}\|_{\th,\e,\l}\le 4|c|^{2\eps}$.
Finally
$$[\psitv]_{\e,\l} = \sup_{x\not= x'}\is{\eta_{c,v}}{g_{c,x,x'}} \le \|\eta_{c,v}\|_{\tel}
\sup_{x\not=x'}\|g_{c,x,x'}\|_{\tel}\le 4|c|^{\eps}(1+C|c|^{\eps})<\8,
$$
and we obtain  $\psi_{c,v}\in \B_{\tel}$. Next notice, using Lemma \ref{3.20'}
$$P_{c,v}(\psi_{c,v})(x)= (T^*_{c,v}\eta_{c,v})(\chi_{cx})
= k'(c,v)\wh{\eta}_{c,v}(cx) = k'(c,v)\psi_{c,v}(x),
$$
but for $|c|$ small enough there exists only one eigenvalue of $P_{c,v}$ close to 1, hence
 $k'(c,v)=k(c,v)$. The last relation is a direct consequence
 of $P_{c,v}\psi_{c,v}=k(c,v)\psi_{c,v}$ and the form of $P_{c,v}$.
\end{proof}
\section{Some technical lemmas}
\label{section3}
\subsection{Some further properties of the stationary measure $\nu$} In the next section very often will appear
expressions of the form $\int_{V}f(c,x)\nu(dx)$, $c\in Z(\ov\mu)\cup\{0\}$ and
we will be interested in their behavior for small values of $|c|$.
We denote $Z_1(\ov\mu)=\{c\in Z(\ov\mu)\cup\{0\}; |c|\le 1 \}$.
We will need the following:
\begin{lem}
\label{p26}
Let $f$ be any continuous function on $Z_1(\ov\mu)\times V$ satisfying
\begin{equation}
\begin{split}
\label{p261}
|f(c,x)| &\le D_{\d,\b} |c|^{\d+\b} \tau(x)^\b,\quad\mbox{ for } \tau(cx)\ge 1,\\
|f(c,x)| &\le D_{\d,\g} |c|^{\d+\g} \tau(x)^\g,\quad\mbox{ for } \tau(cx)\le 1,
\end{split}
\end{equation}
where $\b<\a$, $\g+\d>\a$ and $\d> 0$.
Then
$$ \lim_{c\to 0, c\in Z(\ov\mu)} \frac 1{|c|^\a}\int_{V} f(c,x)\nu(dx) = 0.
$$
\end{lem}
\begin{proof}
Notice that taking the function $f(x)={\bf 1}_{\{\tau(x)\ge 1\}}$, by \eqref{has},
there   exists $C>0$ such that
\begin{equation}
\label{has2} \nu(\tau(x)>t) \le D t^{-\a},
\end{equation} for any $t>0$.
  We divide
the integral into three parts and study each of them independently.
For appropriately small values of $|c|$ we have
\begin{eqnarray*}
 \frac 1{|c|^\a}\int_{\tau(x)\le 1}f(c,x)\nu(dx)
&\le& \frac {D_{\d,\g}}{|c|^\a}\int_{\tau(x)\le 1}|c|^{\g+\d}\tau(x)^\g\nu(dx)  \le
D_{\d,\g}\nu\big[ \tau(x)\le 1 \big]{|c|^{\g+\d-\a }}\\ \frac
1{|c|^\a}\int_{1 < \tau(x)\le \frac 1{|c|}} f(c,x)\nu(dx) &\le&
D_{\d,\g}{|c|^{\g+\d-\a}}\int_{1 < \tau(x)\le \frac 1{|c|}}  \tau(x)^\g \nu(dx)\\ &\le&
D_{\d,\g}{|c|^{\g+\d-\a}} \sum_{n=0}^{\suf{-\log_2 |c|}}   \int_{\frac
1{2^{n+1}|c|} < \tau(x)\le \frac 1{2^{n}|c|}}  \tau(x)^\g \nu(dx)\\ &\le&
D_{\d,\g}{|c|^{\g+\d-\a}} \sum_{n=0}^{\suf{-\log_2 |c|}}  \frac 1{(2^n |c|)^\g}
\nu\Big[\frac 1{2^{n+1}|c|} < \tau(x) \Big]\\ &\le& D_{\d,\g}{|c|^{\g+\d-\a}}
\!\!\!\sum_{n=0}^{\suf{-\log_2 |c|}}  \frac 1{(2^n |c|)^\g} (2^n |c|)^\a \le D_{\d,\g}
|c|^\d\!\!\! \sum_{n=0}^{\suf{-\log_2 |c|}} 2^{n(\a-\g)}\\
\end{eqnarray*}
If $\a\le\g$ then the foregoing sum can be estimated by some constant or  $\big|\log {|c|}\big|$
and the expression converges to zero. On the other hand if $\a>\g$,
the sum is smaller than $D_{\d,\g}|c|^{\g-\a}$ and multiplied by $|c|^\d$ converges
 to zero.
Finally
\begin{multline*}
\frac 1{|c|^\a}\int_{\tau(x) >  \frac 1{|c|}} f(c,x)\nu(dx)
\le D_{\d,\b}{|c|^{\b+\d-\a}} \sum_{n=0}^{\8}   \int_{\frac {2^{n}}{|c|} < \tau(x)\le \frac {2^{n+1}}{|c|}}  \tau(x)^\b \nu(dx)\\
\le D_{\d,\b}{|c|^{\b+\d-\a}} \sum_{n=0}^{\8}  \Big(\frac {2^{n+1} }{|c|}\Big)^\b \nu\Big[
\frac {2^n}{|c|} < \tau(x) \Big]
\le D_{\d,\b}{|c|^{\d-\a}} \sum_{n=0}^{\8} 2^{n\b}\cdot \frac{|c|^\a}{2^{n\a}} \le D_{\d,\b} |c|^\d
\end{multline*}
Of course the same proof gives the second part of the Lemma.
\end{proof}
\subsection{Some properties of the eigenfunction  $\psi_{c,v}$}
Up to now, we have not taken any care about precise values of parameters $\th$, $\e$, $\l$. However,
we will need some further hypotheses, and from now on, we will assume additionally that
\begin{equation}
\label{h6}
\begin{split}
& \mbox{ if $1<\a<2$, then $1+\l+\e>\a$,}\\
& \mbox{ if $\a=2$, then $\l+2\e>1$.}\\
& \mbox{ if $\a>2$, then $\l=1$.}
\end{split}
\end{equation}
It can be easily proved that there exist $\th$, $\e$, $\l$ satisfying all the assumptions
of Propositions \ref{kel}, \ref{lis22} and the condition above.
\begin{lem}
\label{p21}
There exists $D''$ such that
\begin{eqnarray*}
|\psi_{c,v}(x) - \wh{\eta}_v(cx)| &\le& D'' |c|^{2\e}\tau(x)^{ \e},
\quad\mbox{for }\tau(cx)>1,\\
|\psi_{c,v}(x) - \wh{\eta}_v(cx)| &\le& D'' |c|^\e \tau(cx)^{\eta }, \quad\mbox{for } \tau(cx)\le 1,
\end{eqnarray*}
for $\eta =\min \{1,\l + \e\}.$
\end{lem}
\begin{proof}
Let us first estimate the norm $\|\chi_{cx}-1\|_{\tel}$. Let $0<\b <1$.
  We have, since $\eta\le 1$
\begin{eqnarray*}
|\chi_{cx}-1|_\t &\le&  2 \tau(cx)^{\eta }\\ \ [\chi_{cx}-1]_{\e,\l}
&\le& \sup_{y,z}\frac{2 \tau(cx)^{\eta}
\tau(y-z)^{\eta}}{\tau(y-z)^\e (1+\tau(y))^\l(1+\tau(z))^\l} \le 2
\tau(cx)^{\eta}
\end{eqnarray*}
which proves
\begin{equation}
\label{paz1}
\|\chi_{cx}-1\|_{\tel} \le 4 \tau(cx)^{\eta}.
\end{equation}
Therefore, by Proposition \ref{lis23}, and the expression of $\psi_{c,v}$
given by Lemma \ref{3.21}
\begin{eqnarray*}
|\psi_{c,v}(x) - \wh{\eta}_v(cx)| &\le& |\is{\eta_{c,v}-\eta_v}{\chi_{cx}-1}|\le C' |c|^{\e}\| \chi_{cx}-1 \|_{\tel}\le
4C' |c|^{\e}\tau(cx)^{\eta }.
\end{eqnarray*}
For  $\tau(cx)>1$ we need  better estimates. By  Proposition \ref{lis23} and \eqref{p15} we have
$$|\psi_{c,v}(x) - \wh{\eta}_v(cx)| = |\is{\eta_{c,v}-\eta_{v}}{\chi_{cx}}| \le C'|c|^{\e}\|\chi_{cx}\|_{\tel}
\le C' |c|^{\e}\big(1+ 2\tau(cx)^{\e}\big)\le 3C'|c|^{\eps}\tau(cx)^{\eps},
$$ if $\tau(cx) \ge 1$. The result follows with $D''=4C'$.
\end{proof}
\begin{cor}
\label{p35}
If  $\a \le 2$, then
$$
\lim_{c\to 0, c\in Z(\ov\mu)}\frac 1{|c|^\a}\int_{V}\big( \chi_v(cx)-1\big)\big(\psi_{c,v}(x)-\wh{\eta}_v(cx)\big)\nu(dx) = 0.
$$
\end{cor}
\begin{proof} We will apply Lemma \ref{p26} for
$ f(c,x) = \big(
\chi_v(cx)-1\big)\big(\psi_{c,v}(x)-\wh{\eta}_v(cx)\big). $ Let's check
that its hypotheses are satisfied. For $\tau(cx)\ge 1$, by Lemma \ref{p21},
we have $$ |f(c,x)| \le D |c|^{2\e} \tau(x)^{\e }. $$ Next for
$\tau(cx)\le 1$ and $\eta $ as above $$ |f(c,x)| \le D |c|^\e\cdot
\tau(cx)^{1+\eta  }$$ and in view of \eqref{h6}, $1+\eta +\e > \a$, so the hypotheses of Lemma \ref{p26} are fulfilled.
\end{proof}
\subsection{Eigenvalue $k(c,v)$ of $P_{c,v}$}
From Lemma \ref{3.21}
we obtain, if $c\in Z(\ov\mu)\cup \{ 0 \}$
\begin{equation}
\label{kt}
(k(c,v)-1)\nu(\psi_{c,v}) = \nu(\psi_{c,v}(\chi_{c^*v}-1)),
 \end{equation}
for $\psi_{c,v}=\wh\eta_{c,v}\circ c$
The formula  will be crucial in sections \ref{k-euclid} and \ref{k-gen} to describe asymptotic behavior of
the function $c\mapsto k(c,v)$ near zero.
One can easily prove that $\nu(\psi_{c,v})$ goes to 1, hence to understand behavior of $k(c,v)$
near 1 one has to describe the integral above for small $|c|$. For some technical reasons
we will need also speed of convergence of $\nu(\psi_{c,v})$ to 1.

\begin{lem}
\label{p361}
Assume $v$ is fixed. Then there exists $D'''>0$ and $t_3>0$ such that
for $|c|<t_3$, $c\in Z(\ov\mu)\cup\{0\}$, we have
$$|1- \is{\nu}{\psi_{c,v}}| \le D'''|c|^{\min\{1,\l+\e\}}.$$
\end{lem}
\begin{proof}
We use the formula $\psi_{c,v}(x)=\eta_{c,v}(\chi_{cx})$ for $|c|\le t_3$.
Then by Lemma \ref{3.21} and formulae \eqref{paz1} and \eqref{moments} we have, with $C'=\sup_{|c|\le t_3}\|\eta_{c,v}\|_{\tel}$
\begin{eqnarray*}
\big| 1-\nu(\psi_{c,v}) \big| &\le&  \int_{V} \big|\is{\eta_{c,v}}{1-\chi_{cx}}\big|\nu(dx)\le
D' \int_{V} \|{1-\chi_{cx}}\|_{\tel}\nu(dx)\\
&\le& D'|c|^{\min\{1,\l+\e\}}\int_{V} \tau(x)^{\min\{1,\l+\e\}} \nu(dx) \le D'''|c|^{\min\{1,\l+\e\}},
\end{eqnarray*}
with $D'''=D'\int_V \tau(x)^{\min\{1,\l+\e\}}\nu(dx)$.
\end{proof}

\section{Asymptotic expansions of  eigenvalues $k(c,v)$ in the Euclidean case}
\label{k-euclid}
The purpose of this section is give asymptotic expansions
 of the eigenvalues
$k(c,v)$ when $|c|$ goes to 0. First to present main ideas of the proof we will
consider the Euclidean case, then $|x|^2=\sum x_i^2$.
If $\a\le 2$, we take $c\in Z_{\ov\mu}$.
If $\a>2$ we take $c=t\in \R^*_+$.
The main result of this section is the following

\begin{thm}
\label{asymptkt}
\begin{enumerate}
\item If $0<\a<1$ then
$$
\lim_{c\to 0, c\in Z_{\ov\mu}}\frac{k(c,v)-1}{|c|^\a} =  C_\a(v)
$$ where
$$ C_\a(v) = \int_{V} \big( \chi_v(x)-1 \big)\wh{\eta}_v(x)\L(dx)
$$
and $C_\a(c^*v)=|c|^\a C_\a(v)$ for any $c\in Z_{\ov\mu}$.
\item  If $\a=1$ then
$$\lim_{c\to 0, c\in Z_{\ov\mu}}\frac{k(c,v)-1-i\is v{\xi(c)}}{|c|} = C_1(v),$$
for
$$
C_1(v) = \int_{V}\bigg( \big(\chi_v(x)-1\big) \wh{\eta}_v(x)-\frac{i\is vx}{1+|x|^2}\bigg) \L(dx),
$$ and $\xi(c)=\int_{V}\frac{cx}{1+|c|^2|x|^2}\nu(dx)$.
Moreover $|\xi(c)|\le I(v) |c||\log |c||$ for $|c|<1/2$ and
$|\xi(c)|\le I(v) |c|$ for $|c|>1/2$. Furthermore $C_1(c^* v)=|c| C_1(v)+ i\is v{\b(c)}$ with
$\b(c) = \int_V \big( \frac x{1+|x|^2} - \frac x{1+|cx|^2}  \big)\Lambda(dx)$.
\item
If $1<\a<2$
$$
\lim_{c\to 0, c\in Z_{\ov\mu}} \frac {k(c,v)-1 - i\is v{cm}  }{|c|^\a} = C_\a(v),
$$ where
$$ C_\a(v)=\int_{V}\Big( \big(\chi_v(x)-1\big) \wh{\eta}_v(x)-i\is vx \Big) \L(dx) $$
and $C_\a(c^* v)=|c|^\a C_\a(v)$ for any $c\in Z_{\ov\mu}$.
\item If $\a=2$
$$
\lim_{c\to 0, c\in Z_{\ov \mu}} \frac{k(c,v)-1-i\is v{cm}}{|c|^2|\log |c||} = 2C_2(v),
$$
where
$$
C_2(v) = -\frac 14\int_{\Sigma_1}\bigg( \is vw^2 + 2\is vw\eta_v(w^*)\bigg)\s(dw)
$$
and  $C_2(c^*v)=|c|^2 C_2(v)$ for every $c\in Z_{\ov\mu}$.
\item
If $\a>2$ then
$$
\lim_{t\to 0^+}\frac{k(t,v)-1-i\is v{tm}}{t^2} = C_{2+}(v),
$$
$$
C_{2+}(v) = -\frac 12 q(v,v)  -\frac 12 \is vm^2 - q\big(v,(I-z^*)^{-1}z^*v\big).
$$
\end{enumerate}
\end{thm}

To prove the Theorem we shall consider each case separately.
\subsection{Case: $\a<1$}
Let us write
\begin{multline*}
 \frac 1{|c|^\a}\int_{V} \big( \chi_v(cx)-1 \big)\psi_{c,v}(x)\nu(dx) \\
= \frac 1{|c|^\a}\int_{V} \big( \chi_v(cx)-1 \big)\cdot\big(\psi_{c,v}(x)-\wh{\eta}_v(cx)\big)\nu(dx)
 +\frac 1{|c|^\a}\int_{V} \big(\chi_v(cx) - 1 \big)\wh{\eta}_v(cx)\nu(dx)
\end{multline*}
and notice that by Corollary \ref{p35} the first term of the sum above tends to zero. To
describe the second one, observe that the function $f_v=(\chi_v-1)\wh{\eta}_v$ satisfies
\eqref{beta}. In fact the characteristic  function $\wh{\eta}_v$ is bounded by 1, hence also $f_v$ is bounded, and
for $|x|<1$, we have $|f_v(x)|\le 2|x|$. Therefore by \eqref{has} the expression above tends to
the constant $C_\a(v)$. Finally by \eqref{kt} and Lemma \ref{p361}
\begin{equation}
\label{case1}
\begin{split}
\lim_{c\to 0}\frac{k(c,v)-1}{|c|^\a} &= \lim_{c\to 0}\frac{1}{\nu(\psi_{c,v})|c|^\a}\cdot
\int_{V} \Big(\chi_v(cx)-1\Big)\psi_{c,v}(x) \nu(dx) \\
&= \lim_{c\to 0}\frac{1}{\nu(\psi_{c,v})|c|^\a}\cdot
\int_{V} \big(\chi_v(cx)-1\big)\wh{\eta}_v(cx) \nu(dx) = C_\a(v)
\end{split}
\end{equation}
as $c$ goes to 0. The last assertion is an immediate consequence of the homogeneity
of $\eta_v$ given by Lemma \ref{3.6'} and of the homogeneity of $\L$
mentioned in Section 2.
\subsection{Case:  $\a=1$}
\begin{lem}
\label{paz2}
$$\lim_{c\to 0}\frac 1{|c|}\Big( \is{\nu}{(\chi_{c^*v}-1)\psi_{c,v}} - i\is{v}{\xi(c)} \Big)=C_1(v).$$
\end{lem}
\begin{proof}
We have
\begin{equation}
\label{sob1}
\begin{split}
\is{\nu}{(\chi_{c^*v}-1)\psi_{c,v}} &=
\int_{V}\big(  \chi_v(cx)-1\big)\psi_{c,v}(x)\nu(dx) \\
&=\int_{V}\big( \chi_v(cx)-1\big)\cdot \big(\psi_{c,v}(x)-\wh{\eta}_v(cv)\big)\nu(dx)\\
&+ \int_{V}\big( \chi_v(cx) -1\big)\big(\wh{\eta}_v(cv)-1\big)\nu(dx)\\
&+\int_{V}\bigg( \chi_v(cx)-1 - \frac{i\is{v}{cx}}{1+|cx|^2}\bigg)\nu(dx)
+ i\is{v}{\xi(c)} \\&= \sum_{j=0}^2 W_i(c) + i\is v{\xi(c)}
\end{split}
\end{equation}
By Corollary \ref{p35}, $\frac{W_0(c)}{|c|}$ converges to 0, as $c$ goes to 0. Next observe that the function
$ f_1 = \big( \chi_v-1\big)\big(\wh{\eta}_v-1\big)
$ satisfies \eqref{beta}. Indeed $f_1$ is bounded and for $|x|\le 1$, from \eqref{paz1}
$$ |f_1(x)|\le \big|\chi_v(x)-1 \big| \big|  \wh{\eta}_v(x)-1 \big| \le
2|v||x||\is{\eta_v}{\chi_{x}-1}|\le 2|v||x|\|\chi_{x}-1\|_{\tel}\le 8|v||x|^{1+\l+\e}
$$
Similarly one can prove that
$  f_2(x) =\chi_v(x)  - 1 -\frac{i\is vx}{1+|x|^2}$ fulfills \eqref{beta}.
Thus, by \ref{has}
$$\lim_{c\to 0}\bigg(\frac{W_1(c)}{|c|} +\frac{W_2(c)}{|c|}\bigg)
=\int_v\big( f_1(x)+f_2(x) \big) \L(dx) = C_1(v),$$
which finishes the proof.
\end{proof}
\begin{lem}
\label{ksit}
There exists a constant $I(v)$ such that: $$|\xi(c)|\le
\left\{
\begin{array}{cc}
I(v) |c|, & \mbox{for } |c|\ge \frac 12,\\
I(v) |c||\log |c||, & \mbox{for } |c|< \frac 12.
\end{array}
\right.$$
\end{lem}
\begin{proof}
For $|c|\ge 1$ the Lemma is obvious. For $|c|<1$ we write
$$
|\xi(c)|\le \int_{V}\frac{|cx|}{1+|c|^2|x|^2}\nu(dx) = \Bigg(\int_{|x|\le 1} + \int_{1< |x|\le \frac 1{|c|}}+
\int_{\frac 1{|c|} < |x|}\Bigg)\bigg( \frac{|cx|}{1+|c|^2|x|^2}\bigg)\nu(dx).
$$
The first integral is dominated by $|c|$. In view of \eqref{has}
 the third one, divided by $|c|$, converges
 to $\int_{|x|>1}\frac{|x|}{1+|x|^2}\L(dx)$, as $s$ goes to 0.
Finally applying \eqref{has2} we estimate the second integral by
\begin{multline*}
\sum_{k=0}^{|\log |c||}\int_{2^k<|x|\le 2^{k+1}}\frac{|cx|}{1+|c|^2|x|^2}\nu(dx)\le
|c| \sum_{k=0}^{|\log |c||}2^{k+1}\nu[2^k<|x|]\\ \le C|c| \sum_{k=0}^{|\log |c||}2^{k+1}\cdot 2^{-k}\le C |c||\log |c||.
\end{multline*}
\end{proof}
\begin{proof}[Proof of Theorem \ref{asymptkt}, part (2)]
By \eqref{kt} and  Lemmas \ref{p361}, \ref{paz2}, \ref{ksit} we have
\begin{eqnarray*}
\lim_{c\to 0}\frac{k(c,v)-1-i\is v{\xi(c)}}{|c|}
&=& \lim_{c\to 0 }\bigg[\frac{\is{\nu}{(\chi_{cx}-1)\psi_{c,v}}-i\is v{\xi(c)}}{\nu(\psi_{c,v})|c|} +
\frac{i\big(1-\nu(\psi_{c,v})\big)\is v{\xi(c)}}{\nu(\psi_{c,v})|c|}\bigg]\\
&=& C_1(v)
\end{eqnarray*}
\end{proof}

\subsection{Case: $1<\a<2$}
\begin{lem}
\label{u113}
$$
\lim_{c\to 0, c\in Z_{\ov\mu}} \frac 1{|c|^\a}\Big(
\is{\nu}{(\chi_{c^*v}-1)\psi_{c,v}} - i\is v{cm} \Big) = C_\a(v).
$$
\end{lem}
\begin{proof}
In view of Corollary \ref{p35} it is enough to consider
\begin{multline*}
\int_{V}\big(\chi_v(cx) -1\big)\wh{\eta}_v(cx)\nu(dx)
=\int_{V}\big(\chi_v(cx)-1\big)\cdot \big(\wh{\eta}_v(cx) -1\big)\nu(dx)\\
+\int_{V}\big(\chi_v(cx)-1- i\is{v}{cx}\big)\nu(dx)
+i\is v{cm}.
\end{multline*}
Reasoning as in previous cases we prove that the functions $f_1 =
\big( \chi_v-1\big)\cdot \big(\wh{\eta}_v -1\big)$ and $f_2= \chi_v-1- i v^*$
satisfy \eqref{beta}, therefore \eqref{has} implies the Lemma.
\end{proof}
\begin{proof}[Proof of Theorem \ref{asymptkt}, part (3)]
By \eqref{kt}, Lemma \ref{p361} and \eqref{h6}
\begin{eqnarray*}
\lim_{c\to 0}\frac{k(c,v)-1 - i \is v{cm} }{ |c|^\a}
&=& \lim_{c\to 0}\bigg[\frac{\nu(\psi_{c,v})(k(c,v)-1) - i \is v{cm}
}{\nu(\psi_{c,v}) |c|^\a}+ \frac{i \is v{cm} (1- \nu(\psi_{c,v}))}{
\nu(\psi_{c,v})|c|^\a} \bigg] \\
&=&  C_\a(v)
\end{eqnarray*}
\end{proof}
\subsection{Case: $\a=2$}
\begin{lem}
\label{foverg} Suppose we are given two functions on $V$,  $f$ and $h$ such that
 $h(x)=\is x{v_1}\is x{v_2}$ for some $v_1,v_2\in V$,  $\lim_{x\to 0}\frac{f(x)}{h(x)}=C_0$
and $|f(x)|\le C|x|^{1+\eta}$ for some positive constants $C_0$, $C$ and $\eta<1$. Then
$$
\lim_{g\to 0, g\in G_{\ov\mu}} \frac 1{|g|^2|\log|g||} \int_V f(gx)\nu(dx) = C_0\int_{\Sigma_1}h(w)\s(dw),
$$
where $\s$ is the measure on the fundamental domain $\Sigma_1$ defined in \eqref{radial}.

Moreover the function $\wt \Lambda(v)=\int_{\Sigma_1}\is vw^2\sigma(dw)$ is $G^*_{\ov\mu}$-homogeneous,
i.e. $\wt \Lambda(g^* v) = |g|^2\wt \Lambda (v)$ for every $g\in G_{\ov\mu}$.
\end{lem}
\begin{proof}
Fix $\b\in R_{\ov\mu}$ such that $\b>1$ and denote by $U$ the annulus $U=\{x\in V: 1<|x|\le \b\}$.
 Next we fix arbitrary small number $\d>0$. Then there exists $\e$ such that
\begin{equation}
\label{er3}
\bigg| \frac{f(x)}{h(x)} - C_0\bigg|<\d,\qquad\mbox{for }|x|<\e.
\end{equation}
Without any lose of generality we may assume $|v_1|=|v_2|=1$. Given $y_1, y_2\in V$ we define a function
on $V$, $h_{y_1,y_2}(x)=\is x{y_1}\is x{y_2}$. We are going to prove that there exists large $A\in R_{\ov\mu}$
such that
\begin{equation}
\label{er4}
\begin{split}
 \bigg| \frac 1{|g|^2} \int_V {\bf 1}_U(gx)   h_{y_1,y_2}(gx)\nu(dx) - \int_U h_{y_1,y_2}(x)\L(dx)\bigg| &\le\d,\\
 \bigg| \frac 1{|g|^2} \int_V {\bf 1}_U(gx) \big|h_{y_1,y_2}(gx)\big|\nu(dx) - \int_U \big|h_{y_1,y_2}(x)\big|\L(dx)\bigg| &\le\d,\\
\end{split}
\end{equation}
for any $g\in G_{\ov\mu}$ such that $|g|<\frac 1A$ and all $y_1,y_2$ belonging to $S_1$, the unit sphere
in $V$.

Of course the last assertion, by \eqref{has}, is clear for fixed vectors $y_1$ and $y_2$. However,
we will justify that also uniform estimates are valid.
\medskip

Fix $y_1,y_2\in S_1$. Then in view of \eqref{has} there exists $M\in R_{\ov\mu}$ such that
$$
 \bigg| \frac 1{|g|^2} \int_V {\bf 1}_U(gx)   h_{y_1,y_2}(gx)\nu(dx) - \int_U h_{y_1,y_2}(x)\L(dx)\bigg| \le\frac \d2,
$$ and
$$
 \bigg| \frac 1{|g|^2} \int_V {\bf 1}_U(gx) |gx|^2 \nu(dx) - \int_U |x|^2\L(dx)\bigg| \le 1,
$$ for $|g|\le \frac 1M$. Choose $\d'<\frac \d {4 ( 2\b^2\L(U)+1)}$. Define
$B_{y_1,y_2}(\d')$ to be the ball in $V\times V$ centered at $(y_1,y_2)$ of radius $\d'$ and take $(y_1',y_2')\in B_{y_1,y_2}(\d')$.
Notice
$$\big| h_{y_1',y_2'}(x) - h_{y_1,y_2}(x)\big| \le \big|\is{y_1-y_1'}x\is{y_2}x\big| +
\big| \is {y_1'}x\is{y_2-y_2'}x \big| \le 2\d'|x|^2,$$
therefore
\begin{multline*}
 \bigg| \frac 1{|g|^2} \int_V{\bf 1}_U(gx)   h_{y_1',y_2'}(gx)\nu(dx) - \int_U h_{y_1',y_2'}(x)\L(dx)\bigg|\\
 \le
  \bigg| \frac 1{|g|^2} \int_V {\bf 1}_U(gx) \Big(   h_{y_1',y_2'}(gx) -  h_{y_1,y_2}(gx)\Big)\nu(dx)\bigg|
+  \bigg|  \int_U\Big(   h_{y_1',y_2'}(x) -  h_{y_1,y_2}(x)\Big)\L(dx)\bigg| +\frac \d 2\\
\le \frac {2\d'}{|g|^2}\int_V {\bf 1}_U(gx) |gx|^2\nu(dx) + 2\d' \int_U |x|^2\L(dx) + \frac \d2\\
\le 2\d' \big( 2\b^2\L(U)+1 \big) + \frac \d 2 < \d.
\end{multline*}
So, we may find finitely many pairs $\{(y_{i,1},y_{i,2})\}_{1\le i\le N}$ and positive numbers $M_i\in R_{\ov\mu}$
such that the balls $B_{y_{i,1},y_{i,2}}(\d')$ cover $S_1$. Then choosing $A_1=\max_{1\le i \le N} M_i$ we deduce
that the first line of \eqref{er4} is satisfied for $|g|<\frac 1{A_1}$. Next we repeat our argument for $|h_{y_1,y_2}|$
instead of  $h_{y_1,y_2}$, we find $A_2$ and finally choosing $A=\max\{A_1,A_2\}$ we obtain \eqref{er4}.
\medskip

For $|g|<\frac \e A$ (that will be assumed from now), we divide the integral of $f$
into three parts:
\begin{equation}
\label{s1}
\int_V f(gx)\nu(dx) = \int_{|x|\le A} f(gx)\nu(dx) + \int_{A<|x|<\frac \e{|g|}} f(gx)\nu(dx)+\int_{|x|\ge \frac \e{|g|}} f(gx)\nu(dx).
\end{equation}
Notice first that
$$
\bigg|\int_{|x|\le A} f(gx)\nu(dx) \bigg| \le (|C_0|+\d)\int_{|x|\le A}|h(gx)|\nu(dx)
\le (|C_0|+\d)|g|^2A^2
$$
and by \eqref{has}
$$
\lim_{g\to 0,g\in G_{\ov\mu}}\frac 1{|g|^2} \int_{|x|\ge \frac \e{|g|}} f(gx)\nu(dx) = \int_{|x|\ge \e}f(x)\L(dx).
$$
Hence
$$
\lim_{g\to 0,g\in G_{\ov\mu}}\frac 1{|g|^2|\log|g||} \bigg( \int_{|x| \le A} f(gx)\nu(dx)+
\int_{|x|\ge \frac \e{|g|}} f(gx)\nu(dx)\bigg)=0.
$$
Therefore we have to handle with the middle term in  \eqref{s1}. We will prove
\begin{equation}
\label{er2}
\lim_{g\to 0,g\in G_{\ov\mu}}\frac 1{|g|^2|\log|g||}  \int_{A<|x|<\frac\e{|g|} } f(gx)\nu(dx)=
\frac{C_0}{\log\b} \int_U h(x)\L(dx).
\end{equation}
Applying \eqref{er3}, we write
\begin{multline}
\label{er12}
\bigg| \frac 1{|g|^2|\log|g||}  \int_{A<|x|<\frac\e{|g|} } f(gx)\nu(dx)
- \frac{C_0}{\log\b} \int_U h(x)\L(dx)
\bigg|\\
\le
|C_0|\cdot
\bigg| \frac 1{|g|^2|\log|g||} \int_{A<|x|<\frac\e{|g|} } h(gx)\nu(dx)
- \frac{1}{\log\b} \int_U h(x)\L(dx)
\bigg|\\ + \frac \d{|g|^2|\log|g||} \int_{A<|x|<\frac\e{|g|} }\big| h(gx)\big|\nu(dx).
\end{multline}
We estimate the first expression. For this purpose we define $K=\big\lfloor \log\frac{\e}{A|g|}/\log\b \big\rfloor-1$.
For  $ r \in R_{\ov\mu}$ we will denote by $g(r)$ any element of $G_{\ov\mu}$ such that $|g(r)|=r$.
To simplify our notation we define elements of $G_{\ov\mu}$:
$g_n=g(A\b^n)$ and annulus $U_n=\{x:\; A\b^n<|x|\le A\b^{n+1}\}$.
Notice that $|g_n|> A$, therefore applying \eqref{er4} and $G_{\ov\mu}$ homogeneity of $\L$ we obtain
\begin{multline*}
\bigg| \frac 1{|g|^2|\log|g||} \int_{A<|x|<\frac\e{|g|} } h(gx)\nu(dx)
- \frac{1}{\log\b} \int_U h(x)\L(dx)
\bigg|\\
\le
\bigg| \frac 1{|g|^2|\log|g||} \sum_{n=0}^K\int_{U_n } h(gx)\nu(dx)
- \frac{1}{\log\b} \int_U h(x)\L(dx)
\bigg| + \frac 1{|g|^2|\log|g||}\int_{U_{K+1} }\big| h(gx)\big|\nu(dx)\\
=
\bigg| \frac 1{|\log|g||} \sum_{n=0}^K|g_n|^2 \int_V {\bf 1}_U(g_n^{-1}x) h_{\frac{(gg_n)^*v_1}{|gg_n|},\frac{(gg_n)^*v_2}{|gg_n|}}
(g_n^{-1}x)\nu(dx)
- \frac{1}{\log\b} \int_{U} h(x)\L(dx)\bigg| \\
+ \frac {|g_{K+1}|^2}{|\log|g||}\int_V {\bf 1}_U (g_{K+1}^{-1}x)\Big| h_{\frac{(gg_{K+1})^*v_1}{|gg_{K+1}|},
\frac{(gg_{K+1})^*v_2}{|gg_{K+1}|}}(g_{K+1}^{-1}x)\Big|\nu(dx)\\
\le
\bigg| \frac 1{\log|g||} \sum_{n=0}^K \frac 1{|gg_n|^2} \int_U h(gg_nx)\L(dx)
- \frac{1}{\log\b} \int_U h(x)\L(dx)\bigg|
+ \frac {\d(K+1)}{|\log|g||} \\+ \frac 1{|\log|g||}\cdot \bigg( \frac 1{|g g_{K+1}|^2}
 \int_U\big| h(g g_{K+1}x) \big|\L(dx)+\d  \bigg)\\
= \bigg| \bigg( \frac{K+1}{|\log|g||} - \frac 1{\log\b} \bigg)\cdot \int_U h(x)\L(dx) \bigg|
+\frac{\d(K+1)}{|\log|g||} + \frac 1{|\log|g||}\cdot
\bigg( \int_U \big| h(x) \big|\L(dx)+\d  \bigg)
\end{multline*}
The second term in \eqref{er12} can be estimated using exactly the same arguments. Thus, we obtain
$$\frac \d{|g|^2|\log|g||} \int_{A<|x|<\frac\e{|g|} }\big| h(gx)\big|\nu(dx) \le
 \frac {\d (K+2)}{|\log|g||}\cdot \bigg( \int_U\big| h(x) \big|\L(dx)+\d  \bigg).$$
Therefore passing to the limit in \eqref{er12} we obtain
\begin{multline*}
\limsup_{g\to 0, g\in G_{\ov\mu}}\bigg| \frac 1{|g|^2|\log|g||}  \int_{A<|x|<\frac\e{|g|} } f(gx)\nu(dx)
- \frac{C_0}{\log\b} \int_U h(x)\L(dx)
\bigg|\\
\le \frac \d{\log \b}\bigg( 1+\d + \int_U|h(x)|\L(dx) \bigg),
\end{multline*}
but $\d$ can be arbitrary small, hence we obtain \eqref{er2}.

\medskip

Finally to conclude we choose $\b=p$ if $R_{\ov\mu}=\langle p\rangle $. Otherwise, if $R_{\ov\mu}=\R^*_+$,
we compute  the limit as $\b$ tends to 1. For this purpose, given $a\in A_{\ov\mu}$ and $w\in V$ we will write
$aw=|a|\theta(w)$, where $|\theta(w)|=|w|$. Then $\theta(w)$ tends to $w$, if $|a|$ tends to 1.
By \eqref{radial} we write
\begin{eqnarray*}
\frac 1{\log\b} \int_{1<|x|\le \b} h(x)\L(dx) &=& \frac 1{\log \b} \int_{1<|a|\le \b}\int_{\Sigma_1}h(aw)\sigma(dw)\frac{da}{|a|^3}\\
&=& \frac 1{\log \b} \int_{1<|a|\le \b}\int_{\Sigma_1}\big( h(\th(w)) - h(w)\big) \sigma(dw)\frac{da}{|a|}
+ \int_{\Sigma_1}h(w)\sigma(dw).
\end{eqnarray*}
Hence passing with $\b$ to the limit we obtain
$$
\lim_{\b\to 1} \frac 1{\log\b} \int_{1<|x|\le \b} h(x)\L(dx) = \int_{\Sigma_1} h(w)\sigma(dw).
$$
Combining the formula above and \eqref{er2} we prove the first part of the Lemma

\medskip

To prove the last assertion, assume $v=v_1=v_2$ and  notice that the limit \eqref{er2} does not depend on $\b$, hence if
for $\b\in R_{\ov\mu}$ we define
$$
H_\b(v) = \frac {1}{\log\b} \int_{1<x\le |\b|}\is xv^2\L(dx),
$$ then $H_\b$ in fact does not depend on $\b$ and moreover $H_\b(v)=\wt \Lambda(v)$.
Therefore it is enough to prove that
\begin{equation}
\label{er8}
\lim_{\b\to\8,\b\in R_{\ov\mu}} H_{\b}(g^*v) = |g|^2 \lim_{\b\to\8,\b\in R_{\ov\mu}} H_{\b}(v), \quad\mbox{for }g\in G_{\ov\mu},
\end{equation}
because then
$$
\wt \Lambda(g^* v) = \lim_{\b\to\8,\b\in R_{\ov\mu}} H_{\b}(g^*v) = |g|^2 \lim_{\b\to\8,\b\in R_{\ov\mu}} H_{\b}(v)=
|g|^2\wt \Lambda(v), \quad\mbox{for }g\in G_{\ov\mu}.
$$
Assume $|g|>1$. We apply $G_{\ov\mu}$ homogeneity of $\L$ and write
\begin{eqnarray*}
H_\b(g^* v) &=& \frac {1}{\log \b} \int_{1<|x|\le \b} \is v{gx}^2\L(dx)\\
&=& \frac {|g|^2}{\log \b}  \int_{|g|<|x|\le \b|g|} \is vx^2\L(dx) \\&=& \frac{|g|^2\log\big(\b|g|\big)}{\log\b} H_{\b|g|}(v)
-  \frac{|g|^2\log|g|}{\log\b} H_{|g|}(v).
\end{eqnarray*}
 Passing with $\b$ to infinity we obtain \eqref{er8} and finish the proof of the lemma.
\end{proof}
\begin{lem}
\label{alfa2}
We have
$$
\lim_{c\to 0,c\in Z_{\ov\mu}} \frac 1{|c|^2|\log |c||} \Big( \is{\nu}{(\chi_{c^*v}-1)\psi_{c,v}} - i\is v{cm} \Big)= 2 C_2(v).
$$
\end{lem}
\begin{proof}
We begin as in previous cases and write
\begin{eqnarray*}
\is{\nu}{(\chi_{c^*v}-1)\psi_{c,v}}
&=& \int_V\big(\chi_v (cx)-1\big)\cdot \big( \psi_{c,v}(x) -\wh\eta_v(cx) \big)\nu(dx)\\
&&+ \int_V\big(\chi_v (cx)-1\big)\cdot \eta_v\big( \chi_{cx} -1-i (cx)^* \big)\nu(dx)\\
&&+ \int_V\big(\chi_v (cx)-1\big)\cdot \eta_v\big( i (cx)^* \big)\nu(dx)
+ \int_V\big(\chi_v (cx)-1 -i\is v{cx}\big)\nu(dx)\\
&&+i\is v{cm}.
\end{eqnarray*}
The first term, in view of Corollary \ref{p35}, divided by $|c|^2$ goes to zero.
The second one divided by $|c|^2$, by \eqref{has} has a finite limit.
Hence both divided by $|c|^2|\log |c||$ tend to 0. To handle with the third and the fourth
expression we will use Lemma \ref{foverg}. Notice
\begin{eqnarray*}
\lim_{x\to 0}\frac{(\chi_v(x)-1)\eta_v(ix^*)}{\is xv\is x{m_v}} &=&-1,\\
\lim_{x\to 0}\frac{\chi_v(x)-1-i\is vx}{\is xv^2} &=&-\frac 12,
\end{eqnarray*}
where $m_v=\int_V y\eta_v(dy)$ is the mean of $\eta_v$. Hence all the assumptions
of Lemma \ref{foverg} are satisfied, thus
$$
\lim_{c\to 0,c\in Z_{\ov\mu}} \frac 1{|c|^2|\log |c||}
\bigg(  \int_V\big(\chi_v (cx)-1\big)\cdot \eta_v\big( i (cx)^* \big)\nu(dx)
+ \int_V\big(\chi_v (cx)-1 -i\is v{cx}\big)\nu(dx) \bigg) = 2 C_2(v)
$$ and the Lemma follows.
\end{proof}
\begin{proof}[Proof of Theorem \ref{asymptkt}, part (4)]
First we will improve Lemma \ref{p361} and we will show that if $\a=2$ then
\begin{equation}
\label{pimp}
\big|  1 -\is{\nu}{\psi_{c,v}}\big|\le C |c|.
\end{equation}
Indeed, applying Lemma \ref{p21} and \eqref{h6}, we have
\begin{eqnarray*}
\big|  1 -\is{\nu}{\psi_{c,v}}\big|
&\le& \int_V\big| \psi_{c,v} (x)-\wh\eta_v(cx) \big|\nu(dx)
+ \int_V\int_V \big| \chi_y(cx) - 1 \big|\eta_v(dy)\nu(dx)\\
&\le& C|c|^{\l+2\e}\int_V |x|^{\l+\e}\nu(dx)+C|c|\int_V |x|\nu(dx)\int_V|y|\eta_v(dy)\\
&\le& C |c|,
\end{eqnarray*}
which proves \eqref{pimp}. Finally, applying Lemma \ref{alfa2} and \eqref{pimp},
we write
\begin{multline*}
\lim_{c\to 0, c\in Z_{\ov\mu}}\frac{k(c,v)-1-i\is v{cm}}{|c|^2|\log |c||} =
\lim_{c\to 0, c\in Z_{\ov\mu}}\frac{\nu(\psi_{c,v})(k(c,v)-1)-i\is v{cm}}{\nu(\psi_{c,v})|c|^2|\log |c||}\\
-\lim_{c\to 0, c\in Z_{\ov\mu}}\frac{i\is v{cm}(\nu(\psi_{c,v})-1)}{\nu(\psi_{c,v})|c|^2|\log |c||} = C_2(v).
\end{multline*}
\end{proof}
\subsection{Case: $\a>2$}
Here we replace   $Z_{\ov\mu}$ by $\R^*_+$, hence  $c=t\in\R_+^*$. We use
expression of $\psitv$ given by  \eqref{kt}
\begin{lem}
\label{uu13}
$$\lim_{t\to 0}\frac 1{t^2}\Big( \is{\nu}{(\chi_{tv}-1)\psitv} - it\is vm \Big) = C_{2+}^1(v),$$ where
$$C_{2+}^1(v) =-\frac 12 \int_{V}  \is vx^2 \nu(dx) - \int_{V}  \is {v}
{x} {\eta_{v}(x^*)   }\nu(dx).$$
\end{lem}
\begin{proof}
We write
\begin{eqnarray*}
\int_{V}  \big( \chi_{tv}(x)-1 \big)\psitv(x)\nu(dx)&=&
\int_{V}  \big(\chi_{tv}(x) -1 \big)\nu(dx)+
\int_{V}  \big(\chi_{tv}(x) -1 \big)\big( \psitv(x)-1 \big)\nu(dx)\\
&=& W_1(t) + W_2(t).
\end{eqnarray*}
Notice that for any $\d<1$ there exists $C$ such that
$$
\Big| e^{is}-1-is+\frac 12 s^2  \Big| \le C|s|^{2+\d}
$$
for every $s\in \R$. Therefore choosing $\d<\min\{ 1,\a-2 \}$
and applying \eqref{moments} with $2<\theta <\a$ we have
\begin{equation*}
\begin{split}
\lim_{t\to 0} \frac 1{t^2}\Big( W_1(t) -it\is vm \Big)
&= -\frac 12\int_V \is vx^2 \nu(dx) + \lim_{t\to 0} \frac 1{t^2}\int_V
\!\!\bigg( e^{it\is vx} - \! 1-it\is vx + \frac{t^2}2 \is vx^2  \bigg)\nu(dx)\\
&= -\frac 12 \int_V \is vx^2 \nu(dx).
\end{split}
\end{equation*}
To handle  $W_2$ we will prove first  that
\begin{equation}
\label{p62}
\|\chi_x-1-i x^*\|_{\th,\e,1} \le C|x|^{1+\d},
\end{equation}
for some $\d>0$. Indeed,
recall that in view of \eqref{h6} we may assume $\l=1$ and $1<\th<2$. We have
$|\chi_x-1-i x^*|_{\th} \le C|x|^{1+\d}$ and
\begin{multline*}
[\chi_x-1-i\is x{\cdot}]_{\eps,1}
=\sup_{y\not= y'} \frac{\Big|\big(\chi_x(y)  - 1 -i\is xy\big)-\big(\chi_x(y') - 1 -i\is x{y'}\big)\Big|}
{|y-y'|^{\eps}(1+|y|)(1+|y'|)}\\
\le
\sup_{y\not= y'} \min\Bigg\{
\frac{\big|\chi_x(y) - 1 -i\is xy\big|+\big|\chi_x(y') - 1 -i\is x{y'}\big|}{|y-y'|^{\eps}(1+|y|)(1+|y'|)},
\frac{\big|\chi_x(y) - \chi_x(y')\big| + \big|\is x{y-y'}\big|}{|y-y'|^{\eps}(1+|y|)(1+|y'|)}
\Bigg\}\\
\le
C \sup_{y\not= y'}\Bigg[\frac 1{(1+|y|)(1+|y'|)}\min\bigg\{
\frac{|x|^{1+\e}(|y|^{1+\e}+ |y'|^{1+\e} )}{|y-y'|^\e}, |x||y-y'|^{1-\e}\bigg\}\Bigg]\\
\le
C \sup_{y\not= y'}\Bigg[ \frac 1{(1+|y|)(1+|y'|)}\bigg(
\frac{|x|^{1+\e}(|y|^{1+\e}+ |y'|^{1+\e} )}{|y-y'|^\e}\bigg)^{1-\e}\cdot\Big( |x||y-y'|^{1-\e}\Big)^{\e}\Bigg]\\
\le C |x|^{1+\e-\e^2},
\end{multline*}
which proves \eqref{p62}.

Now applying Proposition \ref{lis23} and \eqref{p62} we have
\begin{eqnarray*}
\lim_{t\to 0} \frac{\psitv(x)-1}t &=& \lim_{t\to 0} \frac{\is{\eta_{t,v}}{\chi_{t,v}-1-itx^*}}{t}
+ i\lim_{t\to 0} \is{\eta_{t,v}-\eta_v}{x^*} + i\is{\eta_v}{x^*}\\
&=& i\eta_v(x^*).
\end{eqnarray*}
Therefore we have
$$
\lim_{t\to 0}\frac{W_2(t)}t = \lim_{t\to 0}\int_V \frac{\chi_{tv}(x)-1}t\cdot \frac {\psitv(x)-1}t\nu(dx)
=-\int_V \is vx \eta_v(x^*)\nu(dx)
$$
hence the Lemma.
\end{proof}
\begin{proof}[Proof of Theorem \ref{asymptkt}, part (5)] First we will prove
$$\lim_{t\to 0^+}\frac{\nu(\psi_{t ,v})-1-it\int_{V}\eta_v(x^*)\nu(dx)}t = 0.$$
Applying inequality \eqref{p62} and Proposition \ref{lis23} we have
\begin{multline*}
\nu(\psi_{t,v})-1-it\int_{V}\eta_v(x^*)\nu(dx) =
\int_{V} \Big(\eta_{t,v}(\chi_{tx} - 1-itx^*) +
it (\eta_{t,v}-\eta_{v}) (x^*)\Big)\nu(dx)\\
\le C\int_{V}\| \chi_{tx} - 1  - itx^*\|_{\tel}\nu(dx)
+ t^{1+\e}\int_{V}\|x^*\|_{\tel}\nu(dx)\le C t^{1+\d}
\end{multline*}
Therefore by \eqref{kt} and Lemma \eqref{uu13}
\begin{multline*}
\lim_{t\to 0^+}\frac{k(t,v)-1-i t\is vm }{t^2}
= \lim_{t\to 0^+}\bigg[\frac{\nu(\psi_{t,v})(k(t,v)-1)-it\is vm}{\nu(\psi_{t,v})t^2}
- \frac{i(\nu(\psi_{t,v})-1)t\is vm}{\nu(\psi_{t,v})t^2} \bigg]\\
= C_{2+}^1(v) - \lim_{t\to 0^+}\frac{it\is vm \big(\nu(\psi_{t,v})-1-it\int_{V}\eta_{v}(x^*)\nu(dx)\big)}{\nu(\psi_{t,v})t^2}
+ \lim_{t\to 0^+}\frac{ \is vm \cdot\int_{V} \eta_{v}(x^*)\nu(dx)}{\nu(\psitv)}\\
= -\frac 12 \int_V\is vx^2\nu(dx) -\int_V\is{x-m}{v}\eta_v(x^*)\nu(dx).
\end{multline*}
Finally, since $\kappa(1)<1$ the matrix $I-z^*=\E[I-M^*]$ is invertible,
therefore
\begin{eqnarray*}
\int_V\is{x-m}{v}\eta_v(x^*)\nu(dx)
&=&\int_V\is{x-m}{v}\E\bigg[ \bbis{x}{\sum_{k=1}^\8 M_0^*\ldots M^*_{k-1} v} \bigg]\nu(dx)\\
&=&\int_V\is{x-m}{v} \bbis{x}{\sum_{k=1}^\8 \big(z^*\big)^kv} \nu(dx)\\
&=&\int_V\is{x-m}{v} \bis{x-m}{ \big(I-z^*\big)^{-1}z^* v} \nu(dx)\\
&=& q(v,(I-z^*)^{-1}z^*v)
\end{eqnarray*}
Also
$$ \int_V \is vx^2 \nu(dx) = q(v,v)+\is vm^2,$$
which proves Theorem \ref{asymptkt}
\end{proof}

\subsection{Calculations of $C_\a(v)$ in terms of tails ($0<\a\le2$).}
Observe first that the function $\wt\L$ defined in Introduction  is $G_{\ov\mu}^*$ homogeneous, i.e. $\wt\L(g^* y)=|g|^\a \wt\L(y)$
if $g\in G_{\ov\mu}$.
Indeed for $\a<2$ this follows from \eqref{hom} and for $\a=2$ this was proved in Lemma \ref{foverg}.
As in \cite{BD} we define the polar coordinates $(a(x),\ov x)$ of $x\in V\setminus\{0\}$,
using the decomposition $G_{\ov\mu} = A_{\ov\mu}\ltimes K_{\ov\mu}$. We denote by
$\Sigma_1$ the natural fundamental domain of $A_{\ov\mu}$ on $V\setminus\{0\}$, i.e.
\begin{itemize}
\item $\Sigma_1 = \{ x: 1\le x < p \}$ if $R_{\ov\mu}=\langle p \rangle$,
\item $\Sigma_1=S_1$ the unit sphere, if $R_{\ov\mu}=\R^*_+$.
\end{itemize}
Then we write $x=a(x)\ov x$ with $a(x)\in A_{\ov\mu}$ and $ \ov x\in \Sigma_1$.
Then $r(x)=|a(x)|$ takes values in $R_{\ov\mu}$, and if $R_{\ov\mu}=\R_+$, $r(x)=|x|$.

We will write $\wt \L_s(y)= r^{s-\a}(y)\wt\L(y)$, so that $\wt\L_s(y) = r^s(y)\wt\L(\ov y)$
is well defined by its restriction to $\Sigma_1$, and is $G_{\ov\mu}^*$-homogeneous of degree $s$.
Also we denote $\wt\L^1(y)=\wt\L(\ov y){\bf 1}_{[1,\8)}(r(y))$. We recall
that the tail measure $\Delta_v$ of $\eta_v$ exists i.e.
\begin{equation}
\label{has1}
\D_v=\lim_{|g|\to 0,g\in G_{\ov\mu} }|g|^{-\a}(g^*\eta_v).
\end{equation}
Also, in view of Proposition \ref{2.6}
and Lemma \ref{3.6'}, $\D_v\not= 0$. 
\begin{prop}
\label{5.5}
Assume $R_{\ov\mu}=\R_+^*$, then
\begin{itemize}
\item
 if $\a\in (0,1)\cup (1,2]$ then $C_\a(v)= \a m_\a \D_{ v}(\wt\L^1)$;
\item
if $\a=1$, then $C_1(v)= m_1 \D_{ v}(\wt \L^1)+i \g(v)$, with $\g(v)\in\R$.
\end{itemize}
If $R_{\ov\mu}=\langle p \rangle$, the same formulas are valid, where $\a \D_{v}(\wt\L^1)$
is replaced by $\frac{1-p^{-\a}}{\log p}\D_v(\wt \L^1)$.
\end{prop}
\begin{proof}
If $\a\not= 1$, we have, by definition of $\wt\L$:
$$ C_\a(v) = \int_V \big( \wt\L(y+v)-\wt \L(y) \big) \eta_v(dy).
$$
Indeed if $\a\in(0,1)\cup (1,2)$, this follows immediately from the formulas given in Theorem
\ref{asymptkt} and for $\a=2$ we write
\begin{eqnarray*}
C_2(v) &=& -\frac 14 \int_{\Sigma_1}\Big(   \is wv^2 + 2\is wv \eta_v(w^*) \Big)\s(dw)\\
&=& -\frac 14 \int_V \int_{\Sigma_1}\Big( \is w{v+y}^2 - \is wy^2  \Big) \s(dw)\eta_v(dy)\\
&=&\int_V \big( \wt\L(y+v)-\wt \L(y) \big) \eta_v(dy).
\end{eqnarray*}
If $\a=1$, then
$$ C_1(v)   = \int_V \big( \wt\L(y+v)-\wt \L(y) \big) \eta_v(dy) + i\gamma(v)
$$ with
$$
\g(v) = \int_V \int_V \bigg( -\frac {\is vx}{1+|x|^2} -
\frac {\is yx}{1+|y|^2|x|^2} +\frac {\is {v+y}x}{1+|v+y|^2|x|^2}
\bigg) \L(dx)\eta_v(dy).
$$
We write if $s<\chi$
$$C_{\a,s} =  \int_V \big( \wt\L_s(y+v)-\wt \L_s(y) \big) \eta_v(dy)$$
and we observe: $\lim_{s\to\a^-} C_{\a,s}=C_\a(v)$. On the other
hand $Z^*v$ satisfies
$ Z^*v = M_0^* (Z_1^* v+v)$,
where $Z_1^*$ is another copy of $Z^*$, independent of $M_0^*$.
Since $\wt\L_s$ is $G_{\ov\mu}^*$-homogeneous and $s<\a$,
$$
\E\big[ \wt \L_s(Z^*v) \big] = \kappa(s)\E\big[ \wt \L_s(Z^*v+v) \big],
$$ hence
\begin{eqnarray*}
C_{\a,s} &=& \E\big[ \wt\L_s (Z^* v+v) \big] - \E\big[ \wt\L_s (Z^* v) \big]\\
&=& (1-\kappa(s))\E\big[ \wt\L_s (Z^* v+v) \big].
\end{eqnarray*}
Since $\lim_{s\to\a^-}\frac{1-\kappa(s)}{\a-s}=m_\a$, we need to evaluate
$\lim_{s\to\a^-}(\a-s) \E\big[ \wt\L_s (Z^* v+v) \big]$.
For the sake of brevity, we work with $\lim_{s\to\a^-}(\a-s)\E\big[ \wt\L_s (Z^* v) \big]$
 and we show
that this quantity depends only of the tail of $\eta_v$. This will
give the required result, since the tails of $\eta_v$ and $\d_v*\eta_v$
are the same.

Assume first $R_{\ov\mu}=\R^*_+$ and write $F_v(t) = \int_{|\zeta|\ge t}\wt\L (\ov\zeta)\eta_v(d\zeta)$.
Then $|F_v(t)|\le \sup_{\ov\zeta\in\Sigma_1}|\wt\L(\ov\zeta)|<\8$, hence $F_v$
is a bounded function. Also, if $g\in G_{\ov\mu}$, $|g|=t$:
$$
t^\a F_v(t) = \int_{|\zeta|\ge 1}\wt\L(\ov\zeta) |g|^\a(g^{-1}\eta_v)(d\zeta).
 $$
 Hence, using the convergence of $|g|^\a (g^{-1}\eta_v)$ to $\D_v$ if $|g|\to \8$,
\begin{equation}
\label{has3}
t^\a F_v(t) = \D_v(\wt\Lambda^1)+o(t), \quad \mbox{as }t\to\8.
\end{equation}
By definition of $F_v$:
\begin{eqnarray*}
\E\big[ \wt\L_s (Z^* v) \big] &=& \int_V |y|^s \wt\Lambda(\ov y)\eta_v(dy)\\
&=& \int_V\bigg( \int_{0<t<|y|} s t^{s-1} dt \bigg)\wt \Lambda(\ov y)\eta_v(dy)\\
&=& \int_0^\8 s F_v(t) t^{s-1}dt.
\end{eqnarray*}
Let $r$ be any positive increasing function on $(0,\a)$ satisfying
\begin{equation}
\label{has33}
\lim_{s\to\a^-} r(s)=+\8,\quad
\lim_{s\to\a^-}(\a-s) r^s(s)=0,\quad
\lim_{s\to\a^-} r^{s-\a}(s)=1.
\end{equation}
One can take for example $r(s)=(\a-s)^{-\frac 1{2\a}}$. Then to compute
the required limit we decompose the integral of $F_v$ above according to the
function $r(s)$ and apply \eqref{has3}, which gives the asymptotics of $F_v(t)$:
\begin{multline*}
\lim_{s\to \a^-}(\a-s) \E\big[ \wt\L_s(Z^*v) \big] = \lim_{s\to\a^-}(\a-s)\int_0^{r(s)}s F_v(t) t^{s-1} dt\\
+\lim_{s\to\a^-} (\a-s) \int_{r(s)}^\8 s \Delta_v(\wt\Lambda^1) t^{-\a+ s-1} dt
+ \lim_{s\to\a^-}(\a-s)\int_{r(s)}^\8 o(t) t^{-\a+s-1} dt.
\end{multline*}
Notice that the first and third limit are 0. Indeed, by \eqref{has33},
$$
\lim_{s\to\a^-}\bigg|(\a-s)\int_0^{r(s)}s F_v(t) t^{s-1} dt\bigg|
\le \lim_{s\to\a^-}(\a-s) r^s(s)\sup_{t>0}|F_v(t)|=0.
$$
To compute the third limit take for any $\e>0$, then there exists $s_0$ close to $\a$
such that $|o(t)|<\e$ for $t>r(s_0)$ then, by \eqref{has33}
$$
\lim_{s\to\a^-}\bigg|(\a-s)\int_{r(s)}^\8  o(t) t^{-\a+s-1} dt\bigg|
\le \e  \lim_{s\to\a^-} r^{s-\a}(s) = \e.
$$
Since $\e$ was arbitrary we obtain that the limit above is in fact 0.
 As a result, using \eqref{has33},
 \begin{eqnarray*}
 C_\a(v) &=& \lim_{s\to\a^-} (1-\kappa(s))\E\big[ \wt\L_s(Z^* v) \big] \\
 &=& m_\a \cdot \lim_{s\to\a^-} (\a-s) \int_{r(s)}^\8 s \Delta_v(\wt\Lambda^1) t^{-\a+ s-1} dt\\
 &=& \a m_\a \D_v(\wt\L^1).
 \end{eqnarray*}
 If $R_{\ov\mu}=\langle  p \rangle$, the calculation runs parallel, using the formula
 $$
 \E\big[ \wt\L_s(Z^* v) \big] = \int_V \wt\L_s(\zeta)\eta_v(d\zeta)
 $$
 we decompose $\{ \zeta\in V; |\zeta|>1 \}$ into shells of the form $\{ \zeta\in V; p^k\le \zeta\le p^{k+1} \}$
 and use geometric series instead of the integrals above. Using Theorem 1.4 of \cite{BD} which gives
 a formula for $\D_v$, we get
 $$ C_\a(v) = m_\a\frac{1-p^{-\a}}{\log p}\D_v(\wt\L^1). $$
\end{proof}
\section{Proof of Main Theorem \ref{mthm}}
\label{euclid-proof}
To prove the Theorem, in view of the continuity  theorem, it is enough to justify that the corresponding
characteristic functions converge pointwise to a function, which is continuous
at zero. If $\a<2$, as observed in Section 1, stability follows from the last assertions in
Theorem \ref{asymptkt} (1, 2 and 3).

\subsection{Case $\a<1$}
\label{6.1} Let $\phi_n^\a$ be the characteristic function of the random
variable $c_n S_n^x$.
Then by Lemma \ref{1e4} and Proposition \ref{kel}
we have
$$\phi^\a_{n}(v) = \E\Big[\chi_v(c_nS_n^x) \Big] =
\big(P^n_{c_n,v}(1)\big)(x) = k^n(c_n,v)\big(\pi_{c_n,v}(1)\big)(x)
+ \big(Q^n_{c_n,v}(1)\big)(x).
$$
The second factor tends to 0 as $n$ goes to infinity, because $\|Q_{c_n,v}\|<1$.
 Moreover, by Proposition \ref{lis22}, $\big(\pi_{c_n,v}(1)\big)(x)$ converges to 1.
Therefore it is enough to compute
$$
\lim_{n\to \8}k^n (c_n,v) =
\lim_{n\to \8}\Big( 1 + k(c_n,v )-1 \Big)^{
\frac{1}{k (c_n,v)-1}\cdot n ({k(c_n,v)-1})
} = \lim_{n\to \8}e^{n\cdot \big({k(c_n,v)-1}\big) }
$$
Notice that by \eqref{case1}
$$
\lim_{n\to \8}n\cdot \big({k(c_n,v)-1}\big) =
\lim_{n\to \8}  \frac{{k( c_n,v)-1}}
{|c_n|^\a} = C_\a(v).
$$
This proves pointwise convergence of $\phi_n^\a$ to $\Phi_\a$. Continuity of $\Phi_\a$
at zero follows
easily from the Lebesgue dominated theorem.

\subsection{Case $\a=1$}
Let $\phi_n^1$ be the characteristic function of $c_nS_n^x - n\xi(c_n)$. Then
arguing as above we prove that
\begin{eqnarray*}
\lim_{n\to \8}\phi_n^1(v) &=& \lim_{n\to \8} \E \Big[\chi_v (c_nS_n^x - n\xi(c_n))) \Big]
= \lim_{n\to \8} \Big[ \chi_v(-n \xi(c_n))\big( P_{c_n,v}^n(1)\big)(x)\Big]\\
&=& \lim_{n\to \8}\Big[ \chi_v(-\xi(c_n)) k(c_n,v)\Big]^n
= e^{ \lim_{n\to \8}\big[n\big(\chi_v(-\xi(c_n)) k(c_n,v)-1\big)\big]}
\end{eqnarray*}
Let us compute the limit in the exponent
\begin{multline*}
\lim_{n\to\8}\Big[n\big(\chi_v(-\xi(c_n)) k(c_n,v)-1\big)\Big]\\
= \lim_{n\to\8}\bigg[
\chi_v(-\xi(c_n))    \cdot\frac{k(c_n,v)-1-i\is v{\xi(c_n)}}{|c_n|} + n \chi_v(-\xi(c_n)) (1+i\is v{\xi(c_n)})-n\bigg]\\
= C_1(v) + \lim_{n\to\8}\Big[n\cdot \big(1-i\is v{ \xi(c_n)} + O(\is v{ \xi(c_n)}^2)\big)
\big(1+i\is v{\xi(c_n)} \big)-n \Big]
=C_1(v)
\end{multline*}
To
prove continuity of $\Phi_1$ at zero, it is enough to observe
$$ g_v(x) =  \big(e^{i\is vx} -1 \big)\cdot \wh{\eta}_v(x) - \frac{i\is vx}{1+|x|^2} =
\big(e^{i\is vx} -1 \big)\big(\wh{\eta}_v(x)-1\big) +  e^{i\is vx} -1 - \frac{i\is vx}{1+|x|^2}\le C|x|^{1+\d},
$$ for $|x|<1$ and some constants  $C$ and $\d>0$, independent of $v$,
and next one can apply the Lebesgue dominated theorem.

\subsection{Case $1 < \a < 2$}
Denote by $\phi_n^\a$ the characteristic function of $c_n(S_n^x-nm)$.
We reason as in previous cases and obtain
\begin{eqnarray*}
\lim_{n\to\8} \phi_n^\a(v)&=&
\lim_{n\to\8} \E\Big[ \chi_v(c_n(S_n^x-nm))\Big]
=\lim_{n\to\8} \Big[\chi_v(-c_n m) k(c_n,v) \Big]^n\\
&=&e^{\lim_{n\to\8} \big[n(\chi_v(-c_n m)k(c_n,v) -1)\big]},
\end{eqnarray*}
 and we have
\begin{multline*}
\lim_{n\to\8}\Big[ n\big(\chi_v(-c_n m) k(c_n,v)  - 1\big)\Big]\\
= \lim_{n\to\8}\bigg[ \chi_v(-c_n m) \cdot \frac{k(c_n,v)-1-i\is v{c_n m}}{|c_n|^\a}
+ n \chi_v(-c_n m)(1+i \is v{c_n m}) - n \bigg]\\
=C_\a(v) + \lim_{n\to \8}\Big[ n\big( 1-i\is v{c_n m} + O(n^{-\frac 2\a})\big)\big( 1+ i\is v{c_n m}\big)-n\Big]
 = C_\a(v).
\end{multline*}
To prove that $\Phi_\a$ is continuous at zero and stable, we proceed as before.
\subsection{Case $\a=2$}
Let $\phi_n^2$ be the characteristic function of $c_n(S_n^x-nm)$. Arguing as in previous
case we show
\begin{multline*}
\log\Phi_2(v) = \log \lim_{n\to\8} \phi_n^2(v)= \lim_{n\to\8}\Big[ n\big(\chi_v(-c_n m) k(c_n,v) - 1       \big)  \Big]\\
=\lim_{n\to\8}\bigg[\big(n|c_n|^2|\log|c_n||\big) \cdot \frac{ k(c_n,v)-1-ic_n\is vm }{|c_n|^2|\log |c_n||}\bigg]
+ \lim_{n\to\8} \Big[n \chi_v(-c_n m)(1+ic_n\is vm)-n\Big]\\
= C_2(v).
\end{multline*}
\subsection{Case $\a>2$}
\label{6.4}
We argue  as in the previous case.
Let $\phi_n^{2+}$ the characteristic function of $\frac 1{\sqrt n} (S^x_n-nm)$.
Then
$$
\lim_{n\to\8} \phi_n^{2+}(v)
=e^{\lim_{n\to\8} \big[n(\chi_v(-m/\sqrt n) k(1/\sqrt{n},v)-1)\big]}
$$
An elementary calculation, using the asymptotics of $k(1/\sqrt n,v)$ given in Theorem \ref{asymptkt}, 5) proves
$$\log \Phi_{2+}(v) = \lim_{n\to\8}\Big[ n\big( \chi_v(-m/\sqrt n) k(1/\sqrt n,v) - 1\big)\Big] = C_{2+}(v)+ \frac12 \is vm^2.$$
\subsection{Nondegeneracy of the limit law for $0<\a\le2$}
In order to prove that the limit law  is fully nondegenerate (i.e. its
support is not contained in some lower dimensional subspace of $V$) it is enough
to justify that the function $F_\a(v)=\Re\log \Phi_\a(v)$, defined on $V$,
does not vanishes outside zero.
We use the expression of $C_\a(v)$ given in Proposition \ref{5.5}.


\begin{prop}
 For every $v\in V\setminus\{0\}$, $F_\a(v)=\Re C_\a(v)$ is negative.
\end{prop}
\begin{proof}
If $R_{\ov\mu} = \R^*_+$, the expression of $C_\a(v)$ in Proposition \ref{5.5}
gives $F_\a(v) = \D_v(\Re \wt\L^1)$. The definition of $\wt\L$ gives
$$
\Re \wt\L(y) = \int_{V\setminus\{0\}} \big(\cos\is xy -1\big)\L(dx).
$$
Using Corollary \ref{2.9} and Lemma \ref{3.6'}, we know that supp$\L$
is not contained in a hyperplane. Since $(\cos \is xy -1)\le 0$, we get that for  any
$y\not= 0$, $\Re \wt\L(y)<0$. In particular, for any $y$ with
$|y|\ge 1$: $\Re \wt\L^1(y)<0$.
Since $\D_v$ is $G_{\ov\mu}^*$-homogeneous, nonzero and $G_{\ov\mu}^*$ is not compact,
we have ${\rm supp}\D_v \cap \{ y\in V; |y|\ge 1 \}\not=\emptyset$.
It follows $\D_v(\Re \wt\L^1)<0$ if $v\not=0$.

\medskip

If $R_{\ov\mu}=\langle p \rangle$, a simple modification of the argument above give the same result.
For $\a=2$ we reason analogously.

\end{proof}
\subsection{Nondegeneracy of the limit for $\a>2$}
Notice first that if $G_{\ov\mu}\subset \R^*_+$ then nondegeneracy of the limiting
random variable follows immediately from the formula of its characteristic function. Namely we may
write
\begin{eqnarray*}
-\log\Phi_{2+}(v) &=&
\frac 12 \int_{V}\is {x-m}v^2 \nu(dx) + \int_{V}\is {x-m}{v}\eta_v(x^*) \nu(dx)\\
&=&
\frac 12 \int_{V}\is {x-m}v^2 \nu(dx) + \int_{V}\is {x-m}{v}  \E\bis{x}{ \sum_1^\8 |M_1|\cdots|M_k|v }\nu(dx)\\
&=&
\bigg( \frac 12 +\sum_{n=1}^\8 \kappa^n(1) \bigg)  \int_{V}\is {x-m}v^2  \nu(dx)\\
&=& \frac 12\frac{1+\kappa(1)}{1-\kappa(1)}\ q(v,v),
\end{eqnarray*}
with $\frac{1+\kappa(1)}{1-\kappa(1)} >0$.
If the value above were zero, the support of $\nu$ would be contained in
some hyperplane of $V$ orthogonal to $v$, but this contradicts
to hypothesis {\bf H}.

\medskip

In general we cannot use the foregoing  argument hence we apply ideas of \cite{GH} (see also \cite{HH}).
Define $\sigma_v^2=-\log\Phi_{2+}(v)$, i.e. $\s_v^2$ is equal to the quadratic form $q(v,v)/2+q((I-a^*)^{-1}a^*v,v)$.
Since $e^{-\s_v^2}$ is the characteristic function of a probability measure, $\s_v^2\ge 0$.
Given a function $f$ on $V$ and $y\in V$ we define $f_y(x)=f(x-y)$

\begin{lem}
We have
$$
\sigma_v^2 = \frac 12 \nu \big((v^*_m)^2\big) + \nu(v^*\zeta),
$$
where $\zeta\in \B_{\tel}$
 is uniquely defined by the equations:
\begin{equation}
\label{poisson}
{\nu}({\zeta})=0, \quad
(I-P)(\zeta)=P\big(v^*_m\big).
\end{equation}
\end{lem}
\begin{proof}
Let $h_t$ be the eigenfunction of $P_{t,v}$:
\begin{equation}
\label{eigenf}
P_{t,v}(h_t) = k(t,v)h_t
\end{equation}
such that $\nu(h_t) =1$. The function $t\mapsto h_t=\frac{\pi
_{t,v}(e)}{\nu (\pi _{t,v}(e))} $ is differentiable in $\B_{\tel}$ for
appropriately chosen $\theta $ and $\lambda $ (see below). First we prove that
$P_{t,v}$ is differentiable. Let
$$
M_{t,v}f(x)=i\int \chi _{tv}(gx+b)\langle v,gx+b\rangle f(gx+b)\
\mu (dh). $$
 Then
\begin{equation}\label{M}
 \bigg\|\frac{P_{t+\Delta t,v}f-P_{t,v}f}{\Delta t}-M_{t,v}f\bigg\|
_{\theta , \eps , \lambda} \to 0\ \ \mbox{when}\ \Delta t\to 0
\end{equation}
 for $f\in \B_{\theta ', \eps ,\lambda'}$ with sufficiently
  small $\theta ', \lambda ' $. In particular \eqref{M} applies to
  $h_t$. Using the resolvent we write
 $$
 \pi _{t,v} = \frac{1}{2\pi i}\int _{|z-1|=\delta
 '}(z-P_{t,v})^{-1}\ dz$$
 and we differentiate $\pi _{t,v}(e)$. We need to take triples
 $(\theta ', \eps , \lambda ')$ and $(\theta , \eps , \lambda )$
 in the way that not only \eqref{M} is satisfied but also all the
 assumptions of section \ref{section2} to assure that the resolvent is
 bounded both on both $\B _{\theta ', \eps , \lambda '}$ and
 $\B _{\theta , \eps , \lambda }$. Taking $\eps $ sufficiently
 small, $\lambda '= 5\eps , \theta ' = 9\eps , \lambda = 1+10\eps
 , \theta =1+14 \eps $ will do. Clearly , $h_t \in \B _{\theta ',
 \eps , \lambda'}$. Finally,
 $$
 \frac{h_t - 1}{t}= \frac 1{\nu(\pi_{t,v}(e))}\bigg[
 \frac {\pi_{t,v}(e)-1}t + \nu\Big( \frac {1-\pi_{t,v}(e)}t \Big)
 \bigg]
$$ and
so $\lim _{t\to \8}\frac{h_t -1}{t}=\zeta$ exists.

We apply $\nu$
to both sides of \eqref{eigenf} and we obtain
\begin{equation}
\label{nukt}
\nu(\chi_{tv}h_t) = k(t,v).
\end{equation}
Next differentiating the equation $\nu(h_t)=1$ with respect to $t$ at 0 we obtain $\nu(\zeta)=0$.
Computing the second order term of asymptotic expansion of $k(t,v)$, in view of Theorem
\ref{asymptkt} we have
\begin{eqnarray*}
-\sigma_v^2 -\frac 12 \is vm^2&=&
\lim_{t\to 0}\frac{k(t,v)-1-it\is vm}{t^2}
= \lim_{t\to 0}\frac{\nu(\chi_{tv}h_t) -1-it\is vm}{t^2}\\
&=& \lim_{t\to 0}\nu\bigg(\frac{\big(1+itv^*-\frac{t^2}2 (v^*) ^2\big)h_t-1 -it v^* }{t^2}\bigg)\\
&=& -\frac 12 \nu((v^*)^2) + i \lim_{t\to 0}\nu\bigg(v^* \cdot \frac{h_t-1}{t}\bigg)\\
&=& -\frac 12 \nu((v^*)^2) - \nu\big( v^*\zeta  \big),
\end{eqnarray*}
that gives the required formula for $\sigma_v^2$. To prove that the function $\zeta$
satisfies the Poisson equation \eqref{poisson} we differentiate \eqref{eigenf} at zero, i.e.
applying Theorem \ref{asymptkt}
we write
$$
\lim_{t\to 0}\frac{k(t,v)h_t-1}t = \lim_{t\to 0}\bigg[ \frac{k(t,v)-1}t\cdot h_t + \frac{h_t-1}t \bigg]=i\big(\is vm+ \zeta\big).
$$
On the other side we obtain
$$\lim_{t\to 0}\frac{P_{t,v}(h_t)-1}t = \lim_{t\to 0}\bigg[\frac{(P_{t,v} - P)(h_t)}t + P\Big(\frac {h_t-1}t\Big)\bigg]
= i\big(P\big(v^*\big) + P(\zeta)\big).$$
Comparing both equations we prove \eqref{poisson}.

\medskip

Finally, in order to prove that $\zeta$ is uniquely determined by these two conditions,
assume that  some $\zeta_1$ satisfies $\nu(\zeta_1)=0$ and $P(v^*_m)=(I-P)(\zeta_1)$,
then $(I-P)(\zeta-\zeta_1)=0$, that implies $\zeta=\zeta_1+C$.
Since $\nu(\zeta)=\nu(\zeta_1)$, we get $\zeta=\zeta_1$.
\end{proof}

\begin{lem}
\label{6.9}
Let $u_0$ be the unique solution of $(I-P)u_0=v_m^*$, $\nu(u_0)=0$. Then $2\s_v^2 = \E_{\nu}\big[ u_0(X_1)-Pu_0(X_0) \big]^2$.
In particular, if $\sigma_v^2 = 0$ then $r(P_{t,v})=1$
\end{lem}
\begin{proof} Since $I-P$ is invertible on the space $\{ g:
\langle \nu , g\rangle =0\}$ , the system of equations
$\is{\nu}f=0$ and $(I-P)f=g$, $f\in\B_{\tel}$ has a unique
solution for $g$ such that $\is{\nu}g=0$. Therefore, equation
$(I-P)f=v_m^*$ has unique solution satisfying $\nu(f) =
0$, and we denote this solution by $u_0$. Then $\nu(u_0^2)<\8$.
Indeed, the function $u_0$  belongs to $\B_{\tel}$, therefore
$|u_0(x)|^2\le C (1+|x|)^{2+2\eps}$ and by \eqref{moments} $u_0^2$
is integrable with respect to $\nu$.
\medskip

Notice that $\zeta = Pu_0$. Indeed, it is enough to prove that
$Pu_0$ satisfies \eqref{poisson}. For this purpose we write
$$
(I-P)(Pu_0)=(I-P)(u_0-v_m^*) = (I-P)u_0 - v^*_m + P(v_m^*) = P(v_m^*),
$$
and
$$
\nu P(u_0) = \nu(u_0)=0.
$$
Next we write
\begin{eqnarray*}
2\sigma_v^2 &=& \nu\big((v^*_m)^2  \big) +2\nu\big( v^*_m \zeta \big)
= \nu\big((v^*_m)^2  +2 v^*_m Pu_0\big)
= \is{\nu}{ (u_0 - Pu_0)(u_0 + Pu_0)}\\
&=&\int_H\int_{V}\Big( u_0^2(h\cdot x) - (P u_0)^2(x)\Big)\nu(dx)\mu(dh)
= \int_H\int_{V}\Big( u_0(h\cdot x) - (P u_0)(x)\Big)^2\nu(dx)\mu(dh)\\
&=& \E_{\nu}\Big[ u_0(X_1) - Pu_0(X_0) \Big]^2
= \E_{\nu}\Big[ v_m^*(X_1) + Pu_0(X_1) - Pu_0(X_0) \Big]^2
\end{eqnarray*}
If $\sigma_v^2=0$, then  $$ v_m^*(X_1) = Pu_0(X_0) - P u_0(X_1),\quad \P_{\nu} \mbox{ a.s.}
$$
and
$$
e^{it\is v{X_1}}e^{it Pu_0(X_1)} = e^{it\is vm}e^{it Pu_0(X_0)},
$$ hence, taking the expected value, we have
$$
P_{t,v}(e^{itP u_0})(x) = \int_H e^{it\is v{h\cdot x}}e^{it Pu_0(h\cdot x)}\mu(dh) = e^{it\is vm}e^{it P u_0(x)}.
$$
that proves
$r(P_{t,v})=1$.
\end{proof}
Nondegeneracy of the limit follows immediately from Lemmas \ref{3.13} and \ref{6.9}.
\section{Proof of Theorem \ref{mthme}}
\label{k-gen}
In order to prove Theorem \ref{mthme} we proceed as in the
Euclidean case. However, now we have to handle with general dilations of $V$,
that requires some additional arguments. We omit these parts of the proof
that are similar in both cases. The crucial step is to describe asymptotic expansion
of $k(c,v)$ as $c$ goes to 0.

\subsection{Asymptotic expansion of $k(c,v)$ for $0<\a<2$}

\begin{prop}
\label{prop6.1}
\begin{enumerate}
\item If $0<\a<1$ then
$$
\lim_{|c|\to 0, c\in Z_{\ov\mu}}\frac{k(c,v)-1}{|c|^\a} = C_\a(v)
$$ where
$$ C_\a (v) = \int_{V} \Big(\chi_v(x)-1 \Big)\wh{\eta}_v(x)\L(dx).
$$
In particular $C_\a(c^* v)=|c|^\a C_\a(v)$, if $c\in Z_{\ov\mu}$.
\item
Assume $1\le \a<2$. Let $\xi_1 (c )=  cm_{\a,-} +  \int _{V}\frac{ cx_{\a }}{1+|cx_{\a }|^2} \nu (dx)$. Then
$$
\lim_{|c|\to 0, c\in Z_{\ov\mu}} \frac {k(c,v)-1 - i\langle v,\xi_2 (c)\rangle}{|c|^\a} =  C_\a (v),
$$ where
$$  C_\a(v)=\int_{V}\bigg( \big(\chi_v(x)-1\big)\cdot \wh{\eta}_v(x)-
i\langle v,x_{\a,-}\rangle- \frac{i \langle v,x_{\a }\rangle}{1+|x_{\a
}|^2}\bigg) \L(dx)$$ and
$\lim _{c\to 0}|c|^{-1}|\xi_1 (c)|= m_1$
for $m_1 = \int_V x_1\nu(dx)$.
In particular, if $c\in Z_{\ov\mu}$,  $C_\a(c^* v)=|c|^\a C_\a(v)+i\is v{\b(c)}$, with
$\b(c)=\int_V\big( \frac{x_\a}{1+|x_\a|^2}  - \frac{x_\a}{1+|c x_\a|^2} \big)\L(dx)$.
\end{enumerate}
\end{prop}
\begin{proof}
For $\a <1$ the proof is exactly the same as in  section \ref{k-euclid}.
Assume $1\le \a<2$. First we will prove that
\begin{equation}
\label{lem6.2}
\lim_{|c|\to 0} \frac 1{|c|^\a}\Big( \is{\nu}{(\chi_{c^*v} -1)\psi_{c,v}} - i\is{v}{\xi_1(c)} \Big) = C_\a(v).
\end{equation}
For this purpose we decompose $V=V_{\a,-} \oplus V_{\a } \oplus V_{\a,+}$.
and write $x=x_{\a,-}+x_{\a }+x_{\a,+}$.
Then
\begin{multline*}
\int_{V}\big( \chi_v(cx) - 1 \big)\psi_{c,v}(x) \nu(dx) =
\int_{V}\big(  \chi_v(cx) - 1\big)\big(\psi_{c,v}(x) - \wh{\eta}_v(cx)\big)\nu(dx)\\
+ \int_{V}\big( \chi_v(cx)-1\big)\big(\wh{\eta}_v(cx)-1\big) \nu(dx)
+\int_{V}\bigg( \chi_v(cx)-1- i\is{v}{cx_{\a,-}}
-\frac{i\is{v}{cx_\a}}{1+|cx_\a|^2}\bigg)\nu(dx)\\
+i\is {v}{\xi_1(c)},
\end{multline*}
To handle the first and the second integrals we use the same arguments
as in the proof of Theorem \eqref{asymptkt}, i.e. the first one converges to
0 as $c$ goes to 0, and the second one tends to
$$\int_V \big(\chi_v(x)-1\big)\big( \wh{\eta}_v(x)-1  \big)  \nu(dx) $$
 For the third one
we are going to prove that
$$f(x)= \chi_v(x)-1-i\langle v,x_{\a,-}\rangle- \frac{i
\langle v,x_{\a }\rangle}{1+|x_{\a }|^2}$$ satisfies \eqref{beta}.
 Let $D_+=\min \{ \l_j: \l_j >\a \}$,
$D_-=\max \{ \l_j: \l_j <\a \}$ .
For $\tau(x)\geq 1$ we have
 $|f(x)|\leq C(1+ \sum_{\l_j<\a}|x_{\l_j}|)\le  C\tau(x)^{D_-}$.
 Then for $\tau(x)\leq 1$ we have
$\sum_{\l_j>\a}|x_{\l_j}|\le C \tau(x)^{D_+}$, $|x_\a|\le C\tau(x)^\a$
 and we obtain
\begin{eqnarray*}
|f(x)| &\le& \big| e^{i\is vx} - 1-i\is vx \big| + |\is v{x_{\a,+}}|
+|\is v{x_\a}|\cdot \Big| 1- \frac 1{1+|x_\a|^2}\Big|\\
&\le& C\tau(x)^2 + C \tau(x)^{D_+} + C\tau(x)^\a\cdot  {|x_\a|^2}
\le C \tau(x)^{\min\{2,D_+\}}.
\end{eqnarray*}
Hence \eqref{has} implies \eqref{lem6.2}.
One can easily prove that for $|c|<1/2$ we have
$\big|\int_{V}\frac{|c|^\a x_\a}{1+|c|^{2\a}|x_\a|^2}\nu(dx)\big| < |c|^\a|\log |c||$
(compare proof of Lemma \ref{ksit}) and
$\lim_{|c|\to 0}\frac{\xi_1(c)}{|c|} = m_1$.
Finally by \eqref{h6} and Lemma \ref{p361} we have
\begin{equation*}
\begin{split}
 \lim_{|c|\to 0}\frac{k(c,v)-1-i\is v{\xi_1(c)}}{|c|^\a}
&=\lim_{|c|\to 0}\bigg(\frac{\nu(\psi_{c,v})(k(c,v)-1)-i\is v{\xi_1(c)}}{\nu(\psi_{c,v})|c|^\a}
+ \frac{i \langle v , \xi_1 (a)\rangle (1-\nu
(\psi _{c,v}))} {\nu (\psi _{c,v})|c|^{\a}}\bigg) \\&= C_\a(v)
\end{split}
\end{equation*}
\end{proof}
\subsection{Asymptotic expansion of $k(c,v)$ for $\a>2$ }
In order to get fully nondegenerate laws, we have
to normalize  $S_n^x$ in inhomogeneous way. Let
\begin{eqnarray*}
V_-&=& \vm = \oplus_{\l_j < \frac \a2} V_{\l_j},\\
V_+&=& \vp = \oplus_{\l_j > \frac \a2} V_{\l_j}.
\end{eqnarray*}
We assume that $V_{\frac \a2}=\{0\}$ and so $V=V_-\oplus V_+$.
For $x\in V$ we write $x=x_-+x_+$, where $x_-\in V_-$, $x_+\in V_+$.
Let $c_n\in Z_{\ov\mu}$ be such that
$$ |c_n | = \sup_{\tau(x)\le 1}\tau(c_nx)=\frac 1{n^{\frac 1\a}}.$$
The right normalization in the case $\a>2$ is
$$\frac 1{\sqrt n}\big(S_n^x-nm\big)_-+ \big( c_n S_n^x-d_n \big)_+.$$
We need to modify accordingly the operators $P_{c,v}$ and $T_{c,v}$
and so we consider the following linear transformations:
\begin{eqnarray*}
b_n(x)=\frac 1{\sqrt n}\ x_-,\\
c_n(x)=c_n(x_+),
\end{eqnarray*}
and
$$
a_n(x) = b_n(x) + c_n(x).
$$
Notice that the operators $P_{b_n,v}$ and
$T_{b_n,v}$ are defined both on $\B_{\tel}(V)$
and $\B_{\tel}(V_-)$. If $f(x)=f(x_-)$ then $f\in \B_{\tel}(V_-)$
if and only if $f\in \B_{\tel}(V)$ and
\begin{eqnarray*}
P_{b_n,v}f(x) &=& P_{b_n,v_-}f(x_-),\\
T_{b_n,v}f(x) &=& T_{b_n,v_-}f(x_-).
\end{eqnarray*}
The same holds for $P_{c_n,v}$, $T_{c_n,v}$ and functions $f$
depending only on $x_+$. Therefore, it is convenient to refer to
operators $P_{b_n,v}$, $P_{c_n,v}$, $T_{b_n,v}$, $T_{c_n,v}$ while
they act on $\B_{\tel}(V)$ and to $P_{b_n,v_-}$, $T_{b_n,v_-}$,
$P_{c_n,v_+}$, $T_{c_n,v_+}$ when they are considered on
$\B_{\tel}(V_-)$ or $\B_{\tel}(V_+)$. Clearly, the peripherical
eigenvalues $k(b_n,v_-)$, $k(c_n,v_+)$ of $P_{b_n,v_-}$ and
$P_{c_n,v_+}$ are equal to the peripherical eigenvalues
 $k(b_n,v)$, $k(c_n,v)$
of $P_{b_n,v}$ and $P_{c_n,v}$, respectively. They are closely related to $k(a_n,v)$.
Indeed, we are going to prove that
$$
\lim_{n\to\8} n\big(k(a_n,v)-k(b_n,v)-k(c_n,v)+1\big)=0.
$$
Moreover we describe the asymptotic behavior of $k(b_n,v_-)$
and $k(c_n,v_+)$ restricting our attention to $\B_{\tel}(V_-)$ and
$\B_{\tel}(V_+)$.

\medskip

The content of Section \ref{section2} is needed for both operators
$P_{b_n,v_-}$ and $P_{c_n,v_+}$. The inhomogeneous dilations in
$P_{b_n,v_-}$ imply a slight modification of Propositions
\ref{lis22} and \ref{lis23}. They hold with $|b_n|^\e =\big(\frac
1{\sqrt n}\big)^{\frac \e{\l_{k_0}}}$ instead of $|c|^\e$, where
$\l_{k_0}=\max\{\l_j:\; \l_j<\frac \a2\}$.

\medskip

From now on we assume not only that: $\e<1$, $\th\le 2\l$,
$\l+3\e<\th$, $2\l+3\e<\a$, but also $\l>\l_{k_0}$, $\frac
\a2<\l+2\e< \l_{k_0+1}$ (notice that $\l_{k_0+1}=\min\{ \l_j:\; \l_j>\frac \a2 \}$).

To study $P_{c_n,v_+}$ we need a further decomposition of $V_+$,
i.e.
$$ V_+ = V_{\frac \a 2,\a} \oplus V_\a \oplus V_{\a,+} $$
and for $x_+\in V_+$ we write $x_+=u=u_-+u_\a+u_+$. Let $\nu_-$,
$\nu_+$ be projections of the Poisson kernel $\nu$ on $V_-$ and
$V_+$ and let $\Lambda_+$ be the tail measure \eqref{tailm} for $\nu_+$.

\begin{prop}
\label{prop6.2}
If $\a>2$, then
\begin{equation}
\label{7.4}
\lim_{n\to\8}n\big( k\big(\dt,v\big)-1-i\langle \dt(v),m_0 \rangle \big) =  C_-(v),
\end{equation} where
$$
 C_-(v) =  -\frac 12 q(v_-,v_-) -\frac 12 \is {v_-}{m_0}^2 - q(v_-,(I-z_-^*)^{-1}z_-^*v_-),
$$
 $m_0=m_{\frac \a 2,-}=\int_{V_-}x\nu_-(dx)$ and $z_-=\E\big[ M_n|_{V_-} \big]$.
Moreover
\begin{equation}
\label{7.5}
\lim_{n\to\8} n\big(k(c_n,v)-1-i\is v{\xi_2(c_n)}\big) =  C_+(v),
\end{equation}
where
\begin{eqnarray*}
  C_+(v)&=&\int_{V_+}\bigg( \big(\chi_v(u)-1\big)\cdot \wh{\eta}_v(u)-
i\langle v,u_-\rangle- \frac{i \langle v,u_{\a }\rangle}{1+|u_{\a
}|^2}\bigg) \L_+(du),\\
\xi_2(c_n)&=& \int_{V_+}\bigg( c_n u_- + \frac{c_n u_\a}{1+|c_n u_\a|^2} \bigg) \nu_+(du)
\end{eqnarray*}
and if $V_{\frac \a2,\a}\not=\{0\}$, then
$$
\lim_{n\to\8} \frac{\xi_2(c_n)}{|c_n|^{\l_{k_0+1}}} = \int_{V_+}u_{\l_{k_0+1}}\nu_+(du).
$$
\end{prop}
\begin{proof}
To prove \eqref{7.4} we consider $P_{b_n,v_-}$ and we proceed as in the proof
of Theorem \ref{asymptkt} (5). The crucial estimate is
\begin{equation}
\label{7.6} |x_-|\le (1+\tau(x_-))^{\l_{k_0}},
\end{equation}
where $|\cdot|$ is the Euclidean norm on $V_-$, which
implies
\begin{equation}
\label{7.8}
\int_V |x_-|^2\nu(dx)<\8
\end{equation}
and so the integrals $\int_V\is{v_-}{x_-}^2\nu(dx)$,
$\int_V|\is{v_-}{x_-}|\eta_{v_-}(x_-^*)\nu(dx)$, $\int_V|\eta_{v_-}(x_-^*)|\nu(dx)$
are finite. Moreover in view of \eqref{7.6} we have
$$ \| \chi_{x_-}-1-ix_-^* \|_{\tel} \le C |x_-|^{1+\e-\e^2},
$$ which is also needed.
\medskip

To prove \eqref{7.5} we proceed as in Proposition \ref{prop6.1}, that is we prove
\begin{equation}
\label{7.9} \lim_{n\to\8}\frac 1{|c_n|^\a}\int_{V_+}\big( \chi_{v_+}(c_nu)-1 \big)
\big( \psi_{c,v_+}(u) -\widehat\eta_{v_+}(c_n u) \big)\nu_2(du) =0,
\end{equation}
\begin{equation}
\label{7.10} \lim_{n\to\8}\frac 1{|c_n|^\a}\int_{V_+}\big( \chi_{v_+}(c_nu)-1 \big)\big( \widehat \eta_{v_+}(c_nu)-1\big)
 \big)\nu_2(du) = \int_{V_+}\big( \chi_{v_+}(u)-1 \big)\big( \widehat \eta_{v_+}(u)-1\big)
 \big)\L_2(du),
\end{equation}
\begin{multline}
\label{7.11}
\lim_{n\to\8}\frac 1{|c_n|^\a}\int_{V_+}\bigg( \chi_{v_+}(c_nu)-1-i\is{v_+}{c_n u} - \frac{i\is v{c_n u_\a}}{1+|c_nu_\a|^2}
 \bigg)\nu_2(du) \\
= \int_{V_+}\bigg( \chi_{v_+}(u)-1-i\is{v_+}{ u} - \frac{i\is v{ u_\a}}{1+|u_\a|^2}
 \bigg)\L_2(du).
 \end{multline}
For \eqref{7.9}  we need
\begin{eqnarray}
\label{7.14} |\psi_{c_n,v_+}(u)-\widehat\eta_v(c_n u)|&\le&
C|c_n|^\e\tau(c_n u)^{\l+\e}\ \ \ \mbox{if}\ \tau(c_n u)\le 1\\
\label{7.16} |\psi_{c_n,v_+}(u)-\widehat\eta_v(c_n u)|&\le&
|c_n|^\e\tau(c_n u)^{\e},\ \ \ \mbox{if}\ \tau(c_n u)\ge 1.
\end{eqnarray}
\eqref{7.16} was proved in Lemma \eqref{p21} and \eqref{7.14} follows
from \eqref{7.19} below. Moreover, the assumption $\l+2\e>\frac \a2$
guaranties $\l_{k_0+1}+2\e+\l>\a$, that is used in the calculations.

The function in \eqref{7.10} satisfies \ref{beta}, because
$$
\big |\chi_{v_+}(u)-1\big |\big | \widehat \eta_{v_+}(u)-1\big |
\le \max (1,| u_+ |^2) \le \max (1,\tau ( u_+ )^{2\l_{k_0+1}})$$
and $\a <2\l_{k_0+1}$.
Therefore, \eqref{7.10} follows. Finally,
$\chi_{v_+}(u)-1-i\is{v_+}{ u} - \frac{i\is v{ u_\a}}{1+|u_\a|^2}$
is estimated as in the proof of Proposition \ref{prop6.1}

 \end{proof}
In order to compare $k(a_n, v)$ with $k(b_n,v)$ and $k(c_n,v)$ we need two technical lemmas.
\begin{lem}
\label{ren2}
For every $s\le \l_{k_0+1}$
\begin{equation}
\label{7.19}
\big| \chi_y(a_nx)-1 \big| \le C\Big( \frac 1{\sqrt n}|x_-||y_-| + \frac 1{n^{\frac s\a}}
\tau(x_+)^s\tau(y_+)^s
\Big).
\end{equation}
Moreover
\begin{equation}
\label{7.20}
\big\| \chi_{a_nx}-1 \big\|_{\tel} \le C\Big( \frac 1{\sqrt n}|x_-| + \frac 1{n^{\frac {\l+\eps}\a}} \tau(x_+)^{\l+\eps}
\Big)
\end{equation}
\end{lem}
\begin{proof}
We use the following inequality
\begin{equation}
\label{7.21}
\Big| e^{i\sum_{j=1}^m\a_j}-e^{i\sum_{j=1}^m\b_j}\Big| \le \sum_{j=1}^m\Big| e^{i\a_j}-e^{i\b_j}\Big|,
\end{equation}
which holds for real $\a_j$, $\b_j$, $1\le j\le m$.

In view of \eqref{7.21} we have
\begin{eqnarray*}
\big| \chi_y(a_nx)-1 \big| &\le& \Big|  e^{i\is y {\dt x}} -1 \Big|+\Big|  e^{i\is y {c_n x}} -1\Big|\\
&\le& |y_-||b_n\xm| + \sum_{j>k_0} \Big| e^{i\is {y_{\l_j}}{c_n x_{\l_j}}}- 1  \Big|\\
&\le& \frac 1{\sqrt n} |\ym||\xm| + \sum_{j>k_0} \Big(|y_{\l_j}||c_n x_{\l_j}|  \Big)^{\frac s{\l_j}}\\
&\le& \frac 1{\sqrt n} |\ym||\xm|+ \frac 1{n^{\frac s\a}}\tau(\yp)^s\tau(\xp)^s,
\end{eqnarray*} and \eqref{7.19} follows.

\medskip

For \eqref{7.20} we have
$$ |\ym|\le (1+\tau(y_-))^{\l_{k_0}} \le(1+\tau(y_-))^\t
$$
and
$$
\tau(y_+)^{\l+\e}\le (1+\tau(y_+))^\t,
 $$ hence we obtain the required estimate for $|\chi_{a_n x}-1|_\t$. Applying
 \eqref{7.21} again we have

\begin{equation}
\label{7.22}
\big| \chi_{a_nx}(y)-\chi_{a_n x}(y') \big|\le
\big| \chi_{b_nx}(y)-\chi_{b_n x}(y') \big|+\big| \chi_{c_nx}(y)-\chi_{c_n x}(y') \big|.
\end{equation}
Now
$$
\big| \chi_{b_nx}(y)-\chi_{b_n x}(y') \big| = |\chi_{b_n x}(y-y')|
\le |b_n x_-||y_--y_-'|\le \frac 1{\sqrt n} |x_-||y_--y_-'|.
$$
If $\tau(y_--y_-')\le 1$, then $|y_--y_-'|\le C\tau(y_--y_-')^{\e}$.
If $\tau(y_--y_-')\ge 1$, then $|y_--y_-'|\le C\tau(y_--y_-')^{\l_{k_0}}$.
Since $\l_{k_0}<\l$ in both cases
$$
|y_--y_-'|\le \tau(y_--y_-')^\e(1+\tau(y_-))^\l(1+\tau(y_-'))^\l
$$
and so
$$
\big[ \chi_{b_nx} \big]_{\e,\l}\le \frac C{\sqrt n}|\xm|.
$$
For the second term in \eqref{7.22} we apply \eqref{7.19} with $s=\l+\e<\frac \a2<\l_{k_0+1}$
and we have
$$
\big| \chi_{c_nx}(y)-\chi_{c_n x}(y') \big| = \big| \chi_{c_n x}(y-y')-1 \big|
\le \frac C{n^{\frac{\l+\e}\a}}\tau(\xp)^{\l+\e}\tau(y-y')^{\l+\e},
$$
but $\tau(y-y')^{\l+\e} \le \tau(y-y')^{\e}(1+\tau(y))^\l(1+\tau(y'))^\l$ and
\eqref{7.20} follow.
\end{proof}
\begin{lem}
\label{7.23} If $\l+2\e < \l_{k_0+1}$, then
\begin{equation}
\label{7.24}
\big| \psi_{a_n,v}(x)-\psi_{b_n,v}(x) \big| \le \frac 1{n^{\frac {\l+2\e}\a}}(1+\tau(x))^{\l+2\e}
\end{equation}
and
\begin{equation}
\label{7.25}
\big| \psi_{a_n,v}(x)-\psi_{c_n,v}(x) \big| \le \frac 1{\sqrt n}(1+\tau(x))^{\l+\e}.
\end{equation}
\end{lem}
\begin{proof}
We have
\begin{eqnarray*}
\psi_{a_n,v}(x)-\psi_{b_n,v}(x) &=& \eta_{a_n,v}(\chi_{a_n x}-1) -  \eta_{b_n,v}(\chi_{b_n x}-1)\\
 &=& (\eta_{a_n,v}-\eta_v)(\chi_{a_n x}-1)   - ( \eta_{b_n,v}-\eta_v)(\chi_{b_n x}-1)
 + \eta_v(\chi_{b_n x}(\chi_{c_n x}-1)).
\end{eqnarray*}
Analogously
$$
\psi_{a_n,v}(x)-\psi_{c_n,v}(x) =
 (\eta_{a_n,v}-\eta_v)(\chi_{a_n x}-1)   - ( \eta_{c_n,v}-\eta_v)(\chi_{c_n x}-1) + \eta_v(\chi_{c_n x}(\chi_{b_n x}-1)).
$$
As it was mentioned before
\begin{eqnarray*}
|\eta_{b_n,v}-\eta_v|(\chi_{b_n x}-1) &\le& \Big( \frac 1{\sqrt n} \Big)^{\frac \e{\l_{k_0}}} \| \chi_{b_n x}-1 \|_{\tel},\\
|\eta_{c_n,v}-\eta_v|(\chi_{c_n x}-1) &\le& \frac 1{ n^{\frac \e\a}}\| \chi_{c_n x}-1 \|_{\tel},\\
|\eta_{a_n,v}-\eta_v|(\chi_{a_n x}-1) &\le& \frac 1{ n^{\frac \e\a}}\| \chi_{a_n x}-1 \|_{\tel},\\
\end{eqnarray*}
because $\frac \e{2 \l_{k_0}}>\frac \e\a$.

Now we apply Lemma \ref{7.19} with $s=\l+2\e$ and since $\frac{\l+\e}\a<\frac 12<\frac{\l+2\e}\a$
and $|\xm|\le C (1+\tau(x)^{\l_{k_0}})$, \eqref{7.24} and \eqref{7.25} follow.
\end{proof}
\begin{prop}
\label{ren3}
We have
$$
\lim_{n\to\8} n\big( k(a_n,v)-k(\dt,v)-k(c_n,v)+1
\big)=0
$$
\end{prop}
\begin{proof}
By \eqref{kt} we have
\begin{multline*}
 k(a_n,v)-k(\dt,v)-k(c_n,v)+1\\ = \frac 1{\nu(\psi_{a_n,v})}\nu\Big( \psi_{a_n,v}\big( \chi_{a_n^*v}-1 \big) \Big)
-  \frac 1{\nu(\psi_{b_n,v})}\nu\Big( \psi_{b_n,v}\big( \chi_{b_n^*v}-1 \big) \Big)
-   \frac 1{\nu(\psi_{c_n,v})}\nu\Big( \psi_{c_n,v}\big( \chi_{c_n^*v}-1 \big) \Big)\\
=I_1+I_2+I_3,
\end{multline*}
where
\begin{eqnarray*}
I_1&=& \frac 1{\nu(\psi_{a_n,v})}\nu\Big( \psi_{a_n,v}\big( \chi_{b_n^*v}-1 \big)\big( \chi_{c_n^*v}-1 \big) \Big),\\
I_2&=& \frac 1{\nu(\psi_{a_n,v})}\nu\Big( \psi_{a_n,v}\big( \chi_{b_n^*v}-1 \big) \Big)
- \frac 1{\nu(\psi_{b_n,v})}\nu\Big( \psi_{b_n,v}\big( \chi_{b_n^*v}-1 \big) \Big),\\
I_3&=& \frac 1{\nu(\psi_{a_n,v})}\nu\Big( \psi_{a_n,v}\big( \chi_{c_n^*v}-1 \big) \Big)
- \frac 1{\nu(\psi_{c_n,v})}\nu\Big( \psi_{c_n,v}\big( \chi_{c_n^*v}-1 \big) \Big).
\end{eqnarray*}
Applying Lemma \ref{7.19} with $s=\l+2\e$ we have
\begin{eqnarray*}
|I_1| &\le& \int_V \frac 1{\sqrt n}|\xm|\frac 1{n^{\frac s\a}} \tau(\xp)^s\nu(dx)\\
&\le& \frac 1{n^{\frac 12+\frac s\a}}\int_V \big( 1+\tau(x)^{\l_{k_0}+\l+2\e} \big)\nu(dx)\\
&=& o(\frac{1}{n}),
\end{eqnarray*}
because $\frac s\a > \frac 12$ and $\l_{k_0} + \l+2\e<2\l+2\e<\a$.

For $I_2$ we have
$$
I_2= \frac 1{\nu(\psi_{a_n,v})}\nu\Big( \big(\psi_{a_n,v}-\psi_{b_n,v}\big)\big( \chi_{b_n^*v}-1 \big) \Big)
+ \frac {\nu(\psi_{b_n,v}-\psi_{a_n,v})}{\nu(\psi_{a_n,v})\nu(\psi_{b_n,v})}
\nu\Big( \psi_{b_n,v}\big( \chi_{b_n^*v}-1 \big) \Big).
$$
Therefore by Lemmas \ref{7.19} and \ref{7.23}
\begin{eqnarray*}
|I_2| &\le& \frac 1{n^{\frac{\l+2\e}\a +\frac 12}}\bigg(
\int_V \big( 1+\tau(x)^{\l+2\e} \big)|x_-|\nu(dx) + \int_V \big( 1+\tau(x)^{\l+2\e} \big)\nu(dx)
\cdot \int_V|\xm|\nu(dx)
\bigg)\\
&=& o(\frac{1}{n}),
\end{eqnarray*}
because $\l +2\eps >\frac{\a }{2}$ and $|x_-|\leq (1+\tau (x_-))^{\l_{k_0}}$. For $I_3$ we have
$$
I_3= \frac 1{\nu(\psi_{a_n,v})}\nu\Big( \big(\psi_{a_n,v}-\psi_{c_n,v}\big)\big( \chi_{c_n^*v}-1 \big) \Big)
+ \frac {\nu(\psi_{c_n,v}-\psi_{a_n,v})}{\nu(\psi_{c_n,v})\nu(\psi_{b_n,v})}
\nu\Big( \psi_{c_n,v}\big( \chi_{c_n^*v}-1 \big) \Big)
$$ and so by Lemmas \ref{7.19} and \ref{7.23}, with $s=\l+2\e$
\begin{eqnarray*}
|I_3| &\le& \frac 1{n^{\frac 12 + \frac {\l+2\e}\a}}\bigg(
\int_V\big(1+\tau(x)^{\l+\e}\big)\tau(x)^{\l+2\e}\nu(dx)\!
+  \!\!\int_V\big(1+\tau(x)^{\l+\e}\big)\nu(dx) \int_V\tau(x)^{\l+2\e}\nu(dx)\bigg)\\
&=& o(\frac{1}{n}),
\end{eqnarray*}
because $2\l+3\e<\a$ and both integrals are finite.
\end{proof}

\begin{proof}[Proof of Theorem \ref{mthme}]
For $\a<2$ the Theorem is an immediate consequence of
Proposition  \ref{prop6.1}. To prove existence of limit of appropriately
normed sums $S_n$ we proceed exactly as in paragraphs \ref{6.1} -
\ref{6.4}. Also continuity at 0 of the characteristic function and stability
require only a repetition of the previous arguments, that will be omitted.
 Finally we have to justify nondegeneracy of the
limiting random variable. For this purpose take $v\in V_{\gamma}$
for some nonempty subspace $V_{\gamma}$ of $V$. Notice that
$\is{c_nS_n^x-d_n}{v}=\is{\pi_{\gamma}(c_nS_n^x-d_n)}{v}$,
where $\pi_{\gamma}$ denotes projection onto $V_{\gamma}$. Then
$\pi_{\gamma}(S_n^x)$ are partial sums of $X_n^{\gamma}$ defined by
the recursion
$X_n^{\gamma}=M_n^{\gamma}X_{n-1}^{\gamma}+\pi_{\gamma}(Q_n)$,
where $M_n^{\gamma}$ is the restriction of the action of $M_n$ to
$V_{\gamma}$ (it is well defined, because $V_{\gamma}$ is
invariant under the action of $G_{\ov\mu}$). The law $\mu_\g $ of
$(\pi_{\gamma}(Q_n), M^{\gamma }_n)$ is the projection of $\mu$ on
$\pi_\g(G_{\ov\mu} )$ under the
natural homomorphism.  Since $\mu $ doesn't admit invariant affine
subspaces, there is $\gamma $ such that there is no affine
subspace invariant under $\pi_\g(G_{\ov\mu} )$. Then we have to study $\pi_{\gamma}(S_n)$ on
$V_{\gamma}$ i.e. we reduce the problem in fact to the Euclidean
settings and last part of Theorem \ref{mthm} implies that the
limit $\is{\pi_\g(c_nS_n^x-d_n)}{v}$
is nonzero.

\medskip

If $\a>2$ we proceed as previously, however for the reader convenience and
to underline the role of Proposition \ref{ren3} we will present part of
the proof in more details.

\medskip

Let $\phi_n$ be the characteristic function of $$\dt (S_n^x-nm)_{\frac \a2,-}-(c_n S_n^x-d_n)_{\frac \a2,+}=
a_nS_n^x - nb_nm_0 - d_n.$$
Then
\begin{eqnarray*}
\lim_{n\to\8} \phi_n(v) &=& \lim_{n\to\8} \E\big[ \chi_v (a_n S_n^x-n\dt m_0 - n\xi_2(c_n)) \big]\\
 &=& \lim_{n\to\8} \Big[ \chi_v (-\dt m_0 - \xi_2(c_n)) k(a_n,v) \Big]^n\\
 &=& e^{\lim_{n\to\8} \big[ n \big(\chi_v (-\dt m_0 - \xi_2(c_n)) k(a_n,v) -1\big) \big]}.\\
\end{eqnarray*}
Applying Propositions \ref{prop6.2} and \ref{ren3} we obtain
\begin{eqnarray*}
\log\Big( \lim_{n\to\8} \phi_n(v) \Big)
&=& \lim_{n\to\8}\Big[ n \chi_v (-\dt m_0 - \xi_2(c_n))\big( k(a_n,v) -1-i\is v{\dt m_0}-i\is v{\xi_2(c_n)}\big)\\
&&+ n \chi_v (-\dt m_0 - \xi_2(c_n)) \big( 1+i\is v{\dt m_0}+i\is v{\xi_2(c_n)} \big) -n\Big]\\
&=&
\lim_{n\to\8}n\big( k\big(\dt,v\big)-1-i\langle \dt v,m_0 \rangle \big)
+ \lim_{n\to\8} n\big(k(c_n,v)-1-i\is v{\xi_2(c_n)}\big) \\
&& + \lim_{n\to\8}\big[n \chi_v (-\dt m_0 - \xi_2(c_n)) \big( 1+i\is v{\dt m_0}+i\is v{\xi_2(c_n)} \big) -n\big]\\
&=& C_-(v) +  C_+(v) +\frac 12\is v{m_0}^2.
\end{eqnarray*}
To prove nondegeneracy of the limit we use exactly the same argument as above.
\end{proof}

\section{Local Limit Theorem}
\label{section-local}
\begin{proof}[Proof of Theorem \ref{mthmloc}]
In this section we will study the Euclidean case
and assume $\a\notin\{ 1,2 \}$. Take $x=0$. In view of Theorem 10.7 \cite{Br}
it is enough to prove that
$$ \lim_{n\to \8} n^{\chi} \E \big[ h(S_n-d_n)\big] = p_\a(0)\int_V h(v)dv,
$$
for every function $h\in L^1$ such that the Fourier transform of $h$ is compactly
supported.
By Propositions \eqref{kel}, \eqref{lis22} and Lemma \ref{cont}, using the Fourier
inversion formula
$$
\E \big[ h(S_n^x-d_n)\big] = \frac 1{(2\pi)^d}\int_V \E\big[ e^{i\is v{S_n-d_n}} \big] \widehat h(v)dv
= \frac 1{(2\pi)^d}\int_V  e^{-i\is v{d_n}} P_{v}^n(1)(0) \widehat h(v)dv.
$$
Take $N=[-\d,\d]^d$ and denote by $J$ the support of $\widehat h$. By Lemma \ref{3.13}, $r(P_{v})<1$, if $v\not=0$ hence
using Lemma \ref{cont} with $f=1$
there
exists $\b>0$ such that for $v\in J\setminus N$: $r(P_{v})<1-\b$.
Therefore
$$
\lim_{n\to\8} n^{\chi}\bigg| \int_{J\setminus N} e^{-i\is v{d_n}}(P_{v}^n 1)(0)\widehat h(v)dv
\bigg|\le \lim_{n\to \8} C n^{\chi} (1-\b)^n = 0.
$$
Hence we have reduced the problem to computing the limit
$$
\lim_{n\to\8}\frac{ n^{\chi}}{(2\pi)^d}\int_N e^{-i\is v{d_n}} (P_{v}^n1)(0)\widehat h(v)dv
= \lim_{n\to\8}\frac{ n^{\chi}}{(2\pi)^d}\int_N e^{-i\is v{d_n}} \big[
k^n(v)\pi_v(1)(0) + Q_{v}^n(1)(0)
\big]\widehat h(v)dv
$$
Next by Proposition \ref{kel} there exists $\g>0$ such that $\|Q_{v}\|\le 1-\g$
for $v\in N$ hence
$$
 \lim_{n\to\8}{ n^{\chi}}\bigg|\int_N e^{-i\is v{d_n}} (Q_{v}^n1)(0) \widehat h(v)dv\bigg|
\le \lim_{n\to\8} C n^{\chi}(1-\g)^n = 0.
$$
To handle with the remaining term we take a similarity $c_n$ such that $|c_n|=n^{-\frac{\chi}d}$ change variables $v\mapsto  c_nv$ and obtain
\begin{multline}
\label{lebf}
 \lim_{n\to\8}{n^{\chi}}\int_N e^{-i\is v{d_n}}
k^n(v)\pi_{v}(1)(0) \widehat h(v)dv \\
= \lim_{n\to\8} \int_{\{|v|<\d n^{\chi}\}} \Big( e^{-i\is {c_n v}m}
k\big(c_n v\big)\Big)^n \pi_{ c_n v}(1)(0) \widehat h\big( c_n v \big)dv,
\end{multline}
where $m=0$ if $\a<1$. Now we are going to use the Lebesgue theorem.
For this purpose we need that for every $v\in V$ there are $\d>0$ and $D>0$ such that
\begin{equation}
\label{zamiana}
\big|e^{i\is v{m}}k(v)
\big| \le e^{- D |v|^{\frac d{\chi}}}.
\end{equation}
Indeed assume first $\a>2$.
Then by Theorem \ref{asymptkt}, for small values of $|v|$ we have
\begin{eqnarray*}
e^{-i\is vm}k(v) &=&
\Big( 1-i\is vm - \frac{1}2\big( \is vm^2 +o(|v|) \big)\Big) \cdot \Big(
1 + i\is vm +\big( C_{2+}(v)+o(|v|) \big)\Big)\\
&=& 1+  C_{2+}(v) -\frac 12 \is vu^2 + o(|v|)
\end{eqnarray*}
Moreover nondegeneracy of the limit in Theorem \ref{mthm} implies
$C_{2+}(v)-\frac 12 \is vm^2<0$, that gives \eqref{zamiana} in this case.

\medskip

If $\a<2$, then by Theorem \ref{asymptkt}
$$
k(v) = 1+i\is v{m} + |v|^\a \big(  C_\a(\ov v) + o(|v|) \big),
$$
with $\Re C_\a(\ov v)<0$. Therefore
\begin{eqnarray*}
\Big| e^{-i\is v{m}}  k(v)\Big|^2 &=& \Big| \big( 1-i\is v{m} + o(|v|) \big)
\big( 1+i\is v{m} +|v|^\a (C_\a(\ov v)+o(|v|))  \big)\Big|\\
&=& \Big| 1+ |v|^\a \big(C_\a(\ov v)+o(|v|)\big) \Big|^2= 1 + |v|^\a\big( 2\Re C_\a(\ov v) + o(|v|) \big) + O(|v|^{2\a})\\
&\le& e^{-D |v|^\a},
\end{eqnarray*}
that proves \eqref{zamiana}.
Therefore we may use the Lebesgue dominated theorem and pass in \eqref{lebf} to the limit under the integral. Then
reasoning as in the proof of Theorem \ref{mthm} we obtain that the limit above is equal to
$$ \widehat h(0) \cdot \int_V  \Phi_{\a}(v) dv =
 (2\pi)^d p_\a(0)\int_V h(v)dv, \qquad \mbox{if }\a\in(0,1)\cup (1,2)
$$
and
$$ \widehat h(0) \cdot \int_V \Phi_{2+}(v) dv=
(2\pi)^d p_\a(0) \int_V h(v)dv, \qquad \mbox{if }\a>2,
$$ that proves the Theorem.
\end{proof}
Theorem \ref{mthmloc} can be interpreted as a local limit theorem for a random
walk on a homogeneous space. This interpretation brings up some aspects already
encountered for the case of groups in \cite{V}, namely that the degree of the monomial part of the asymptotics
depends of the measure, hence is not determined by the geometry.

One is led to consider the Markov chain on $\wt V = V\times \R^d$ defined by
$x_n = X_n^x$, $y_n = y+S_{n-1}^x =  y_{n-1}+x_{n-1}$. We denote by $\wt P$ its transition
kernel. Then $\wt P$ is a fibered Markov kernel over $P$ (see \cite{GH}), with
typical fiber $\R^d$. Clearly the 'vertical translations' $(x,y)\mapsto(x,y+r)=(x,y)\circ r$
with $r\in \R^d$ commutes with $\wt P$. We denote also by $\circ$ the convolution operation between measures on
$\wt V$ and on $\R^d$, and by $\l$ a Lebesgue measure on $\R^d$. Since
$\wt V=V\times \R^d$, we can identify the measure $\nu$ on $V$ with a measure $\wt\nu$ on $\wt V$.
Then we observe that $\wt\nu\circ\l$ is a $\wt P$ stationary measure. Furthermore $h=(g,b)\in H$ acts
on $V\times \R^d$ by $h(x,y)=(gx+b,y+x)$. This is an affine action of $H$ which is part
of a natural action of a larger group $\wt H$ on $\wt V$ considered as a homogeneous space as follows.

Let $T$ be the real Lie group of $2d\times 2d$ matrices of the form
$
\xi = \left( \begin{array}{cc}g\ & 0\\ u\ & I\end{array} \right),
$  where $g\in G$, $u\in {\rm End}(V)$.  Then $T$ acts on $\R^{2d}$ and
we consider the corresponding semidirect product $\wt H= T\ltimes \R^{2d}$.
Then $\wt V = \R^{2d}$ is a homogeneous space of $\wt H$, i.e. $\wt V = \wt H /\, T$.
The action  of $\wt h=(\xi,\eta)$ ($\xi\in T$, $\eta=(b,c)\in \R^{2d}$) on $\wt v=(x,y)\in \wt V$
is given by $x'=gx+b$, $y'=y+ux+c$. Hence this action commutes with the 'vertical translations'
 on $\wt V$. We recover the $H$-action on  $V$ as a factor of the $\wt H$-action by the vertical translations.

 In particular we denote by $\wt\mu$ the push forward of $\mu$ by the map $(g,b)\mapsto(\xi,\eta)$
 with
$
\xi = \left( \begin{array}{cc}g\ & 0\\ I\ & I\end{array} \right)
$,
$
\eta = \left( \begin{array}{c}b\\ 0\end{array} \right)
$. Then we can write $\wt P(\wt v,\cdot)=\wt\mu*\d_{\wt v}$, hence
$\wt P^n(\wt v,\cdot)=\wt\mu^n*\d_{\wt v}$. We know that $X_n^x$ converges in law to $\nu$.
Furthermore the theorem tells us that if $\mu_n$ denotes the law of $y_n=y+ S_n^x$, then the sequence of measures
$n^{\chi}(\mu_n\circ \d_{-d_n})$ converges weakly to $p_{\a}(0)\l$. Then,  following the analysis of \cite{GH}
for local limit asymptotics in the context of fibered Markov kernels we get
\begin{cor}
With the above notations, for any $\wt v\in\wt V$ we have the weak convergence:
$$
\lim_{n\to\8} n^{\chi}(\wt\mu^n*\d_{\wt v})\circ \d_{-d_n} = p_\a(0)\wt \nu\circ \l.
$$
\end{cor}
\appendix
\section{On the structure of closed subgroups of $G=D\times K$}
\label{app1}
Here $V=\oplus_{j=1}^l V_{\l_j}$,
$D$ is the one parameter subgroup of elements $\g_a$, which act
on $V_{\l_j}$ by multiplication by $a^{\l_j}$, $\tau(x)=\sum_{j=1}^l|x_j|^{\frac 1{\l_j}}$
 and $K=\{g\in {\rm GL}(V);\; |g|=1\}$.
We denote $\wt G = K\times(\R_+^*)^l$ and we observe that $\wt G$ is a linear
algebraic group, which contains $G$ as a closed subgroup. We denote by
$\wt\gg$ the Lie algebra of $\wt G$, by $\gg$ and $\kk$ the Lie algebras
of $G$ and $K$, respectively.
\begin{prop}
Let $G_1$ be a closed noncompact subgroup of $G$, $R_1$ its
projection on $D$, $Z_1$ its center, $K_1=G_1\cap K$. Then $R_1 = D$
or $R_1 = \iss{p}$, $p\in D$, $|p|>1$. There exists $Y_1\in \wt
\gg$, $|\exp Y_1|>1$,
such that ${\rm Ad} G_1(Y_1)=Y_1$, $\exp  {Y_1}\in G_1$. 
In particular
the subgroup $\exp\R Y_1$ commutes with $G_1$ and $\exp \Z Y_1$ is a central subgroup of $G_1$.
Moreover
\begin{itemize}
\item if $R_1=D$, then $G_1 = A_1\times K_1$, with $A_1=\exp\R Y_1$ and $Z_1=(Z_1\cap K)\times A_1$.
\item if $R_1=\iss{p}$, then $G_1$ contains $A_1\times K_1$
as a normal subgroup of finite index, with $A_1=\exp\Z Y_1$
and $G_1$ is the semidirect product of $\langle  g_1\rangle$ and $K_1$, with $g_1\in G_1$ and $g_1^r = \exp Y_1$ for
some $r\in\N$. The center $Z_1$ of $G_1$  is the product of
$Z_1\cap K$ by a cyclic subgroup  $\iss{z_1}$ $(|z_1|>1)$ such that $\iss{z_1}/ \iss{z_1}\cap A_1$ is finite.

\end{itemize}
\end{prop}
\begin{proof}
We observe that the projection map $\pi$ of $G$ on $D$ has compact kernel
hence $\pi$ is proper. It follows that $R_1=\pi(G_1)$ is closed, hence
$\pi(G_1)$ is either $\{1\}$ or $D$ or $\iss{p}$ ($c>1$). Since $G_1$
is non compact $\pi(G_1)=\{1\}$ is excluded, hence the first assertion.
On the other hand $K_1$ is normal in $G_1$. Every element $X$ of $\gg$
can be written as $X=(\l ,\t)$ with $\l\in\R$ and $\t \in \kk$,
an antisymmetric matrix. The quadratic form $q$ on $\gg$ defined by $q(X)=\l^2-{\rm Tr}\t^2$
is positive definite and ${\rm Ad}G$-invariant.

\medskip

If $R_1=D$, the result
is proved as follows. Indeed, the Lie algebra $\gg_1$ is ${\rm Ad}G_1$-invariant
and  $\gg_1$ contains an element $Y_1$ with $\pi(Y_1)\not= 0$. It follows
$G_1=(\exp \R Y_1)\ltimes K_1$, $\gg_1 = \R Y_1 \ltimes \kk_1$.
Since $\kk_1$ and $\gg_1$ are ${\rm Ad}G_1$-invariant, the same is true for the orthogonal
line $\kk_1^{\perp}\subset \gg_1$. On the other hand, we have for any $t\in\R$, $g\in G_1$:
$ge^{t Y_1}g^{-1}e^{-tY_1}\in K_1$, hence ${\rm Ad} g (Y_1)\in Y_1+\kk_1$. Then
the affine hyperplane $Y_1+\kk_1$ of $\gg_1$ is ${\rm Ad}G_1$-invariant, hence
the point of intersection  $Y$ of $\kk_1^{\perp}$ and $Y_1+\kk_1$ satisfies
${\rm Ad} G_1(Y)=Y$. It follows $G_1=A_1\times K_1$, $Z_1=A_1\times (Z_1\cap K)$.

\medskip

 If $R_1=\iss{p}$, we
consider $y\in G_1$ with $\pi(y)=p$, and the Zariski closure $L$
of the subgroup $\iss{y}$. Then $L$ is a closed abelian Lie group with a finite
number of connected components. Let $L^0$ be the connected component
of $e$ in $L$. Let $r\in \N$ with $y^r=\exp Y_0\in L^0$
with $Y_0\in \wt \gg$, hence $\exp \R Y_0 \subset L^0$, ${\rm Ad}y(Y_0)=Y_0$. Since for any $g\in G_1$,
$n\in\Z$, $g y^n g^{-1} y^{-n}\in K_1$ and $K_1$ is algebraic, we have
also, $gzg^{-1}z^{-1}\in K_1$ for any $z\in L$. In particular for any $t\in\R$,
$g e^{t Y_0} g^{-1} e^{-t Y_0}\in K_1$. It follows ${\rm Ad}g(Y_0)-Y_0\in\kk_1$,
hence the affine hyperplane $Y_0 +\kk_1$ of $\R Y_0 + \kk_1$ is ${\rm Ad}G_1$-invariant.
Since $\kk_1$ and $\R Y_0 + \kk_1$ and are ${\rm Ad}G_1$-invariant, we can repeat
the argument used if $R_1=D$: the point $Y_1=Y_0+U$, $U\in \kk_1$ of intersection
of $\kk_1^{\perp}$ and $Y_0 + \kk_1$ satisfies ${\rm Ad}G_1(Y_1)=Y_1$.
In particular $[Y_1,\kk_1]=[Y_1,U]=\{0\}$, $\exp Y_1 = \exp Y_0\exp U$. Then, using
$y^r = \exp Y_0\in G_1$, $\exp U\in K_1$, we get $\exp Y_1\in G_1$, hence $\exp \Z Y_1$ is
a central subgroup of $G_1$.
Since $G_1=\iss{y}\rtimes K_1$ and $y^r = \exp Y_0$,
we conclude that $A_1\times K_1$ has finite index in $G_1$.
Furthermore ${\rm Ad}{y} (Y_1)=Y_1$ and ${\rm Ad}{y}(Y_0)=Y_0$ imply
${\rm Ad}{y}(U)=U$, hence $(y \exp \frac Ur)^r = y^r\exp U = \exp Y_1$.
Since $|y\exp \frac Ur|=p$ and $g_1= y \exp \frac Ur\in G_1$ we conclude $G_1 =\langle \exp {\iss{g_1}} \rangle\ltimes K_1$,
with $g^r_1=\exp Y_1$.
Using the above we can write $Z_1=\iss{z}\times (Z_1\cap K_1)$. Also $A_1\times (Z_1\cap K_1)$
is a subgroup of finite index in $Z_1$. This follows from the fact that $\pi$
defines an isomorphism of $\iss{z}$ onto a cyclic subgroup of $D$, which contains $\pi(A_1)$
as a finite index subgroup. Then, for some $n\in N$, $z^n=u a_1$ with $a_1\in A_1$,
$u\in (Z_1\cap K)^0$. We can write $u^{-1}=v^n$ with $v\in (Z_1\cap K)^0$,
hence $(zv)^n = z^n u^{-1}=a_1$. Then with $z_1 =zv$ we have $Z_1=\iss{z_1}\times(Z_1\cap K_1)$,
$\iss{z_1}\cap A_1\supset \iss{a_1}$, hence $\iss{z_1}/ \iss{z_1}\cap A_1$ is finite.
\end{proof}
\bibliographystyle{alpha}
\newcommand{\etalchar}[1]{$^{#1}$}

\end{document}